\documentclass[english]{article}
\usepackage[utf8]{inputenc}
\usepackage{babel}
\usepackage[T1]{fontenc}

\usepackage[
backend=biber,
style=alphabetic,
]{biblatex}
\usepackage{csquotes}
\usepackage[all]{xy}
\usepackage{verbatim}
\usepackage{units}
\usepackage{textcomp}
\usepackage[margin=0.75in]{geometry}

\usepackage{amsmath}
\usepackage{amsthm}
\newtheorem{thm}{Theorem}[section]
\newtheorem*{thm*}{Theorem}
\newtheorem{lem}[thm]{Lemma}
\newtheorem{cor}[thm]{Corollary}
\newtheorem{prop}[thm]{Proposition}

\theoremstyle{definition}
\newtheorem{defn}[thm]{Definition}
\newtheorem*{ex*}{Example}

\theoremstyle{remark}
\newtheorem{rmk}[thm]{Remark}

\usepackage{amssymb}
\usepackage{stmaryrd}
\usepackage{setspace}
\usepackage{wasysym}
\usepackage{enumitem}
\usepackage{hyperref}
\hypersetup{
    colorlinks,
    citecolor=black,
    filecolor=black,
    linkcolor=black,
    urlcolor=black
}
\usepackage[nottoc,numbib]{tocbibind}

\title{Comparison of integral structures on the space of modular forms of full level $N$}
\author{Anthony Kling}
\date{}

\begin{document}
\maketitle

\begin{abstract}
Let $N\geq3$ and $r\geq1$ be integers and $p\geq2$ be a prime such that $p\nmid N$. One can consider two different integral structures on the space of modular forms over $\mathbb{Q}$, one coming from arithmetic via $q$-expansions, the other coming from geometry via integral models of modular curves. Both structures are stable under the Hecke operators; furthermore, their quotient is finite torsion. Our goal is to investigate the exponent of the annihilator of the quotient. We will apply methods due to Brian Conrad to the situation of modular forms of even weight and level $\Gamma(Np^{r})$ over $\mathbb{Q}_{p}(\zeta_{Np^{r}})$ to obtain an upper bound for the exponent. We also use Klein forms to construct explicit modular forms of level $p^{r}$ whenever $p^{r}>3$, allowing us to compute a lower bound which agrees with the upper bound. Hence we are able to compute the exponent precisely.
\end{abstract}

\tableofcontents{}

\section{Introduction}

Let $k\geq1$ be an integer and let $\Gamma$ be a congruence subgroup of $\rm{SL}_{2}(\mathbb{Z})$ of level $N$ \emph{i.e. }a subgroup of ${\rm SL}_{2}(\mathbb{Z})$ containing the kernel, $\Gamma(N)$, of the usual reduction mod $N$ map 

\[{\rm SL}_{2}(\mathbb{Z})\rightarrow{\rm{SL}_{2}(\mathbb{Z}/N\mathbb{Z})}.\] Consider the $\mathbb{Q}$-vector space $M_{k}(\Gamma,\mathbb{Q})$ consisting of modular forms of weight $k$ and level $\Gamma$ over $\mathbb{Q}$. We associate to each modular form $f$ in $M_{k}(\Gamma,\mathbb{Q})$ and each cusp $c$ of $\Gamma$, the $q$-expansion of $f$ at $c$, denoted $f_{c}$, which is a power series in $\mathbb{Q}[[q]]$.

We concern ourselves with $\mathbb{Z}$-structures on the vector space $M_{k}(\Gamma,\mathbb{Q})$ \emph{i.e. }$\mathbb{Z}$-submodules $M'$ of $M_{k}(\Gamma,\mathbb{Q})$ such that the natural map $M'\otimes_{\mathbb{Z}}\mathbb{Q}\rightarrow M_{k}(\Gamma,\mathbb{Q})$ is an isomorphism. The space $M_{k}(\Gamma,\mathbb{Q})$ naturally has two different $\mathbb{Z}$-structures. Define the first $\mathbb{Z}$-structure to be \[M_{k}(\Gamma,\mathbb{Z})=\left\{ f\in M_{k}(\Gamma,\mathbb{Q}):f_{\infty}\in\mathbb{Z}[[q]]\right\} \] which consists of modular forms in $M_{k}(\Gamma,\mathbb{Q})$ whose $q$-expansion at the cusp $\infty$ has integral coefficients. One can show the Hecke operators of $M_{k}(\Gamma,\mathbb{Q})$ preserve integrality at the cusp $\infty$ by explicitly computing the $q$-expansion at $\infty$ under the Hecke operators (see \cite[§4.9.2]{[Kato04]}). 

The second $\mathbb{Z}$-structure we consider is $M_{k,\mathbb{Z}}$, which consists of modular forms in $M_{k}(\Gamma,\mathbb{Q})$ whose $q$-expansions at \emph{all} cusps have integral coefficients. This structure is also stable under the Hecke operators (cf. \cite[Theorem 1.2.2]{[Con07]}). There is an obvious containment $M_{k}(\Gamma,\mathbb{Z})\subseteq M_{k,\mathbb{Z}}$ with quotient that is a torsion $\mathbb{Z}$-module. Our aim is to study and determine the annihilator. 

To better understand and work with $M_{k,\mathbb{Z}}$, we realize $M_{k,\mathbb{Z}}$ as the global sections of some line bundle $\underline{\omega}_{\mathfrak{X}}^{\otimes k}$ on $\mathfrak{X}$, the moduli space parameterizing $\Gamma$-level structures over $\mathbb{Z}$. More precisely, for an integer $N\geq1$,\cite[§3]{[KM85]} considers four different moduli problem parameterizing $\Gamma(N)$-, $\Gamma_{1}(N)$-, balanced $\Gamma_{1}(N)$-, and $\Gamma_{0}(N)$-structures on elliptic curves. The definition of $\Gamma(N)$-structures is given in Appendix \ref{Ap - mod curve} since this will be our focus. When $N\geq3$, the moduli problems parameterizing $\Gamma(N)$-, $\Gamma_{1}(N)$-, and balanced $\Gamma_{1}(N)$-structures are represented by a regular, flat two-dimensional scheme $\mathfrak{Y}(\Gamma)$ over $\mathbb{Z}$, by \cite[{}5.5.1]{[KM85]}. For the rest of the introduction, we let $\Gamma$ denote one of these level structures. The scheme $\mathfrak{Y}(\Gamma)$ extends to an arithmetic surface $\mathfrak{X}(\Gamma)$ over $\mathbb{Z}$, known as the modular curve, so that $M_{k,\mathbb{Z}}$ is identified with $H^{0}(\mathfrak{X}(\Gamma),\underline{\omega}_{\mathfrak{X}(\Gamma)}^{\otimes k})$. From this, we can show that $H^{0}(\mathfrak{X}(\Gamma),\underline{\omega}_{\mathfrak{X}(\Gamma)}^{\otimes k})$ and $M_{k}(\Gamma,\mathbb{Z})$ are both finitely-generated $\mathbb{Z}$-modules of the same rank. Thus the quotient $M_{k}(\Gamma,\mathbb{Z})/H^{0}(\mathfrak{X}(\Gamma),\underline{\omega}_{\mathfrak{X}(\Gamma)}^{\otimes k})$ is torsion.

Let $p$ be a prime. If $p$ does not divide the level $N$, then the $p$-adic valuation of the annihilator is zero (see Remark \ref{Trivial Ann}). Thus we focus on primes $p$ which divide the level and instead can work over $\mathbb{Z}_{p}$. Let $c$ be a cusp of $\mathfrak{X}(\Gamma)$. Define $\nu_{p}(f_{c})$ to be the minimal $p$-adic valuation among all the coefficients of $f_{c}$. Then our general goal can be restated as follows: we seek to compute the smallest integer $e\geq0$ such that $\nu_{p}(p^{e}f_{c})\geq0$ for all cusps $c$ of $\mathfrak{X}(\Gamma)$ not equal to $\infty$ and all $f\in M_{k}(\Gamma,\mathbb{Z}_{p})$. 

The problem of computing and bounding $e$ has a long history, which we briefly recall.  Computing a bound for $e$ for arbitrary weight $k$ and level $\Gamma_{0}(p)$ was done in \cite[§3.19, §3.20]{[DeRa]}, where they obtained an upper bound of \[e\leq\frac{kp}{p-1}.\] Their methods involved using intersection theory on $\mathfrak{X}_{0}(p)$ and the fact that $\mathfrak{X}_{0}(p)_{/\mathbb{F}_{p}}$ is reduced with two irreducible components. 

In \cite{[Edi06]}, Edixhoven investigates the situation of weight
$k=2$ and level $\Gamma_{0}(N)$ \emph{cusp forms }where ${\rm ord}_{p}(N)=1$. In particular, he establishes the existence of a non-zero global section of $\Omega|_{C_{0}'}(-e\cdot S)$ where $\Omega$ is the relative dualizing sheaf of $\mathfrak{X}_{0}(N)$, $C_{0}'$ is essentially\footnote{More precisely, since $\mathfrak{X}_{0}(N)$ is not necessarily regular, \cite{[Edi06]} works instead with some finite cover $X'\rightarrow\mathfrak{X}_{0}(N)$ by appending extra level structure to $\Gamma_{0}(N)$. Let $C_{0}$ denote the irreducible component of $\mathfrak{X}_{0}(N)_{/\mathbb{F}_{p}}$ containing the cusp $0$. Then $C_{0}'$ is the inverse image of $C_{0}$ under $X'\rightarrow\mathfrak{X}_{0}(N)$.} the irreducible component containing the cusp $0$, and $S$ is the divisor given by the sum of the supersingular points. By computing $\deg(\Omega|_{C_{0}'}(-e\cdot S))$ in terms of $e$ and using the inequality \[0<\deg(\Omega|_{C_{0}'}(-e\cdot S)),\] Edixhoven is able to conclude $e<1+\frac{2}{p-1}$ which forces $e=1$ whenever $p\geq3$. This agrees with the bound provided by \cite{[DeRa]}. In the case $p=2$, Edixhoven separately concludes $e=2$. 

In \cite[Appendix B]{[Con17]}, Conrad investigates $e$ in the general situation ${\rm ord}_{p}(N)=r\geq1$ and arbitrary weight $k$. He begins by developing intersection theory on the regular proper Artin stack $\mathfrak{X}_{0}(N)$ over $\mathbb{Z}$. It is worth noting the stack $\mathfrak{X}_0(N)$ in \cite{[Con07]} which Conrad considers does not agree with the one in \cite{[DeRa]}.  In particular, $\mathfrak{X}_0(N)$ as defined in \cite{[DeRa]} is a Deligne-Mumford stack instead of merely Artin (see \cite{[Ces17]} for more on this issue). We will not need to concern ourselves with this issue since we only work with the arithmetic surface $\mathfrak{X}(N)$ which parameterizes $\Gamma(N)$-structures.  Furthermore, Conrad's expression for the exponent readily holds for $\mathfrak{X}(N)$ as we show in Section \ref{sec:The-Exponent}. Let $T$ denote the $r\times r$ matrix obtained by removing the column and row corresponding to the irreducible component containing the cusp $\infty$ from the intersection matrix of $\mathfrak{X}_{0}(N)$. Conrad provides an expression for the upper bound in terms of $k,p,r,$ and the entries of $T^{-1}$. As pointed out in \cite[footnote 5]{[CNS23]}, this bound is incorrect for $r>1$ due to a typo in the values of the multiplicities of the components of $\mathfrak{X}_0(N)_{\mathbb{F}_p}$ used in the calculation of $T^{-1}$. In our situation, since $\mathfrak{X}(N)$ is reduced, these multiplicities are all equal to $1$. For general $r$, it is not clear how to obtain a uniform description of the entries of $T^{-1}$, and this prevents Conrad from establishing explicit bounds for general $r$.  

Lastly, \cite{[CNS23]} uses an automorphic approach to bound the exponent in the situation $\mathfrak{X}_0(N)$.  This approach appears quite powerful; in particular \cite[Theorem 4.6]{[CNS23]} provides bounds depending on the cusp at which the q-expansion is being taken.  These bounds are also shown to be sharp in a few cases via explicit computations (see \cite[Example 4.8]{[CNS23]}).  We will still follow the approach in \cite[Appendix B]{[Con17]} as this yields an upper bound which we show is sharp in all cases except two: when the level is exactly $2$ or $3$.  Consequently, this provides an exact computation for the exponent.

Let $N\geq3$ and $r\geq1$ be an integers and $p\geq2$ be a prime such that $p\nmid N$. Fix a primitive root of unity $\zeta_{Np^{r}}$ and let $\pi$ be a uniformizer of $\mathbb{Z}_{p}[\zeta_{Np^{r}}]$. We consider the situation of the modular curve $\mathfrak{X}(Np^{r})$ over $\mathbb{Z}_{p}[\zeta_{Np^{r}}]$, which parameterizes $\Gamma(Np^{r})^{{\rm can}}$-structures (see the paragraph preceding Theorem \ref{KM5.1.1}, or \cite[§9]{[KM85]}). In this situation, the special fiber of $\mathfrak{X}(Np^{r})$ consists of $p^{r}+p^{r-1}$ many irreducible components. We also restrict ourselves to modular forms of even weight. Applying Conrad's method in this situation, we obtain a formula in terms of the entries of $T^{-1}$ and $\deg(\underline{\omega}^{\otimes2}|_{\Lambda})$, where $\Lambda$ is any irreducible component of the special fiber of $\mathfrak{X}(Np^{r})$. We make the following two significant computations:

\begin{itemize}
\item We explicitly compute the entries of $T^{-1}$ for all $r$ and $p$.
\item We explicitly compute $\deg(\underline{\omega}^{\otimes2}|_{\Lambda})$. Note that Conrad in \cite[Appendix B]{[Con17]} bounds $\deg(\underline{\omega}^{\otimes2}|_{\Lambda})$ by $\deg(\underline{\omega}^{\otimes2})$. In our situation,

\begin{align*}
\deg(\underline{\omega}_{\mathfrak{X}(Np^{r})}^{\otimes2}) & = \frac{2}{24}{\rm SL}_{2}(\mathbb{Z}/N\mathbb{Z}){\rm SL}_{2}(\mathbb{Z}/p^{r}\mathbb{Z})\\
 & = \frac{{\rm SL}_{2}(\mathbb{Z}/N\mathbb{Z})}{12}(p^{3r}-p^{3r-2})
\end{align*}

by \cite[{}10.13.12]{[KM85]}. In Theorem \ref{Deg of Res Modular Sheaf }, we show 

\[
\deg(\underline{\omega}_{\mathfrak{X}(Np^{r})}^{\otimes2}|_{\Lambda})=\frac{{\rm SL}_{2}(\mathbb{Z}/N\mathbb{Z})}{12}(p-1)p^{2r-1}
\]

independent of $\Lambda$, which is smaller than $\deg(\underline{\omega}_{\mathfrak{X}(Np^{r})}^{\otimes2})$ by a factor of $p^{r-1}(p+1)$. This savings leads to a marked improvement in our upper bound.
\end{itemize}

In Theorem \ref{Entries of Inv Int Matrix}, we calculate the entries of $T^{-1}$ for all $r$. First we describe the entries of $T$ by explicitly computing the intersection number between each irreducible component of the special fiber, using \cite[{}13.8.5]{[KM85]}, then we compute the self-intersection numbers.

We compute 
\[
T=\deg{\rm S}(N)\left(\begin{array}{cc} M(p^{r})_{\hat{1},\hat{1}} & \mathbf{1}_{p^{r}-1\times p^{r-1}}\\ 
\mathbf{1}_{p^{r-1}\times p^{r}-1} & p^{2}M(p^{r-1}) \end{array}\right)
\]
where $M(p^{r})$ is a $p^{r}\times p^{r}$ circulant matrix, dependent on $p$ and $r$, $M(p^{r})_{\hat{1},\hat{1}}$ is the matrix $M(p^{r})$ with the first row and column removed, $\mathbf{1}_{n\times m}$ is the $n\times m$ matrix consisting of all 1's, and $\deg{\rm S}(N)$ is the degree of supersingular locus in $\mathfrak{X}(Np^{r})$. Applying the Woodbury Matrix Identity (see Appendix \ref{apx: Inv via WMI}) we obtain a formula for the entries of $T^{-1}$ involving the entries of $(M(p^{r})_{\hat{1},\hat{1}})^{-1}$ and $p^{-2}M(p^{r-1})^{-1}$, as well as their row sums and the total sum of all entries. 

Since $M(p^{r})$ is circulant, we have a description of its eigenvalues and corresponding eigenvectors, allowing us to diagonalize $M(p^{r})$. In turn, this allows us to explicitly compute the entries of $M(p^{r})^{-1}$. Consequently, using Proposition \ref{JCP16}, we also obtain the entries of the inverse of $M(p^{r})_{\hat{1},\hat{1}}$. After more careful calculations, we obtain an exact expression for each entry of $T^{-1}$. 

In Theorem \ref{Deg of Res Modular Sheaf }, we compute $\deg(\underline{\omega}_{\mathfrak{X}(Np^{r})}^{\otimes2}|_{\Lambda})$. Using the Kodaira-Spencer isomorphism (Theorem \ref{K-S Iso}) and the adjunction formula (Theorem \ref{Adj Formula}), we are able to identify $\underline{\omega}_{\mathfrak{X}(Np^{r})}^{\otimes2}$ with the relative dualizing sheaf \[\Omega_{\mathfrak{X}(Np^{r})/\mathbb{Z}_{p}[\zeta_{Np^{r}}]}(-\mathfrak{C}(Np^{r}))\] twisted by minus the cuspidal divisor $\mathfrak{C}(Np^{r})$, tensored with relative dualizing sheaf $\Omega_{\mathfrak{X}(Np^{r})/\mathfrak{X}(N)}$ of the map ${\rm pr}:\mathfrak{X}(Np^{r})\rightarrow\mathfrak{X}(N)$ induced by forgetting the $\Gamma(p^{r})$-level structure. 

Investigating $\Omega_{\mathfrak{X}(Np^{r})/\mathfrak{X}(N)}$ amounts to understanding the different of the morphism ${\rm pr}$, which amounts to understanding the different of 
\[
{\rm pr}_{x}:{\cal O}_{\mathfrak{X}(N),{\rm pr}(x)}\rightarrow{\cal O}_{\mathfrak{X}(Np^{r}),x}
\] where $x\in\mathfrak{X}(Np^{r})$ is a codimension 1 point. Let $d_{x}$ denote the valuation of the different of ${\rm pr}_{x}$ in ${\cal O}_{\mathfrak{X}(Np^{r}),x}$. We split our analysis of ${\rm pr}_{x}$ into two cases. If $x$ is a closed point of the generic fiber, then we can compute $d_{x}$ over $\mathbb{C}$. This amounts to understanding the ramification of the analogous map of Riemann surfaces $X(Np^{r})\rightarrow X(N)$, which is done in \cite[Prop. 1.37]{[Shi71]}. If $x$ is a generic point of the special fiber, then we show all $d_{x}$ contributions are the same which taken together contribute nothing to the different of ${\rm pr}$.  We ultimately conclude 
\[
\Omega_{\mathfrak{X}(Np^{r})/\mathfrak{X}(N)}\simeq{\cal O}_{\mathfrak{X}(Np^{r})}((p^{r}-1)\mathfrak{C}(Np^{r}))
\]
and consequently 
\[
\underline{\omega}_{\mathfrak{X}(Np^{r})}^{\otimes2}\simeq\Omega_{\mathfrak{X}(Np^{r})/R}(\mathfrak{C}(Np^{r})).
\]
This identification allows us to directly compute $\deg(\underline{\omega}_{\mathfrak{X}(Np^{r})}^{\otimes2}|_{\Lambda})$ by computing the quantities 
\[
\deg(\mathfrak{C}(Np^{r})|_{\Lambda})\mbox{ and }\deg(\Omega_{\mathfrak{X}(Np^{r})/R}|_{\Lambda}).
\] The relative dualizing sheaf $\Omega_{\mathfrak{X}(Np^{r})/R}$ enjoys strong functoriality properties, enabling us to identify 
\[
\Omega_{\mathfrak{X}(Np^{r})/R}|_{\Lambda}\simeq\Omega_{\Lambda/\mathbb{F}_{p^{\varphi(N)}}}\otimes{\cal O}_{\mathfrak{X}(Np^{r})}(-\Lambda)|_{\Lambda}.
\] This allows us to compute $\deg(\Omega_{\mathfrak{X}(Np^{r})/R}|_{\Lambda})$ in terms of the genus of an Igusa curve, and the self-intersection number of $\Lambda$. Combined with our calculation of the entries of $T^{-1}$, we arrive at the following upper bound.

\begin{thm}
Let $k\geq1$, $N\geq3$, and $r\geq1$ be integers and $p\geq2$ be a prime such that $p\nmid N$. The exponent $e$ of $\pi$ in the annihilator of \[M_{2k}(\Gamma(Np^{r}),\mathbb{Z}_{p}[\zeta_{Np^{r}}])/H^{0}(\mathfrak{X}(Np^{r}),\underline{\omega}_{\mathfrak{X}(Np^{r})}^{\otimes2k})\] is bounded above by 
\[
e\leq2kp^{r-1}(pr-r+1).
\]
\end{thm}

For $(r_{1},r_{2})\in\mathbb{Q}^{2}-\mathbb{Z}^{2}$, let $\kappa_{(r_{1},r_{2})}(\tau)$ denote the Klein form associated to $(r_{1},r_{2})$ (see Definition \ref{Klein form}). A lower bound for $e$ is obtained by explicitly constructing a modular form of level $p^{r}$ out of a product of Klein forms. Let $\left\{ m(t)\right\} _{t=1}^{N-1}$ be a family of integers. \cite[Thm. 2.6]{[EKS11]}, which builds upon results in \cite[§2.1, §2.4]{[KL81]}, provides a criterion on $\left\{ m(t)\right\} _{t=1}^{N-1}$ for when a product of Klein forms
\[
\kappa(\tau)=\prod_{t=1}^{N-1}\kappa_{(t/N,0)}(N\tau)^{m(t)}
\] is a nearly holomorphic modular form of level $\Gamma_{1}(N)$ and weight $-\sum_{t=1}^{N-1}m(t)$. We are able to choose a family $\left\{ m(t)\right\} _{t=1}^{p^{r}-1}$ such that $\kappa(\tau)$ is a modular form of level $\Gamma_{1}(p^{r})$ and weight $2$ with integral $q$-expansion at $\infty$. The inclusion $\Gamma(Np^{r})\le\Gamma_{1}(p^{r})$ induces a map $X(Np^{r})\rightarrow X_{1}(p^{r})$ between the corresponding modular curves over $\mathbb{C}$. Pulling back $\kappa(\tau)$ under this map and taking the $k$th power $\kappa(\tau)^{k}$, we can view $\kappa(\tau)^{k}$ as a modular form of level $\Gamma(Np^{r})$ and weight $2k$. By explicitly computing the $q$-expansion of $\kappa(\tau)$ at the cusp $0$, we obtain the following lower bound. 

\begin{thm}
For $k\geq1$, and $p^{r}>3$, $e$ is bounded below by $2kp^{r-1}(pr-r+1)$. Consequently, $e$ is equal to $2kp^{r-1}(pr-r+1)$.
\end{thm}

We are also able to obtain an upper bound in the situation of cusp forms. Although Edixhoven only considers $\mathfrak{X}_{0}(Np)$ with $p\nmid N$, his method can be adapted to yield an upper bound in our situation as well. However, this upper bound is worse than the one we obtain by a factor of $p^{r}/(pr-r+1)$, yet coincides for $r=1$ (see Remark \ref{Compare Edi Bound}). 

We now briefly summarize the contents of this paper.  In Section \ref{sec:The-Exponent}, we formulate our problem of computing $e$ geometrically, and the resulting formula for it in terms of the intersection theory of $\mathfrak{X}(Np^{r})$, and the degree of $\underline{\omega}_{\mathfrak{X}(Np^{r})}^{\otimes2}|_{\Lambda}$. We follow \cite[Appendix B]{[Con17]} to provide an expression that calculates the exponent in the situation of $\mathfrak{X}(Np^{r})$. We also provide material on intersection theory for arithmetic surfaces, most of which comes from \cite[§8, §9]{[Liu02]}.

Section \ref{sec:Intersection-Matrix} is devoted to describing the intersection matrix of $\mathfrak{X}(Np^{r})$, and then computing the entries of $T^{-1}$. In Section \ref{sec:Computing-deg}, we compute $\deg(\underline{\omega}_{\mathfrak{X}(Np^{r})}^{\otimes2}|_{\Lambda})$. Combined with the work in Section \ref{sec:Intersection-Matrix}, this culminates in Theorem \ref{Upper Bound} where we provide an upper bound for the exponent. Lastly, we construct explicit modular forms of level $p^{r}$ in Section \ref{sec:A-Lower-Bound} which provides a lower bound for the exponent that agrees with the upper bound. 

\subsection{\label{sub:Notation}Notation}

If $f:S\rightarrow T$ is a morphism of schemes, for any $T$-scheme $X$ we let $X_{/S}$ denote the base change of $X$ along $f$. When $S={\rm Spec}(A)$ and $T={\rm Spec}(B)$, we will often abuse notation and write $X_{/A}$ for $X_{/S}$. If $X$ is a scheme over a DVR, we denote the special fiber of $X$ by $\bar{X}$. 

Let $N\geq3$ be an integer. We denote the ``compactified'' regular integral model of the modular curve of full level $N$ by $\mathfrak{X}(N)$ over the cyclotomic integers $\mathbb{Z}[\zeta_{N}]$, as presented in \cite{[KM85]}. Refer to Appendix \ref{Ap - mod curve} for more on the modular curve. We denote the cuspidal locus of $\mathfrak{X}(N)$ by $\mathfrak{C}(\mathfrak{X}(N))$  or sometimes $\mathfrak{C}(N)$. We let $\underline{\omega}$ denote the modular sheaf of $\mathfrak{X}(N)$ (see the paragraph proceeding Theorem \ref{X(N) is arith surf}). For an integer $k\geq1$, the global sections of $\underline{\omega}^{\otimes2k}$ define the space of modular forms of weight $2k$ and level $\Gamma(N)$. 

Starting in §\ref{sub: geo exp} and onward, we will usually consider modular forms of level $\Gamma(Np^{r})$ where $p\geq2$ is a prime number such that $p\nmid N$ and $r\geq1$. We will also be working over the DVR $\mathbb{Z}_{p}[\zeta_{Np^{r}}]$ which has uniformizer $\pi=1-\zeta_{p^{r}}$.


\section{\label{sec:The-Exponent}The exponent}

\subsection{Formulation of the exponent}

Let $\mathfrak{X}(N)$ denote the ``compactified'' regular integral model of the modular curve of full level $N\geq3$ and let $\underline{\omega}$ denote the modular sheaf over $\mathfrak{X}(N)$ following the notation in §\ref{sub:Notation}. We will consider two different sub $\mathbb{Z}[\zeta_{N}]$-modules of the space of modular forms $H^{0}(\mathfrak{X}(N)_{/\mathbb{Q}(\zeta_{N})},\underline{\omega}^{\otimes k})$ over $\mathbb{Q}(\zeta_{N})$. The first structure is $H^{0}(\mathfrak{X}(N),\underline{\omega}^{\otimes k})$, the space of modular forms over $\mathbb{Z}[\zeta_{N}]$ (see Definition \ref{def:mod form}), which has the following description in terms of $q$-expansions.

\begin{lem}
\label{Description of H0(N)}  We have 
\[
H^{0}(\mathfrak{X}(N),\underline{\omega}^{\otimes k})=\left\{ f\in H^{0}(\mathfrak{X}(N)_{/\mathbb{Q}(\zeta_{N})},\underline{\omega}^{\otimes k}):f_{c}\in\mathbb{Z}[\zeta_{N}][[q^{1/N}]]\mbox{ for all }c\in\mathfrak{C}(N)\right\} .
\]
\end{lem}
\begin{proof}
This comes from the $q$-expansion principle, as stated in Proposition
\ref{q-exp principle}.
\end{proof}

Let $R$ be a $\mathbb{Z}[\zeta_{N}]$-algebra contained in $\mathbb{Q}(\zeta_{N})$. The second structure is defined as 
\[ M_{k}(N,R):=\left\{ f\in H^{0}(\mathfrak{X}(N)_{/\mathbb{Q}(\zeta_{N})},\underline{\omega}^{\otimes k}):f_{\infty}\in R[[q^{1/N}]]\right\} \]
which are modular forms over $\mathbb{Q}(\zeta_{N})$ whose $q$-expansion at the cusp $\infty$ has coefficients in $R$. 

\begin{lem}\label{M2 base change} 
The usual map \[M_{k}(N,\mathbb{Z}[\zeta_{N}])\otimes_{\mathbb{Z}[\zeta_{N}]}\mathbb{Z}[\zeta_{N},1/N]\rightarrow M_{k}(N,\mathbb{Z}[\zeta_{N},1/N])\] is an isomorphism.
\end{lem}

\begin{proof}
By Proposition \ref{Base change}, we have 
\[H^{0}(\mathfrak{X}(N)_{/\mathbb{Z}[\zeta_{N}]},\underline{\omega}_{\mathfrak{X}(N)_{/\mathbb{Z}[\zeta_{N}]}}^{\otimes k})\otimes_{\mathbb{Z}[\zeta_{N}]}\mathbb{Q}(\zeta_{N})=H^{0}(\mathfrak{X}(N)_{/\mathbb{Q}(\zeta_{N})},\underline{\omega}_{\mathfrak{X}(N)_{/\mathbb{Q}(\zeta_{N})}}^{\otimes k}).\]
In particular, the coefficients of any $q$-expansion of a modular form over $\mathbb{Q}(\zeta_{N})$ have bounded denominators. Therefore we can write the $q$-expansion of any $f\in M_{k}(N,\mathbb{Z}[\zeta_{N},1/N])$ at $\infty$ as $f_{\infty}=\frac{1}{N^{m}}f'$ for some $f'\in\mathbb{Z}[\zeta_{N}][[q^{1/N}]]$ and integer $m\geq0$, \emph{i.e. }$f'\in M_{k}(N,\mathbb{Z}[\zeta_{N}])$. Then $f$ is the image of $f'\otimes\frac{1}{N^{m}}$.
\end{proof}

\begin{defn}
Let $M$ be a module over a ring $R$. Then \textbf{annihilator }of $M$ is the ideal 
\[
{\rm Ann}_{R}(M):=\left\{ r\in R:rm=0\mbox{ for all }m\in M\right\} .
\]
The \textbf{annihilator }of an element $m\in M$ is the ideal 
\[ {\rm Ann}_{R}(m):=\left\{ r\in R:rm=0\right\} .\]
\end{defn}

The next proposition showcases some properties of our two $\mathbb{Z}[\zeta_{N}]$-modules necessary to discuss the annihilator of the quotient $M_{k}(N,\mathbb{Z}[\zeta_{N}])/H^{0}(\mathfrak{X}(N),\underline{\omega}^{\otimes k})$. 

\begin{prop}\leavevmode\label{Prop of int structs}
\begin{enumerate}
    \item[a.] Both $H^{0}(\mathfrak{X}(N),\underline{\omega}^{\otimes k})$ and $M_{k}(N,\mathbb{Z}[\zeta_{N}])$ are finitely generated $\mathbb{Z}[\zeta_{N}]$-modules of the same rank.
    
    \item[b.] We have \[H^{0}(\mathfrak{X}(N),\underline{\omega}^{\otimes k})\otimes_{\mathbb{Z}[\zeta_{N}]}\mathbb{Z}[\zeta_{N},1/N]=M_{k}(N,\mathbb{Z}[\zeta_{N}])\otimes_{\mathbb{Z}[\zeta_{N}]}\mathbb{Z}[\zeta_{N},1/N].\]
\end{enumerate}
\end{prop}

\begin{proof}
Since $\mathfrak{X}(N)$ is projective over $\mathbb{Z}[\zeta_{N}]$, $H^{0}(\mathfrak{X}(N),\underline{\omega}^{\otimes k})$ is a finitely generated $\mathbb{Z}[\zeta_{N}]$-module by \cite[Theorem 5.3.2]{[Liu02]}. Since $M_{k}(N,\mathbb{Z}[\zeta_{N}])$ is a submodule of a finitely generated module over a noetherian ring, namely of $H^{0}(\mathfrak{X}(N),\underline{\omega}^{\otimes k})$, it is also finitely generated. 

Now we show the second claim. By Proposition \ref{Base change}, we have 
\[
H^{0}(\mathfrak{X}(N),\underline{\omega}^{\otimes k})\otimes_{\mathbb{Z}[\zeta_{N}]}\mathbb{Z}[\zeta_{N},1/N]=H^{0}(\mathfrak{X}(N)_{/\mathbb{Z}[\zeta_{N},1/N]},\underline{\omega}^{\otimes k})
\]
while, by Lemma \ref{M2 base change}, we have 
\[
M_{k}(N,\mathbb{Z}[\zeta_{N}])\otimes_{\mathbb{Z}[\zeta_{N}]}\mathbb{Z}[\zeta_{N},1/N]=M_{k}(N,\mathbb{Z}[\zeta_{N},1/N]).
\]
Since modular forms over $\mathbb{Z}[\zeta_{N},1/N]$ are determined by their $q$-expansion at $\infty$ (see \cite[{}1.6.2]{[Kat72]}), we get an equality 
\[
H^{0}(\mathfrak{X}(N)_{/\mathbb{Z}[\zeta_{N},1/N]},\underline{\omega}^{\otimes k})=M_{k}(N,\mathbb{Z}[\zeta_{N},1/N])
\]
as desired. In particular, we have an equality of finite dimensional vector spaces 
\[
H^{0}(\mathfrak{X}(N),\underline{\omega}^{\otimes k})\otimes_{\mathbb{Z}[\zeta_{N}]}\mathbb{Q}(\zeta_{N})=M_{k}(N,\mathbb{Z}[\zeta_{N}])\otimes_{\mathbb{Z}[\zeta_{N}]}\mathbb{Q}(\zeta_{N})
\]
which shows both our $\mathbb{Z}[\zeta_{N}]$-modules have the same rank. 
\end{proof}

Our description of $H^{0}(\mathfrak{X}(N),\underline{\omega}^{\otimes k})$ in Lemma \ref{Description of H0(N)} shows that it is contained in $M_{k}(N,\mathbb{Z}[\zeta_{N}])$. Having established that $H^{0}(\mathfrak{X}(N),\underline{\omega}^{\otimes k})$ and $M_{k}(N,\mathbb{Z}[\zeta_{N}])$ have the same rank, the quotient
\[
M_{k}(N,\mathbb{Z}[\zeta_{N}])/H^{0}(\mathfrak{X}(N),\underline{\omega}^{\otimes k})
\]
is a torsion $\mathbb{Z}[\zeta_{N}]$-module. Our goal will be to investigate the annihilator of this quotient. 

By localizing, we can focus our attention on investigating a single prime in the annihilator. The following lemma, which is \cite[Tag080S]{[Stacks]}, shows how the annihilator behaves under flat base change.

\begin{lem}\label{Ann Base Change}
Let $R\rightarrow S$ be a flat ring map. Let $M$ be an $R$-module and $m\in M$. Then 
\[
{\rm Ann}_{S}(m\otimes1)={\rm Ann}_{R}(m)S
\]
for any $m\otimes1\in M\otimes_{R}S$. If $M$ is finite over $R$, then 
\[
{\rm Ann}_{S}(M\otimes_{R}S)={\rm Ann}_{R}(M)S.
\]
\end{lem}

\begin{cor}\label{Exp of Ann Localized}
Let $M$ be a finitely generated torsion module over a Dedekind domain $R$ and let $\mathfrak{p}$ be a prime in $R$. The exponent of $\mathfrak{p}$ appearing in ${\rm Ann}_{R}(M)$ is equal to the exponent of $\mathfrak{p}$ appearing in ${\rm Ann}_{R_{\mathfrak{p}}}(M\otimes_{R}R_{\mathfrak{p}})$.
\end{cor}

\begin{proof}
We factor the annihilator as a product of distinct primes ${\rm Ann}_{R}(M)=\mathfrak{p}^{e}\mathfrak{q}_{1}^{e_{1}}\dots\mathfrak{q}_{m}^{e_{m}}$. By Lemma \ref{Ann Base Change}, we have 
\[
{\rm Ann}_{R_{\mathfrak{p}}}(M\otimes_{R}R_{\mathfrak{p}})={\rm Ann}_{R}(M)R_{\mathfrak{p}}=\mathfrak{p}^{e}\mathfrak{q}_{1}^{e_{1}}\dots\mathfrak{q}_{m}^{e_{m}}R_{\mathfrak{p}}=\mathfrak{p}^{e}R_{\mathfrak{p}}.
\] Hence the exponent of $\mathfrak{p}$ appearing in ${\rm Ann}_{R_{\mathfrak{p}}}(M\otimes_{R}R_{\mathfrak{p}})$ is precisely $e$.
\end{proof}

\begin{prop}
Let $\mathfrak{p}\in\mathbb{Z}[\zeta_{N}]$ be a prime lying over the prime $p\in\mathbb{Z}$. The exponent of $\mathfrak{p}$ appearing in the annihilator of $M_{k}(N,\mathbb{Z}[\zeta_{N}])/H^{0}(\mathfrak{X}(N),\underline{\omega}^{\otimes k}))$ and $M_{k}(N,\mathbb{Z}_{p}[\zeta_{N}])/H^{0}(\mathfrak{X}(N)_{/\mathbb{Z}_{p}[\zeta_{N}]},\underline{\omega}^{\otimes k})$ are the same.
\end{prop}

\begin{proof}
For convenience, we let 
\[
M=M_{k}(N,\mathbb{Z}[\zeta_{N}]),\;N=H^{0}(\mathfrak{X}(N),\underline{\omega}^{\otimes k})
\]
and
\[
M_{p}=M_{k}(N,\mathbb{Z}_{p}[\zeta_{N}]),\;N_{p}=H^{0}(\mathfrak{X}(N)_{/\mathbb{Z}_{p}[\zeta_{N}]},\underline{\omega}^{\otimes k}).
\]
Consider the exact sequence $0\rightarrow N\rightarrow M\rightarrow M/N\rightarrow0.$ Since $\mathbb{Z}_{p}[\zeta_{N}]$ is flat over $\mathbb{Z}[\zeta_{N}]$, the above sequence remains exact after tensoring with $\mathbb{Z}_{p}[\zeta_{N}]$. Thus we have a commutative diagram 
\[
\xymatrix{
0 \ar[r] & N\otimes\mathbb{Z}_{p}[\zeta_{N}]\ar[r]\ar[d]^{f} 
& M\otimes\mathbb{Z}_{p}[\zeta_{N}]\ar[r]\ar[d]^{g} 
& (M/N)\otimes\mathbb{Z}_{p}[\zeta_{N}]\ar[r]\ar[d]^{h} 
& 0\\
0\ar[r] & N_{p}\ar[r] & M_{p}\ar[r] & M_{p}/N_{p}\ar[r] & 0
}
\]
By Proposition \ref{Base change}, $f$ and $g$ are isomorphisms; consequently $h$ is an isomorphism. By Corollary \ref{Exp of Ann Localized}, the exponent of $\mathfrak{p}$ appearing in $M_{p}/N_{p}\simeq(M/N)\otimes\mathbb{Z}_{p}[\zeta_{N}]$ agrees with that of $M/N$.
\end{proof}

\begin{rmk}\label{Trivial Ann}
If $p\nmid N$, then $N$ is invertible in $\mathbb{Z}_{p}[\zeta_{N}]$ so by Proposition \ref{Prop of int structs} we have 
\[
H^{0}(\mathfrak{X}(N),\underline{\omega}^{\otimes k})\otimes_{\mathbb{Z}[\zeta_{N}]}\mathbb{Z}_{p}[\zeta_{N}]=M_{k}(N,\mathbb{Z}[\zeta_{N}])\otimes_{\mathbb{Z}[\zeta_{N}]}\mathbb{Z}_{p}[\zeta_{N}].
\]
Hence the $p$-adic valuation of the annihilator of $M_{k}(N,\mathbb{Z}[\zeta_{N}])/H^{0}(\mathfrak{X}(N),\underline{\omega}^{\otimes k})$ is trivial in this case. Therefore we restrict our attention to primes $p$ dividing the level.
\end{rmk}

Let $N\geq3$ and $r\geq0$ be integers and let $p\geq2$ be a prime such that $p\nmid N$. Let $\pi$ be a uniformizer of $\mathbb{Z}_{p}[\zeta_{Np^{r}}]$. We seek to compute the smallest integer $e\geq0$ such that \[
\pi^{e}H^{0}(\mathfrak{X}(Np^{r})_{/\mathbb{Z}_{p}[\zeta_{Np^{r}}]},\underline{\omega}^{\otimes k})\subseteq M_{k}(Np^{r},\mathbb{Z}_{p}[\zeta_{Np^{r}}]).
\]

\subsection{\label{sub: geo exp}Geometric interpretation of the exponent}

For convenience, we let $\mathfrak{X}=\mathfrak{X}(Np^{r})_{/\mathbb{Z}_{p}[\zeta_{Np^{r}}]}$. We will formulate an explicit, geometric description of the exponent $e$ by first providing an algebraic description coming from $q$-expansions. 

Let $c$ be a cusp of $\mathfrak{X}$, and let $f\in H^{0}(\mathfrak{X},\underline{\omega}^{\otimes k})$ be a non-zero modular form. Let $f_{c}\in\mathbb{Z}_{p}[\zeta_{Np^{r}}][[q^{1/N}]]$ denote the $q$-expansion of $f$ at the cusp $c$. Define $\nu_{\pi}(f_{c})$ to be the minimal $\pi$-adic valuation among all the coefficients of $f_{c}$ \emph{i.e.} 
\[
\nu_{\pi}\left(\sum_{n=0}^{\infty}a_{n}q^{n/(Np^{r})}\right)=\min_{n\geq0}\left\{ \nu_{\pi}(a_{n})\right\} .
\]
Note that this minimum exists since the denominators of the coefficients $a_{n}$ are bounded (see proof of Lemma \ref{M2 base change}). Thus, by Lemma \ref{Description of H0(N)}, we may describe 
\[
H^{0}(\mathfrak{X},\underline{\omega}^{\otimes k})=\left\{ f\in H^{0}(\mathfrak{X}_{/\mathbb{Q}_{p}(\zeta_{Np^{r}})},\underline{\omega}^{\otimes k}):\nu_{\pi}(f_{c})\geq0\mbox{ for all cusps }c\mbox{ of }\mathfrak{X}\right\} .
\]
Consequently, we seek to compute the smallest integer $e\geq0$ such that $\nu_{\pi}(\pi^{e}f_{c})\geq0$ for all cusps $c$ of $\mathfrak{X}$ not equal to $\infty$ and all $f\in M_{k}(Np^{r},\mathbb{Z}_{p}[\zeta_{Np^{r}}])$. 

Next we will provide a geometric interpretation of $\nu_{\pi}(f_{c})$ which does not rely on $q$-expansions. Let $\eta_{\Lambda}$ denote the generic point of an irreducible component $\Lambda$ of the special fiber $\bar{\mathfrak{X}}$. The stalk ${\cal O}_{\mathfrak{X},\eta}$ is a DVR so it has a valuation which we denote by $\nu_{\eta}$. Viewing $f$ inside the stalk $\underline{\omega}_{\eta}^{\otimes k}$, we can write $f=f_{\eta}\cdot\omega_{{\rm can},\eta}$ where $\omega_{{\rm can},\eta}$ is the canonical generator of $\underline{\omega}_{\eta}^{\otimes k}$ (see Definition \ref{def: q-exp}). We define
\[
\nu_{\Lambda}(f):=\nu_{\eta}(f_{\eta}).
\]
\noindent Indeed, $\nu_{\Lambda}$ is a valuation, independent of the choice of local generator. The following result is stated in \cite[Th\'{e}or\'{e}me 3.10(ii)]{[DeRa]}. 

\begin{prop}\label{Geom v. Alg Val}
Let $c\in\mathfrak{X}$ be a cusp and let $\Lambda$ be an irreducible component of $\bar{\mathfrak{X}}$ with generic point $\eta$ on which $c$ lies. For any non-zero $f\in H^{0}(\mathfrak{X}(N),\underline{\omega}^{\otimes k})$, we have $\nu_{\pi}(f_{c})=\nu_{\Lambda}(f)$. Furthermore, $\pi$ is a uniformizer of ${\cal O}_{\mathfrak{X},\eta}$.
\end{prop}

\begin{proof}
Since $c\in\overline{\left\{ \eta\right\} }$, the stalk ${\cal O}_{\mathfrak{X},\eta}$ is a localization of ${\cal O}_{\mathfrak{X},c}$. Furthermore, the map ${\cal O}_{\mathfrak{X},c}\rightarrow{\cal O}_{\mathfrak{X},\eta}$ is injective since all stalks are regular local rings so are integral domains in particular. The induced map on completions $\hat{{\cal O}}_{\mathfrak{X},c}\rightarrow\hat{{\cal O}}_{\mathfrak{X},\eta}$ is also injective by \cite[Tag00MB, Tag07N9]{[Stacks]}.

We will show $\pi$ is also a uniformizer of ${\cal O}_{\mathfrak{X},\eta}$. Consider the exact sequence 
\[
0\rightarrow\mathfrak{m}_{\eta}\rightarrow{\cal O}_{\mathfrak{X},\eta}\rightarrow\kappa(\eta)\rightarrow0
\]
where $\mathfrak{m}_{\eta}$ is the maximal ideal of ${\cal O}_{\mathfrak{X},\eta}$ and $\kappa(\eta)$ is the residue field of $\eta$ which is of characteristic $p$. Since $(\pi^{m})=(p)$ when $m=p^{r-1}(p-1)$, we have $\pi^{m}=0$ in $\kappa(\eta)$ so $\pi=0$ in $\kappa(\eta)$. Therefore $\pi\in\mathfrak{m}_{\eta}$. It remains to show $(\pi)$ is a maximal ideal in ${\cal O}_{\mathfrak{X},\eta}$. Consider the quotient ${\cal O}_{\mathfrak{X},\eta}/\pi{\cal O}_{\mathfrak{X},\eta}={\cal O}_{\bar{\mathfrak{X}},\eta}$. Since $\eta$ corresponds to a minimal prime and $\bar{\mathfrak{X}}$ is reduced, the stalk ${\cal O}_{\bar{\mathfrak{X}},\eta}$ is a field. Therefore $(\pi)$ is a maximal ideal in ${\cal O}_{\mathfrak{X},\eta}$ hence $\pi$ is a uniformizer in ${\cal O}_{\mathfrak{X},\eta}$. 

Let $f\in H^{0}(\mathfrak{X}(N),\underline{\omega}^{\otimes k})$ be non-zero. Using the local generator $\omega_{{\rm can},c}$ of $\underline{\hat{\omega}}^{\otimes k}$ (see Definition \ref{def: q-exp}), we write $f=f_{c}\omega_{{\rm can},c}$ where $f_{c}\in\hat{{\cal O}}_{\mathfrak{X},c}$ is the $q$-expansion of $f$ at $c$. Write $f_{c}=\pi^{\nu_{\pi}(f_{c})}f_{c}'$ where $f_{c}'\notin\pi\hat{{\cal O}}_{\mathfrak{X},c}$. The map 
\[
{\cal O}_{\mathfrak{X},c}/\pi{\cal O}_{\mathfrak{X},c}\rightarrow{\cal O}_{\mathfrak{X},\eta}/\pi{\cal O}_{\mathfrak{X,\eta}}
\]
is the same as the map ${\cal O}_{\bar{\mathfrak{X}},c}\rightarrow{\cal O}_{\bar{\mathfrak{X}},\eta}$, which is also a localization map as $c\in\overline{\left\{ \eta\right\} }$.

By \cite[10.9.1(2)]{[KM85]}, $\mathfrak{X}$ is smooth at the cusps so, in particular, ${\cal O}_{\bar{\mathfrak{X}},c}$ is a domain. Therefore ${\cal O}_{\bar{\mathfrak{X}},c}\rightarrow{\cal O}_{\bar{\mathfrak{X}},\eta}$ is injective so the map on completions 
\[
{\cal \hat{O}}_{\mathfrak{X},c}/\pi{\cal \hat{O}}_{\mathfrak{X},c}\rightarrow{\cal \hat{O}}_{\mathfrak{X},\eta}/\pi{\cal \hat{O}}_{\mathfrak{X,\eta}}
\]
is injective. Since $f_{c}'\notin\pi\hat{{\cal O}}_{\mathfrak{X},c}$, we can conclude $f_{c}'\notin\pi\hat{{\cal O}}_{\mathfrak{X},\eta}$ \emph{i.e. }$f_{c}'$ is a unit in $\hat{{\cal O}}_{\mathfrak{X},\eta}$. Thus $\nu_{\pi}(f_{c})$ coincides with $\nu_{\Lambda}(f)$, the valuation in ${\cal O}_{\mathfrak{X},\eta}$. 

\end{proof}

With this geometric description of the valuation, we can reformulate our integrality condition as follows. Let $\Lambda_{0},\dots,\Lambda_{n}$ be the irreducible components of the special fiber $\bar{\mathfrak{X}}$ where $\Lambda_{0}$ contains the cusp $\infty$. Let $\mu_{\Lambda_{i}}$ denote the multiplicity of $\Lambda_{i}$. We desire to find the smallest integer $e\geq0$ such that $\nu_{\Lambda_{i}}(\pi^{e}f)\geq0$ for all $i=1,\dots,n$ and all $f\in M_{k}(Np^{r},\mathbb{Z}_{p}[\zeta_{Np^{r}}])$. 

By Proposition \ref{Geom v. Alg Val}, $\pi$ is a uniformizer of ${\cal O}_{\mathfrak{X},\eta_{i}}$, where $\eta_{i}$ is the generic point of $\Lambda_{i}$. Hence
\begin{align*}
\nu_{\Lambda_{i}}(\pi^{e}f) & =  e\cdot\nu_{\Lambda_{i}}(\pi)+\nu_{\Lambda_{i}}(f)\\
& =  e\cdot\mu_{\Lambda_{i}}+\nu_{\Lambda_{i}}(f).
\end{align*}
Since $\bar{\mathfrak{X}}$ is reduced, $\mu_{\Lambda_{i}}=1$. Hence the condition $\nu_{\Lambda_{i}}(\pi^{e}f)\geq0$ is equivalent to $e\geq-\nu_{\Lambda_{i}}(f)$. Therefore we have the following expression for the exponent: 

\begin{equation}
e=\max_{\substack{f\in M_{k}(Np^{r},\mathbb{Z}_{p}[\zeta_{Np^{r}}])}
}\left\{ -\nu_{\Lambda_{1}}(f),-\nu_{\Lambda_{2}}(f),\dots,-\nu_{\Lambda_{n}}(f)\right\} \label{eq: exponent initially}
\end{equation}
which illustrates we can investigate $e$ by trying to understand the quantities $\nu_{\Lambda_{i}}(f)$.

\subsection{\label{sub:int_theory}Intersection theory on arithmetic surfaces}

In this section we briefly recall some important intersection theory facts from \cite[§9]{[Liu02]}. By Theorem \ref{X(N) is arith surf}, the modular curve $\mathfrak{X}(N)$ is an arithmetic surface (see Definition \ref{arith surf}). We begin with a more general situation with $\mathfrak{X}$ an arithmetic surface over a Dedekind domain $R$ and ${\cal L}$ a line bundle on $\mathfrak{X}$.

Since $\mathfrak{X}\rightarrow S={\rm Spec}(R)$ is regular, we have, by \cite[{}7.2.16]{[Liu02]}, an isomorphism between Cartier divisors and Weil divisors ${\rm Div}(\mathfrak{X})\rightarrow Z^{1}(\mathfrak{X})$ given by 
\[
D\mapsto\left[D\right]:=\sum_{\substack{x\in\mathfrak{X}\\
{\rm dim}{\cal O}_{\mathfrak{X},x}=1
}
}{\rm mult}_{x}(D)[\overline{\left\{ x\right\} }]\in Z^{1}(\mathfrak{X})
\]
where ${\rm mult}_{x}(D):={\rm mult}_{{\cal O}_{\mathfrak{X},x}}(D_{x})$. Furthermore, this map respects principal divisors and effective divisors. Let $Z$ be a closed subscheme of $\mathfrak{X}$ and let $\xi_{1},\dots,\xi_{n}$ be the generic points of $Z$. The Weil divisor associated to $Z$ is given by 
\[
\left[Z\right]:=\sum_{i=1}^{n}{\rm length}({\cal O}_{\mathfrak{X},\xi_{i}})[\overline{\left\{ \xi_{i}\right\} }].
\]
\noindent Let $D$ be an irreducible Weil divisor of $\mathfrak{X}$. According to \cite[{}8.3.4]{[Liu02]}, $D$ is either an irreducible component of a closed fiber, or the closure of a closed point of the generic fiber $X_{\eta}$. 

\begin{defn}
If $D$ is an irreducible component of a closed fiber, or equivalently $\pi(D)$ is a point, then $D$ is called \textbf{vertical}. If $D$ is instead the closure of a closed point of $X_{\eta}$, or equivalently $\pi(D)=S$, then $D$ is called \textbf{horizontal}.

Recall a Weil divisor is a formal $\mathbb{Z}$-linear sum of irreducible closed subsets of codimension 1. In general, a Weil divisor is horizontal (resp. vertical) if all its irreducible components are horizontal (resp. vertical).
\end{defn}

Let $D$ and $E$ be two effective divisors of $\mathfrak{X}$ with no common irreducible component.  As in \cite[§9.1.1]{[Liu02]}, we define the \textbf{(local) intersection number }of $D$ and $E$ at a point $x\in\mathfrak{X}$ by 
\[
i_{x}(D,E)={\rm length}_{{\cal O}_{\mathfrak{X},x}}{\cal O}_{\mathfrak{X},x}/({\cal O}_{\mathfrak{X}}(-D)_{x}+{\cal O}_{\mathfrak{X}}(-E)_{x}).
\]
By definition, we immediately have that $i_{x}(D,E)$ is symmetric and bilinear. If $x\notin{\rm Supp}(D)\cap{\rm Supp}(E)$, then $i_{x}(D,E)=0$.

We will now establish the intersection number between a general divisor of $\mathfrak{X}$ and a vertical divisor of a fixed closed fiber. For a fixed closed point $s\in S$, let ${\rm Div}_{s}(\mathfrak{X})$ denote the set of divisors of $\mathfrak{X}$ with support in $\mathfrak{X}_{s}$; such divisors are vertical divisors. A $\mathbb{Z}$-basis of ${\rm Div}_{s}(\mathfrak{X})$ consists of all irreducible components of $\mathfrak{X}_{s}$.

\begin{thm}\label{Props of Int Num}
Let $s\in S$ be a closed point. Then there
exists a unique bilinear map of $\mathbb{Z}$-modules 
\[
i_{s}:{\rm Div}(\mathfrak{X})\times{\rm Div}_{s}(\mathfrak{X})\rightarrow\mathbb{Z}
\]
such that:

\begin{enumerate}
    \item[a.] If $D\in{\rm Div}(\mathfrak{X})$ and $E\in{\rm Div}_{s}(\mathfrak{X})$ have no common component, then \[i_{s}(D,E)=\sum_{x\in\mathfrak{X}_{s}}i_{x}(D,E)[k(x):k(s)]\] where the sum is over all closed points.
    \item[b.]  $i_{s}$ restricted to ${\rm Div}_{s}(\mathfrak{X})\times{\rm Div}_{s}(\mathfrak{X})$ is symmetric.
    \item[c.] $i_{s}(D,E)=i_{s}(D',E)$ if $D\sim D'$ are linearly equivalent.
    \item[d.] If $0<E\leq\mathfrak{X}_{s}$, then $i_{s}(D,E)=\deg_{k(s)}{\cal O}_{\mathfrak{X}}(D)|_{E}$. $x\in\mathfrak{X}_{s}$.
    \item[e.] If $D$ is principal, then $i_{s}(D,E)=0$.
\end{enumerate}
\end{thm}

\begin{proof}
This is proved in \cite[Theorem 9.1.12]{[Liu02]} besides part (e), which we prove.  By \cite[Corollary 9.1.10]{[Liu02]}, there exists a principal divisor $F$ such that $D+F$ and $E$ have no common component. Thus we may assume $D$ and $E$ have no common component. By (a), it suffices to show if $D$ is principal, then the local intersection number $i_{x}(D,E)$ is zero. By definition, 
\[
i_{x}(D,E)={\rm length}_{{\cal O}_{\mathfrak{X},x}}{\cal O}_{\mathfrak{X},x}/({\cal O}_{\mathfrak{X}}(-D)_{x}+{\cal O}_{\mathfrak{X}}(-E)_{x}).
\]
Since $D$ is principal, ${\cal O}_{\mathfrak{X}}(-D)\simeq{\cal O}_{\mathfrak{X}}$. Thus 
\begin{alignat*}{1}
i_{x}(D,E) & ={\rm length}_{{\cal O}_{\mathfrak{X},x}}{\cal O}_{\mathfrak{X},x}/({\cal O}_{\mathfrak{X},x}+{\cal O}_{\mathfrak{X}}(-E)_{x})\\
 & ={\rm length}_{{\cal O}_{\mathfrak{X},x}}(0)=0.\tag*{\qedhere}
\end{alignat*}
\end{proof}

\begin{defn}
Let $\pi:\mathfrak{X}\rightarrow S$ be an arithmetic surface, and let $s\in S$ be closed. For any $(D,E)$ in ${\rm Div}(\mathfrak{X})\times{\rm Div}_{s}(\mathfrak{X})$, we call $i_{s}(D,E)$ the \textbf{intersection number }as defined in Theorem \ref{Prop of int structs}. More generally, if $E$ is a vertical divisor, we define 
\[
D.E:=\sum_{s\in S}i_{s}(D,E)[s]
\]
where the sum is over all closed points of $S$, which is a $0$-cycle on $S$. We call $E^{2}=E.E$ the \textbf{self-intersection number}of $E$. If $D.E$ is only concentrated at a single point $s$, as in the case when $S$ is the spectrum of a DVR, we identify $D.E$ with the integer $i_{s}(D,E)$.
\end{defn}

\begin{prop}\label{Special Fiber as Div}
Let $\pi:\mathfrak{X}\rightarrow S$ be an arithmetic surface and fix a closed point $s\in S$. Let $\Gamma_{1},\dots\Gamma_{n}$ be the irreducible components of $\mathfrak{X}_{s}$. We have an equality of Weil divisors 
\[
[\mathfrak{X}_{s}]=\sum_{i=1}^{n}{\rm mult}_{\mathfrak{X}_{s}}(\Gamma_{i})\Gamma_{i}=[\pi^{*}s]
\]
where $[\pi^{*}s]$ is the Weil divisor associated to the Cartier divisor $\pi^{*}s$.
\end{prop}

\begin{proof}
The first equality follows from the definition of $\left[\mathfrak{X}_{s}\right]$. Let $t\in{\cal O}_{S,s}$ and $t_{i}\in{\cal O}_{\mathfrak{X},\xi_{i}}$ be uniformizers where $\xi_{i}$ is the generic point of $\Gamma_{i}$. Let $\nu_{i}$ denote the normalized valuation of $K(\mathfrak{X})$ associated to ${\cal O}_{\mathfrak{X},\xi_{i}}$ and let $\xi$ denote the generic point of $\mathfrak{X}$. 

In $K(\mathfrak{X})={\cal O}_{\mathfrak{X},\xi}$, we can write $t=t_{i}^{\nu_{i}(t)}u_{i}$ for some $u_{i}\in{\cal O}_{\mathfrak{X},\xi_{i}}^{\times}$. Then
\[
[\pi^{*}s]=\sum_{i=1}^{n}\nu_{i}(t)\Gamma_{i}.
\]
Note that
\begin{alignat*}{1}
{\rm mult}_{\mathfrak{X}_{s}}(\Gamma_{i}) & ={\rm length}({\cal O}_{\mathfrak{X}_{s},\xi_{i}})\\
 & ={\rm length}({\cal O}_{\mathfrak{X},\xi_{i}}/(t))\\
 & ={\rm length}({\cal O}_{\mathfrak{X},\xi_{i}}/(t_{i}^{\nu_{i}(t)}))\\
 & =\nu_{i}(t).\tag*{\qedhere}
\end{alignat*}
\end{proof}




Lastly we state some results on the intersection between a horizontal divisor and a closed fiber which won't be used until §\ref{sub:Computing deg(omega|lambda)}. The following is \cite[Proposition 9.1.30]{[Liu02]}.

\begin{prop}\label{Int Horiz with Xs}
Let $\pi:\mathfrak{X}\rightarrow S$ be an arithmetic surface. Let $\eta$ be the generic point of $S$ and $s\in S$ a closed point. Then for any closed point $P\in\mathfrak{X}_{\eta}$, we have 
\[
\overline{\left\{ P\right\} }.\mathfrak{X}_{s}=\left[K(P):K(S)\right]
\]
where $\overline{\left\{ P\right\} }$ is the Zariski closure of $\left\{ P\right\}$ in $\mathfrak{X}$, endowed with the reduced closed subscheme structure.
\end{prop}

\begin{cor}\label{Int of Ratl Pt}
Let $\pi:\mathfrak{X}\rightarrow S$ be an arithmetic surface and let $P\in\mathfrak{X}_{\eta}$ be a $K(S)$-rational point. Then $\overline{\left\{ P\right\} }\cap\mathfrak{X}_{s}$ is reduced to a single point $z\in\mathfrak{X}_{s}(k(s))$ and $\overline{\left\{ P\right\} }.\mathfrak{X}_{s}=1$. 

In particular, $\overline{\left\{ P\right\} }$ intersects exactly one irreducible component $\Gamma$ of $\mathfrak{X}_{s}$. Moreover, $\Gamma$ has multiplicity $1$ and $\overline{\left\{ P\right\} }.\Gamma=1$. 
\end{cor}

\begin{proof}
Since $P$ is $K(S)$-rational, by Proposition \ref{Int Horiz with Xs} we have $\overline{\left\{ P\right\} }.\mathfrak{X}_{s}=1$. Let $\Gamma_{1},\dots,\Gamma_{n}$ denote the irreducible components of $\mathfrak{X}_{s}$ with multiplicities $\mu_{1},\dots,\mu_{n}$ respectively. Then 
\[
1=\overline{\left\{ P\right\} }.\mathfrak{X}_{s}=\sum_{i=1}^{n}\mu_{i}\overline{\left\{ P\right\} }.\Gamma_{i}.
\]
Since $\overline{\left\{ P\right\} }$ and $\Gamma_{i}$ are effective divisors with no common components, $\overline{\left\{ P\right\} }.\Gamma_{i}\geq0$ which forces $\mu_{i}\overline{\left\{ P\right\} }.\Gamma_{i}=1$ for some $i$. Hence $\mu_{i}=\overline{\left\{ P\right\} }.\Gamma_{i}=1$. Moreover, $\mu_{j}\geq1$ so $\overline{\left\{ P\right\} }.\Gamma_{j}=0$ for $j\ne i$ i.e. $\overline{\left\{ P\right\} }$ does not intersect $\Gamma_{j}$.
\end{proof}

\subsection{\label{sub:description_of_e}A more explicit description of the exponent}

In this section, we will provide a more explicit description of the exponent by relating the quantities $\nu_{\Lambda}(f)$ to intersection numbers and the degree of a line bundle, due to Conrad in \cite[Appendix B]{[Con17]}. Let $\mathfrak{X}$ be an arithmetic surface over a DVR $R$. Let ${\cal L}$ be a line bundle on $\mathfrak{X}$ and let $\bar{\mathfrak{X}}$ denote the special fiber of $\mathfrak{X}$. 

\begin{defn}
Let $f\in H^{0}(\mathfrak{X},{\cal L})$ be a non-zero global section and let $\left\{ U_{i}\right\} $ be a trivializing open cover of ${\cal L}$ so that ${\cal L}(U_{i})={\cal O}_{X}(U_{i})e_{i}$ for some generator $e_{i}\in{\cal L}(U_{i})$. We can write $f|_{U_{i}}=f_{i}e_{i}$ for some $f_{i}\in{\cal O}_{\mathfrak{X}}(U_{i})$. Then the system $\left\{ (U_{i},f_{i})\right\} $ is an effective Cartier divisor of $\mathfrak{X}$, which we denote by ${\rm div}_{{\cal L}}(f)$ or simply ${\rm div}(f)$ if the line bundle is clear from context. 
\end{defn}

We can equivalently define ${\rm div}_{{\cal L}}(f)$ as a Weil divisor. Let $Z$ be a prime divisor of $\mathfrak{X}$ with generic point $\eta$. Then ${\cal O}_{\mathfrak{X},\eta}$ is a DVR with valuation which we denote by $\nu_{Z}$. We can write the image of $f$ in the stalk of $x$ as $f=f_{\eta}e_{\eta}$ for some $f_{\eta}\in{\cal O}_{X,\eta}$ and generator $e_{\eta}\in{\cal L}_{\eta}$. Define $\nu_{Z}(f):=\nu_{Z}(f_{\eta})$ which agrees with our valuation defined in Section \ref{sub: geo exp} for $Z=\Lambda$ an irreducible component of $\bar{\mathfrak{X}}$ and is independent of choice of $e_{\eta}$. According to \cite[Tag02SE]{[Stacks]}, the Weil divisor associated to $f$ is equal to 
\[
{\rm div}_{{\cal L}}(f)=\sum_{Z}\nu_{Z}(f)[Z]
\]
where the sum is over prime divisors $Z$ of $\mathfrak{X}$. Decomposing the divisor ${\rm div}_{{\cal L}}(f)$ into its horizontal and vertical components, we can write 
\begin{equation}
{\rm div}_{{\cal L}}(f)=H_{f}+\sum_{\Lambda}\nu_{\Lambda}(f)\left[\Lambda\right]\label{eq: div(f)=00003DH+V}
\end{equation}
where $H_{f}$ is some effective horizontal divisor and the sum is over the irreducible components of $\bar{\mathfrak{X}}$. 

\begin{prop}\label{twist by div(f)}
Let $f\in H^{0}(\mathfrak{X},{\cal L})$ be a non-zero global section. Then ${\cal L}\simeq{\cal O}_{\mathfrak{X}}({\rm div}_{{\cal L}}(f))$.
\end{prop}

\begin{proof}
This is \cite[Exercise 7.1.13]{[Liu02]} (see also immediately before \cite[Proposition B.2.2.10]{[Con17]}).
\end{proof}

Let $\Lambda_{0},\dots,\Lambda_{n}$ denote the irreducible components of $\bar{\mathfrak{X}}$. By Proposition \ref{twist by div(f)} and Theorem \ref{Props of Int Num}d, we have

\begin{align*}
\deg_{k}({\cal L}|_{\Lambda_{j}}) & = \deg({\cal O}_{\mathfrak{X}}({\rm div}_{{\cal L}}(f))|_{\Lambda_{j}})\\
& = {\rm div}_{{\cal L}}(f).\Lambda_{j}.
\end{align*}
Decomposing ${\rm div}_{{\cal L}}(f)$ into its horizontal and vertical components, as in (\ref{eq: div(f)=00003DH+V}), we get 
\begin{align*}
& = \left(H_{f}+\sum_{i=0}^{n}\nu_{\Lambda_{i}}(f)\left[\Lambda_{i}\right]\right).\Lambda_{j}\\
& = H_{f}.\Lambda_{j}+\sum_{i=0}^{n}\nu_{\Lambda_{i}}(f)\Lambda_{i}.\Lambda_{j}.
\end{align*}
Using the equation 
\begin{equation}
\deg({\cal L}|_{\Lambda_{j}})=H_{f}.\Lambda_{j}+\sum_{i=0}^{n}\nu_{\Lambda_{i}}(f)\Lambda_{i}.\Lambda_{j}\label{main eqn}
\end{equation}
we will provide an explicit expression for the quantities $\nu_{\Lambda_{i}}(f)/\mu_{\Lambda_{i}}$ for $i=1,\dots,n$ where $\mu_{\Lambda_{i}}$ is the multiplicity of $\Lambda_{i}$. 

Let $M=(\Lambda_{i}.\Lambda_{j})_{i,j=0,\dots,n}$ denote the intersection matrix of $\mathfrak{X}$. Since $\mathfrak{X}$ is an arithmetic surface, \cite[{}9.1.23]{[Liu02]} (see also \cite[III, 3.4]{[Lan88]}) says $M$ is negative semi-definite and moreover the kernel of $M$ is one-dimensional whenever $\bar{\mathfrak{X}}$ is connected. Multiplication by $M$ induces an exact sequence 
\[
0\rightarrow\ker M\rightarrow\mathbb{Q}^{n+1}\overset{M}{\longrightarrow}\mathbb{Q}^{n+1}.
\]
Let $\vec{\mu}=(\mu_{\Lambda_{0}},\mu_{\Lambda_{1}},\dots,\mu_{\Lambda_{n}})$ which is a non-zero vector in $\mathbb{Q}^{n+1}$. Note that 
\[
(M\vec{\mu})_{j}=\sum_{i=0}^{n}\mu_{\Lambda_{i}}(\Lambda_{i}.\Lambda_{j})=\left(\sum_{i=0}^{n}\mu_{\Lambda_{i}}\Lambda_{i}\right).\Lambda_{j}=\bar{\mathfrak{X}}.\Lambda_{j}=0
\]
where the last equality is due to $\bar{\mathfrak{X}}$ being principal. Therefore we can write $\ker(M)={\rm span}_{\mathbb{Q}}\left\{ \vec{\mu}\right\} $. The following lemma describes the image of $M$.

\begin{lem}\label{im(M)}
We have 
\[
{\rm im}(M)=\left\{ \vec{a}=(a_{0},\dots,a_{n})\in\mathbb{Q}^{n+1}:\sum_{j=0}^{n}\mu_{\Lambda_{j}}a_{j}=0\right\}.
\]
\end{lem}
\begin{proof}
This is shown in the paragraph preceding \cite[Remark B.2.3.1]{[Con17]}.  We will provide the proof here for convenience.  Let $V$ denote the space of vectors $\vec{v}=(v_{i})\in\mathbb{Q}^{n+1}$ such that $\sum_{j=0}^{n}\mu_{\Lambda_{j}}v_{j}=0$. Let $\vec{a}=(a_{j})\in{\rm im}(M)$ so $\vec{a}=M\vec{b}$ for some $\vec{b}=(b_{j})\in\mathbb{Q}^{n+1}$. Then 
\[
a_{j}=(M\vec{b})_{j}=\sum_{i=0}^{n}b_{i}(\Lambda_{i}.\Lambda_{j})
\]
so
\begin{align*}
\sum_{j=0}^{n}\mu_{\Lambda_{j}}a_{j} & =  \sum_{j=0}^{n}\mu_{\Lambda_{j}}\sum_{i=0}^{n}b_{i}(\Lambda_{i}.\Lambda_{j})=\sum_{j=0}^{n}\mu_{\Lambda_{j}}\left(\sum_{i=0}^{n}b_{i}\Lambda_{i}\right).\Lambda_{j}\\
& =  \left(\sum_{i=0}^{n}b_{i}\Lambda_{i}\right).\left(\sum_{j=0}^{n}\mu_{\Lambda_{j}}\Lambda_{j}\right)=\left(\sum_{i=0}^{n}b_{i}\Lambda_{i}\right).\bar{X}=0.\\
\\
\end{align*}
Therefore ${\rm im}(M)\subseteq V$. Since ${\rm im}(M)$ and $V$ are both of dimension $n$, we get ${\rm im}(M)=V$.
\end{proof}

In order to isolate the terms $\nu_{\Lambda_{i}}(f)$ appearing in (\ref{main eqn}), we would ideally invert the matrix $M$. In light of $M$ being not invertible, we will instead invert an $n\times n$ submatrix of $M$ to obtain an expression for each $\nu_{\Lambda_{i}}(f)$, excluding $\nu_{\Lambda_{0}}(f)$, which suffices for our purpose.

Let $W=\left\{ (x_{0},\dots,x_{n})\in\mathbb{Q}^{n+1}:x_{0}=0\right\} $ and let ${\rm pr}:\mathbb{Q}^{n+1}\rightarrow W$ denote the projection map 
\[
(x_{0},x_{1},\dots,x_{n})\mapsto(0,x_{1},\dots,x_{n}).
\]
Consider the restriction ${\rm pr}|_{{\rm im}(M)}:{\rm im}(M)\rightarrow W$. Given any $a_{1},\dots,a_{n}\in\mathbb{Q}$, we let 
\[
a_{0}=-\frac{1}{\mu_{\Lambda_{0}}}\sum_{i=1}^{n}\mu_{i}a_{i}
\]
which forces $\vec{a}:=(a_{0},a_{1},\dots,a_{n})\in{\rm im}(M)$ by Lemma \ref{im(M)}. Hence ${\rm pr}|_{{\rm im}(M)}$ is surjective. Since ${\rm im}(M)$ and $W$ are both of dimension $n$, the map ${\rm pr}|_{{\rm im}(M)}$ is an isomorphism of vector spaces. 

Next we consider the restriction $M|_{W}:W\rightarrow{\rm im}(M)$. If $M\vec{w}=0$ where $\vec{w}=(0,w_{1},\dots,w_{n})$, then $\vec{w}=\alpha\vec{\mu}$ for some $\alpha\in\mathbb{Q}$. Coordinatewise, this means 
\[
(0,w_{1},\dots,w_{n})=\vec{w}=\alpha\vec{\mu}=(\alpha\mu_{\Lambda_{0}},\alpha\mu_{\Lambda_{1}},\dots,\alpha\mu_{\Lambda_{n}}).
\]
Comparing the first coordinates and noting $\mu_{\Lambda_{i}}\ne0$ for each $i$, we must have $\alpha=0$. Hence $\vec{w}=0$ so $M|_{W}$ is injective and therefore an isomorphism. Define $T:W\rightarrow W$ as the composition of our isomorphisms
\[
T:W\overset{M|_{W}}{\longrightarrow}{\rm im}(M)\overset{{\rm pr}|_{{\rm im}(M)}}{\longrightarrow}W.
\]
After identifying $W\simeq\mathbb{Q}^{n}$ by forgetting the first coordinate, we can identify $T$ as the lower right $n\times n$ submatrix of $M$. The following is \cite[Proposition B.2.3.2]{[Con17]}.

\begin{prop}\label{T relation}
For any $\vec{a}\in{\rm im}(M)$ and $\vec{b}=(b_{0},b_{1},\dots,b_{n})\in\mathbb{Q}^{n+1}$ such that $M\vec{b}=\vec{a}$, we have 
\begin{equation}
\vec{b}-\frac{b_{0}}{\mu_{\Lambda_{0}}}\vec{\mu}=T^{-1}({\rm pr}(\vec{a})).\label{eq:T}
\end{equation}
\end{prop}

We now apply Proposition \ref{T relation} to our specific situation involving the quantities $\nu_{\Lambda_{i}}(f)$. 

\begin{rmk}\label{rmk on exp}
We make a quick remark about the term $\nu_{\Lambda_{0}}(f)$ that appears in the following theorem. Recall the equation for the exponent $e$ in (\ref{eq: exponent initially}). We claim that 
\[
e=\max_{\substack{1\leq i\leq n\\
f\in M_{k}(Np^{r},\mathbb{Z}_{p}[\zeta_{Np^{r}}])
}
}\left\{ -\nu_{\Lambda_{i}}(f)\right\} =\max_{\substack{1\leq i\leq n\\
f\in M_{k}(Np^{r},\mathbb{Z}_{p}[\zeta_{Np^{r}}])\\
\nu_{\Lambda_{0}}(f)=0
}
}\left\{ \nu_{\Lambda_{0}}(f)-\nu_{\Lambda_{i}}(f)\right\} .
\]
Indeed, the maximum on the left hand side must occur at some $f\in M_{k}(Np^{r},\mathbb{Z}_{p}[\zeta_{Np^{r}}])$ with $\nu_{\Lambda_{0}}(f)=0$. Otherwise if $f=\pi^{a}g$ with $g\in M_{k}(Np^{r},\mathbb{Z}_{p}[\zeta_{Np^{r}}])$, $\nu_{\Lambda_{0}}(g)=0$, and $a>0$, so 
\[
-\nu_{\Lambda_{i}}(g)=a-\nu_{\Lambda_{i}}(f).
\]
Hence $-\nu_{\Lambda_{i}}(g)\geq-\nu_{\Lambda_{i}}(f)$ so the maximum must occur over such $g$. 

Furthermore, the differences $\nu_{\Lambda_{0}}(f)-\nu_{\Lambda_{i}}(f)$ are visibly invariant under arbitrary scaling of $f$, so we have
\begin{align*}
e & = \max_{\substack{1\le i\leq n\\
f\in H^{0}(\mathfrak{X}_{/\mathbb{Q}_{p}(\zeta_{Np^{r}})},\underline{\omega}^{\otimes k})\\
\nu_{\Lambda_{0}}(f)=0
}
}\left\{ \nu_{\Lambda_{0}}(f)-\nu_{\Lambda_{i}}(f)\right\} .\\
 & = \max_{\substack{1\le i\leq n\\
f\in H^{0}(\mathfrak{X},\underline{\omega}^{\otimes k})
}
}\left\{ \nu_{\Lambda_{0}}(f)-\nu_{\Lambda_{i}}(f)\right\} 
\end{align*}
In particular, this shows the exponent can be computed using Theorem \ref{Equation for v(f)} below, which provides a formula for the differences $\nu_{\Lambda_{0}}(f)-\nu_{\Lambda_{i}}(f)$ in terms of geometric data. 
\end{rmk}

\begin{thm}\label{Equation for v(f)}
Let $\mathfrak{X}=\mathfrak{X}(Np^{r})_{/\mathbb{Z}_{p}[\zeta_{Np^{r}}]}$ and $f\in H^{0}(\mathfrak{X}(N),\underline{\omega}^{\otimes 2k})$ be non-zero and let $\Lambda_{0},\dots,\Lambda_{n}$ denote the irreducible components of $\bar{\mathfrak{X}}$. We have 
\begin{equation}
\nu_{\Lambda_{i}}(f)-\nu_{\Lambda_{0}}(f)=\sum_{j=1}^{n}(\deg(\underline{\omega}^{\otimes 2k}|_{\Lambda_{j}})-H_{f}.\Lambda_{j})c^{i,j}\label{eq: val deg int}
\end{equation}
where $c^{i,j}$ is the $(i,j)$ entry of $T^{-1}$ and $T$ is the matrix obtained by removing the first row and column of the intersection matrix of $\mathfrak{X}$.
\end{thm}

\begin{proof}
We will write out equation (\ref{eq:T}) entry-wise in the generality of Proposition \ref{T relation}, with $\mathfrak{X}$ an arithmetic surface, and ${\cal L}$ a line bundle on $\mathfrak{X}$. Then we apply our results to the situation of the modular curve. 

Let $\vec{a}=(a_{0},\dots,a_{n})$ and $\vec{b}=(b_{0},\dots,b_{n})$ such that $M\vec{b}=\vec{a}$. Then equation \ref{eq:T} gives 
\[
\left(\begin{array}{c}
0\\
b_{1}-\frac{b_{0}}{\mu_{\Lambda_{0}}}\mu_{\Lambda_{1}}\\
\vdots\\
b_{n}-\frac{b_{0}}{\mu_{\Lambda_{0}}}\mu_{\Lambda_{n}}
\end{array}\right)=T^{-1}\left(\begin{array}{c}
0\\
a_{1}\\
\vdots\\
a_{n}
\end{array}\right).
\]

Recall $T=(\Lambda_{i}.\Lambda_{j})_{i,j=1,\dots,n}$ is the $n\times n$ lower right submatrix of the intersection matrix of $X$. Write $T^{-1}=(c^{i,j})_{i,j=1,\dots,n}$ where $c^{i,j}$ is the $(i,j)$-entry of the inverse of $T$. Then we have 
\[
T^{-1}\left(\begin{array}{c}
0\\
a_{1}\\
\vdots\\
a_{n}
\end{array}\right)=\left(\begin{array}{c}
a_{1}c^{1,1}+a_{2}c^{1,2}+\cdots+a_{n}c^{1,n}\\
a_{1}c^{2,1}+a_{2}c^{2,2}+\cdots+a_{n}c^{2,n}\\
\vdots\\
a_{1}c^{n,1}+a_{2}c^{n,2}+\cdots+a_{n}c^{n,n}
\end{array}\right).
\]
Therefore
\[
b_{i}-\frac{b_{0}}{\mu_{\Lambda_{0}}}\mu_{\Lambda_{i}}=\sum_{j=1}^{n}a_{j}c^{i,j}.
\]
Let $f\in H^{0}(\mathfrak{X},{\cal L})$ be non-zero and let $\vec{b}=(\nu_{\Lambda_{0}}(f),\dots,\nu_{\Lambda_{n}}(f))$. The $j$th coordinate of $\vec{a}=M\vec{b}$ is precisely 
\[
a_{j}=\sum_{i=0}^{n}\nu_{\Lambda_{i}}(f)\Lambda_{i}.\Lambda_{j}=\deg({\cal L}|_{\Lambda_{j}})-H_{f}.\Lambda_{j}
\]
by equation (\ref{main eqn}). Thus we have 
\begin{equation}
\nu_{\Lambda_{i}}(f)-\frac{\nu_{\Lambda_{0}}(f)}{\mu_{\Lambda_{0}}}\mu_{\Lambda_{i}}=\sum_{j=1}^{n}\left(\deg({\cal L}|_{\Lambda_{j}})-H_{f}.\Lambda_{j}\right)c^{i,j}.\label{eq: gen inverted main eqn}
\end{equation}
Now we take $\mathfrak{X}=\mathfrak{X}(Np^{r})_{/\mathbb{Z}_{p}[\zeta_{Np^{r}}]}$ to be our modular curve and ${\cal L}=\underline{\omega}^{\otimes 2k}$ to be the modular sheaf. Recall $\bar{\mathfrak{X}}$ is reduced so $\mu_{\Lambda_{i}}=1$ for each $i$. Then equation (\ref{eq: gen inverted main eqn}) becomes 
\[
\nu_{\Lambda_{i}}(f)-\nu_{\Lambda_{0}}(f)=\sum_{j=1}^{n}(\deg(\underline{\omega}^{\otimes2k}|_{\Lambda_{j}})-H_{f}.\Lambda_{j})c^{i,j}
\]
as desired.
\end{proof}
\noindent Thus Theorem \ref{Equation for v(f)} expresses the quantities $\nu_{\Lambda_{i}}(f)$ in terms of $\deg(\underline{\omega}^{\otimes 2k}|_{\Lambda_{j}})$, $H_{f}.\Lambda_{j}$, and the entries of the inverse of $T$. 

\section{\label{sec:Intersection-Matrix}Intersection matrix}

Throughout this chapter, we will use the following notation:
\begin{itemize}
    \item $N\geq3$ and $r\geq1$ will be integers and $p\geq2$ will be a prime such that $p\nmid N$
    \item $\mathbb{F}_{q}$ will denote the residue field $\mathbb{F}_{p}(\zeta_{N})$ of $\mathbb{Z}_{p}[\zeta_{Np^{r}}]$ where $q=p^{\ell}$ where $\ell$ is the order of $p$ in $(\mathbb{Z}/N\mathbb{Z})^{\times}$. 
    \item $\mathfrak{X}$ will denote the modular curve $\mathfrak{X}(Np^{r})_{/\mathbb{Z}_{p}[\zeta_{Np^{r}}]}$.
\end{itemize}

\subsection{Intersection matrix for the modular curve}

In this section, we will obtain an explicit description of the intersection matrix of $\mathfrak{X}$. First we recall the description of the irreducible components of the special fiber $\bar{\mathfrak{X}}$ as in Theorem \ref{Special Fiber Comps}.

\begin{thm}
The special fiber of $\mathfrak{X}(Np^{r})$ is the disjoint union, with crossings at the supersingular points of $\mathfrak{X}(N)_{/\mathbb{F}_{q}}$, of the exotic Igusa curves ${\rm ExIg}(p^{r},r,N)$ over $\mathfrak{X}(N)_{/\mathbb{F}_{q}}$, indexed by 
\[
(\mathbb{Z}/p^{r}\mathbb{Z})^{\times}/{\rm HomSurj}((\mathbb{Z}/p^{r}\mathbb{Z})^{2},\mathbb{Z}/p^{r}\mathbb{Z}).
\]
Furthermore, $\bar{\mathfrak{X}}(Np^{r})$ is reduced.
\end{thm}

In the paragraph proceeding Corollary \ref{Irr Comps Same Num Int}, a complete list of representatives in
\[
(\mathbb{Z}/p^{r}\mathbb{Z})^{\times}/{\rm HomSurj}((\mathbb{Z}/p^{r}\mathbb{Z})^{2},\mathbb{Z}/p^{r}\mathbb{Z})
\]
was given by 
\[
\begin{cases}
\Lambda_{(1,-a)} & a\in\mathbb{Z}/p^{r}\mathbb{Z}\\
\Lambda_{(-pb,1)} & b\in\mathbb{Z}/p^{r-1}\mathbb{Z}
\end{cases}.
\]

We will identify an irreducible component of $\bar{\mathfrak{X}}$ by it's index $\Lambda$. By \cite[{}13.8.5]{[KM85]}, the local intersection number at a supersingular point $s$ between two distinct irreducible components $\Lambda_{1}$ and $\Lambda_{2}$ is precisely 
\begin{equation}
i_{s}(\Lambda_{1},\Lambda_{2})=\left[\#\left((\mathbb{Z}/p^{r}\mathbb{Z})^{2}/(\ker\Lambda_{1}+\ker\Lambda_{2})\right)\right]^{2}.\label{KM intersection formula}
\end{equation}

\begin{rmk}
In \cite[{}13.8.5]{[KM85]}, we require $s$ to be a $\mathbb{F}_{q}$-rational supersingular point. However, when $q=p$, the supersingular points may not be $\mathbb{F}_{q}$-rational. Indeed, \cite[p. 96]{[KM85]}, shows that the supersingular points are $\mathbb{F}_{p^{2}}$-rational. However, we can instead compute the intersection numbers by \'{e}tale base change. Consider the ring $R'=\mathbb{Z}_{p}[\zeta_{p^{r}},\zeta_{p^{2}-1}]$. By \cite[IV, §4, Prop 16 \& 17]{[Ser79]}, $R'$ is the ring of integers of the local field $\mathbb{Q}_{p}(\zeta_{p^{r}},\zeta_{p^{2}-1})$ so $R'$ is a DVR. Since $p$ has order $2$ in $(\mathbb{Z}/(p^{2}-1)\mathbb{Z})^{\times}$, the residue field of $R'$ is $\mathbb{F}_{p^{2}}$. The map $\mathbb{Z}_{p}[\zeta_{p^{r}}]\rightarrow R'$ is unramified by \cite[IV, §4, Prop. 16]{[Ser79]}, noting that $p\nmid(p^{2}-1)$, and flat because $R'$ is torsion-free over the DVR $\mathbb{Z}_{p}[\zeta_{p^{r}}]$ (see \cite[Tag0539]{[Stacks]}). Therefore ${\rm Spec}(R')\rightarrow{\rm Spec}(\mathbb{Z}_{p}[\zeta_{p^{r}}])$ is \'{e}tale.

Let $\Lambda$ and $\Lambda'$ be two irreducible components of $\bar{\mathfrak{X}}$. Let $\mathfrak{X}'=\mathfrak{X}_{/R'}$ and let $s'$ be a supersingular point of $\mathfrak{X}'_{/\mathbb{F}_{p^{2}}}$ which maps to $s$. By \cite[{}9.1.5, {}9.1.6]{[Liu02]}, we have 
\[
i_{s'}(f^{*}\Lambda,f^{*}\Lambda')=i_{s}(\Lambda,\Lambda').
\]
Note that \cite[{}9.2.15]{[Liu02]} uses a desingularization of $\mathfrak{X}'$ \emph{i.e. }a proper birational morphism $\mathfrak{X}''\rightarrow\mathfrak{X}'$ where $\mathfrak{X}''$ is regular. Since $R$ and $R'$ are both DVRs, they are both regular, excellent, and noetherian. Therefore $\mathfrak{X}'$ is an arithmetic surface (by Proposition \ref{Base Change Modular}), so we can take $\mathfrak{X}''=\mathfrak{X}'$.
\end{rmk}

From equation (\ref{KM intersection formula}), the local intersection number doesn't depend on the supersingular point. Let ${\rm S}(N)$ (resp. $S(Np^{r})$) denote the supersingular locus of $\mathfrak{X}(N)_{/\mathbb{F}_{q}}$ (resp. $\bar{\mathfrak{X}}$). By \cite[{}12.7.2]{[KM85]}, the map ${\rm Ig}(p^{r},N)\rightarrow\mathfrak{X}(N)_{/\mathbb{F}_{q}}$ is totally ramified at the supersingular points so $\deg{\rm S}(Np^{r})=\deg{\rm S}(N)$. Therefore to obtain the global intersection number $\Lambda_{1}.\Lambda_{2}$ we multiply $i_{s}(\Lambda_{1},\Lambda_{2})$ by $\deg{\rm S}(N)$. We will now compute the local intersection number between each pair of irreducible components of $\bar{\mathfrak{X}}$, beginning with the case of distinct pairs.

\begin{prop}\label{int nums for distinct comps}
Let $\nu_{p}$ denote the $p$-adic valuation normalized so $\nu_{p}(p)=1$ and let $s$ be a supersingular point. We have
\begin{align*}
i_{s}(\Lambda_{(1,-a)},\Lambda_{(-pb,1)}) & = 1\\
i_{s}(\Lambda_{(1,-a)},\Lambda_{(1,-a')}) & = p^{2\nu_{p}(a'-a)}\\
i_{s}(\Lambda_{(-pb,1)},\Lambda_{(-pb',1)}) & = p^{2\nu_{p}(b'-b)+2}
\end{align*}
for $a,a'\in\mathbb{Z}/p^{r}\mathbb{Z}$ distinct and $b,b'\in\mathbb{Z}/p^{r-1}\mathbb{Z}$ distinct. 
\end{prop}

\begin{proof}
Since each group homomorphism $\Lambda:(\mathbb{Z}/p^{r}\mathbb{Z})^{2}\rightarrow\mathbb{Z}/p^{r}\mathbb{Z}$ corresponding to an irreducible component is a surjective, we know $\#\ker\Lambda=p^{r}$. By definition, 
\[
\Lambda_{(1,-a)}(a,1)=a\Lambda_{(1,-a)}(1,0)+\Lambda_{(1,-a)}(0,1)=a-a=0.
\]
Therefore 
\[
\ker\Lambda_{(1,-a)}={\rm span}_{\mathbb{Z}/p^{r}\mathbb{Z}}\left\{ \left(\begin{array}{c} a\\
1
\end{array}\right)\right\} .
\]
Similarly, since 
\[
\Lambda_{(-pb,1)}(1,pb)=\Lambda_{(-pb,1)}(1,0)+pb\Lambda_{(-pb,1)}(0,1)=-pb+pb=0
\]
we have 
\[
\ker\Lambda_{(-pb,1)}={\rm span}_{\mathbb{Z}/p^{r}\mathbb{Z}}\left\{ \left(\begin{array}{c}
1\\
pb
\end{array}\right)\right\} .
\]
By equation (\ref{KM intersection formula}), we have 
\[
i_{s}(\Lambda_{1},\Lambda_{2})=(p^{2r}/\#(\ker\Lambda_{1}+\ker\Lambda_{2}))^{2}.
\]
We will now compute $\#(\ker\Lambda_{1}+\ker\Lambda_{2})$ by considering the following three cases.

\textbf{Case 1}:
We have 
\[
\ker\Lambda_{(1,-a)}+\ker\Lambda_{(-pb,1)}={\rm span}_{\mathbb{Z}/p^{r}\mathbb{Z}}\left\{ \left(\begin{array}{c}
a\\
1
\end{array}\right),\left(\begin{array}{c}
1\\
pb
\end{array}\right)\right\} .
\]
Note that 
\[
\det\left(\begin{array}{cc}
a & 1\\
1 & pb
\end{array}\right)=pab-1
\]
is invertible in $\mathbb{Z}/p^{r}\mathbb{Z}$. Therefore $\left\{ (a,1),(1,pb)\right\}$ is a basis for $(\mathbb{Z}/p^{r}\mathbb{Z})^{2}$ hence
\[
\ker\Lambda_{(1,-a)}+\ker\Lambda_{(-pb,1)}=(\mathbb{Z}/p^{r}\mathbb{Z})^{2}.
\]
We conclude $i_{s}(\Lambda_{(1,-a)}.\Lambda_{(-pb,1)})=1$.

\textbf{Case 2}: We have
\begin{align*}
\ker\Lambda_{(1,-a)}+\ker\Lambda_{(1,-a')} & = {\rm span}_{\mathbb{Z}/p^{r}\mathbb{Z}}\left\{ \left(\begin{array}{c}
a\\
1
\end{array}\right),\left(\begin{array}{c}
a'\\
1
\end{array}\right)\right\} \\
 & = {\rm span}_{\mathbb{Z}/p^{r}\mathbb{Z}}\left\{ \left(\begin{array}{c}
a\\
1
\end{array}\right),\left(\begin{array}{c}
a'-a\\
0
\end{array}\right)\right\} 
\end{align*}
For any $0\leq i\leq r$, consider the injective group homomorphism
\[
m_{i}:\mathbb{Z}/p^{r-i}\mathbb{Z}\rightarrow\mathbb{Z}/p^{r}\mathbb{Z}
\]
given by multiplication by $p^{i}$. This gives us a filtration 
\[
\mathbb{Z}/p^{r}\mathbb{Z}={\rm im}(m_{0})\supset{\rm im}(m_{1})\supset\cdots\supset{\rm im}(m_{r-1})\supset{\rm im}(m_{r})=(0)
\]
which is exhaustive and separated. Thus for any non-zero $c\in\mathbb{Z}/p^{r}\mathbb{Z}$, there exists a unique smallest $i$ such that $c\in{\rm im}(m_{i})$. We let $\nu_{p}(c)$ denote the quantity $i$. We claim $c=p^{\nu_{p}(c)}u_{c}$ for some $u_{c}\in(\mathbb{Z}/p^{r}\mathbb{Z})^{\times}$. Indeed, if $u_{c}$ is not a unit of $(\mathbb{Z}/p^{r}\mathbb{Z})^{\times}$, then $u_{c}\in(p)$ so $u_{c}=pu_{c}'$ for some $u_{c}'\in\mathbb{Z}/p^{r}\mathbb{Z}$. Hence $c=p^{\nu_{p}(c)+1}u_{c}'\in{\rm im}(m_{\nu_{p}(c)+1})$, contradicting the minimality of $\nu_{p}(c)$. 

We will show the sequence of abelian groups 
\begin{equation}
\begin{gathered}\xymatrix{0\ar@{->}[r] & \mathbb{Z}/p^{\nu_{p}(c)}\mathbb{Z}\ar@{->}[rr]^{m_{r-\nu_{p}(c)}} &  & \mathbb{Z}/p^{r}\mathbb{Z}\ar@{->}[rr]^{\cdot c} &  & {\rm span}_{\mathbb{Z}/p^{r}\mathbb{Z}}(c)\ar@{->}[r] & 0}
\end{gathered}
\label{eq:exact of mi}
\end{equation}
is exact. Exactness at the second and fourth term are clear. Let $x\in\mathbb{Z}/p^{\nu_{p}(c)}\mathbb{Z}$. Then 
\[
c\cdot m_{r-\nu_{p}(c)}(x)=c\cdot p^{r-\nu_{p}(c)}x=u_{c}p^{\nu_{p}(c)}p^{r-\nu_{p}(c)}x=u_{c}p^{r}x=0
\]
so ${\rm im}(m_{r-\nu_{p}(c)})\subseteq\ker(\cdot c)$. On the other hand, let $d\in\ker(\cdot c)$ so $d\cdot c=0$. Then 
\[
u_{c}u_{d}p^{\nu_{p}(c)}p^{\nu_{p}(d)}=0.
\]
Since $u_{c},u_{d}\in(\mathbb{Z}/p^{r}\mathbb{Z})^{\times}$, we have $p^{\nu_{p}(c)+\nu_{p}(d)}=0$. Hence $\nu_{p}(d)\geq r-\nu_{p}(c)$ or equivalently ${\rm im}(m_{d})\subseteq{\rm im}(m_{r-\nu_{p}(c)})$. Therefore $\ker(\cdot c)\subseteq{\rm im}(m_{r-\nu_{p}(c)})$ allowing us to conclude exactness at the third term.

Let $i=\nu_{p}(a'-a)$. Consider the sequence of abelian groups
\begin{equation}
0\rightarrow\mathbb{Z}/p^{i}\mathbb{Z}\overset{\phi}{\longrightarrow}\mathbb{Z}/p^{r}\mathbb{Z}\times\mathbb{Z}/p^{r}\mathbb{Z}\overset{\psi}{\longrightarrow}{\rm span}_{\mathbb{Z}/p^{r}\mathbb{Z}}\left\{ \left(\begin{array}{c}
a\\
1
\end{array}\right),\left(\begin{array}{c}
a'-a\\
0
\end{array}\right)\right\} \rightarrow0\label{eq:exact of phi psi}
\end{equation}
where $\phi(x)=(0,m_{r-i}(x))$ and 
\[
\psi(c,d)=c\left(\begin{array}{c}
a\\
1
\end{array}\right)+d\left(\begin{array}{c}
a'-a\\
0
\end{array}\right)=\left(\begin{array}{c}
ca+d(a'-a)\\
c
\end{array}\right).
\]

One can show exactness similar to the above sequence (\ref{eq:exact of mi}). Indeed, showing ${\rm im}(\phi)\subseteq\ker(\psi)$ is clear. Conversely, if $\psi(c,d)=(0,0)$, then $c=0$ and consequently $d(a'-a)=0$. Writing $p^{\nu_{p}(d)}p^{i}u_{d}u_{a'-a}=0$, we similarly conclude $\nu_{p}(d)\geq r-i$ so $\ker(\psi)\subset{\rm im}(\phi)$. Using exactness of (\ref{eq:exact of phi psi}), we get
\begin{align*}
\#{\rm span}_{\mathbb{Z}/p^{r}\mathbb{Z}}\left\{ \left(\begin{array}{c}
a\\
1
\end{array}\right),\left(\begin{array}{c}
a'-a\\
0
\end{array}\right)\right\} & = \frac{\#(\mathbb{Z}/p^{r}\mathbb{Z}\times\mathbb{Z}/p^{r}\mathbb{Z})}{\#\mathbb{Z}/p^{i}\mathbb{Z}}\\
& = p^{2r-i}.
\end{align*}
We conclude $i_{s}(\Lambda_{(1,-a)}.\Lambda_{(1,-a')})=p^{2\nu_{p}(a'-a)}$.

\textbf{Case 3}:  Lastly, we have
\begin{align*}
\ker\Lambda_{(-pb,1)}+\ker\Lambda_{(-pb',1)} & = {\rm span}_{\mathbb{Z}/p^{r}\mathbb{Z}}\left\{ \left(\begin{array}{c}
1\\
pb
\end{array}\right),\left(\begin{array}{c}
1\\
pb'
\end{array}\right)\right\} \\
 & = {\rm span}_{\mathbb{Z}/p^{r}\mathbb{Z}}\left\{ \left(\begin{array}{c}
1\\
pb
\end{array}\right),\left(\begin{array}{c}
0\\
p(b'-b)
\end{array}\right)\right\} 
\end{align*}
Arguing as above by creating a sequence similar to (\ref{eq:exact of phi psi}) and replacing $a$ (resp. $a'$) with $pb$ (resp. $pb'$), we get
\[
\#\left(\ker\Lambda_{(-pb,1)}+\ker\Lambda_{(-pb',1)}\right)=p^{r}\cdot p^{r-(\nu_{p}(b'-b)-1)}.
\]
We conclude 
\[
i_{s}(\Lambda_{(-pb,1)}.\Lambda_{(-pb',1)})=(p^{2r}/p^{2r-\nu_{p}(b'-b)+1})^{2}=p^{2\nu_{p}(b'-b)+2}.\qedhere
\].
\end{proof}

To finish our calculation of intersection numbers, we will now compute each self-intersection. First we introduce the following lemma.

\begin{lem}\label{sum of p^(vp)}
Let $r\geq1$. We have 
\[
\sum_{\substack{a'\in\mathbb{Z}/p^{r}\mathbb{Z}\\
a'\ne0
}
}p^{2\nu_{p}(a')}=p^{2r-1}-p^{r-1}.
\]
\end{lem}
\begin{proof}
We will group the index $a'\in\mathbb{Z}/p^{r}\mathbb{Z}$ based on its $p$-adic valuation $m=\nu_{p}(a')$ and then sum over $m$. Observe that there are precisely $\varphi(p^{r-m})$-many elements in $\mathbb{Z}/p^{r}\mathbb{Z}$ with $p$-adic valuation equal to $m$. Thus
\begin{align*}
\sum_{\substack{a'\in\mathbb{Z}/p^{r}\mathbb{Z}\\
a'\ne0
}
}p^{2\nu_{p}(a')} & = \sum_{m=0}^{r-1}\varphi(p^{r-m})p^{2m}\\
 & = \sum_{m=0}^{r-1}p^{r-m-1}(p-1)p^{2m}\\
 & = \sum_{m=0}^{r-1}p^{r+m-1}(p-1)\\
 & = (p-1)p^{r-1}\sum_{m=0}^{r-1}p^{m}\\
 & = (p-1)p^{r-1}\frac{1-p^{r}}{1-p}\\
 & = p^{2r-1}-p^{r-1}
\end{align*}
as desired. 
\end{proof}

\begin{prop}\label{self-int num}
For any $\Lambda\in(\mathbb{Z}/p^{r}\mathbb{Z})^{\times}\backslash{\rm HomSurj}((\mathbb{Z}/p^{r}\mathbb{Z})^{2},\mathbb{Z}/p^{r}\mathbb{Z})$, the self-intersection number is 
\[
\Lambda.\Lambda=-\deg{\rm S}(N)\cdot p^{2r-1}.
\]
\end{prop}

\begin{proof}
By \cite[Proposition 9.1.21]{[Liu02]}, we have
\[
\Lambda.\Lambda=-\frac{1}{\mu_{\Lambda}}\sum_{\Lambda'\ne\Lambda}\mu_{\Lambda'}(\Lambda'.\Lambda)
\]
where $\mu_{\Lambda'}$ is the multiplicity of $\Lambda'$ for any irreducible component $\Lambda'$. Since $\bar{\mathfrak{X}}$ is reduced, $\mu_{\Lambda}=1$. 

We compute the self-intersection for the two possible cases of $\Lambda$. First we have
\begin{align*}
\Lambda_{(1,-a)}.\Lambda_{(1,-a)} & = -\sum_{\substack{a'\in\mathbb{Z}/p^{r}\mathbb{Z}\\
a'\ne a
}
}\Lambda_{(1,-a')}.\Lambda_{(1,-a)}-\sum_{b\in\mathbb{Z}/p^{r-1}\mathbb{Z}}\Lambda_{(-pb,1)}.\Lambda_{(1,-a)}\\
 & = -\deg{\rm S}(N)\sum_{\substack{a'\in\mathbb{Z}/p^{r}\mathbb{Z}\\
a'\ne a
}
}p^{2\nu_{p}(a'-a)}-\deg{\rm S}(N)\sum_{b\in\mathbb{Z}/p^{r-1}\mathbb{Z}}1\\
 & = -\deg{\rm S}(N)\sum_{\substack{a'\in\mathbb{Z}/p^{r}\mathbb{Z}\\
a'\ne0
}
}p^{2\nu_{p}(a')}-\deg{\rm S}(N)\cdot p^{r-1}
\end{align*}
where we have used Proposition \ref{int nums for distinct comps} to calculate the intersection numbers. Using Lemma \ref{sum of p^(vp)}, we have 
\begin{align*}
\Lambda_{(1,-a)}.\Lambda_{(1,-a)} & = -\deg{\rm S}(N)(p^{2r-1}-p^{r-1})-\deg{\rm S}(N)\cdot p^{r-1}\\
 & = -\deg{\rm S}(N)\cdot p^{2r-1}.
\end{align*}
Next we consider the case $\Lambda=\Lambda_{(-pb,1)}$. We have
\begin{alignat*}{1}
\Lambda_{(-pb,1)}.\Lambda_{(-pb,1)} & =-\sum_{\substack{b'\in\mathbb{Z}/p^{r-1}\mathbb{Z}\\
b'\ne b
}
}\Lambda_{(-pb',1)}.\Lambda_{(-pb,1)}-\sum_{a\in\mathbb{Z}/p^{r}\mathbb{Z}}\Lambda_{(1,-a)}.\Lambda_{(1,-pb)}\\
 & =-\deg{\rm S}(N)\sum_{\substack{b'\in\mathbb{Z}/p^{r-1}\mathbb{Z}\\
b'\ne b
}
}p^{2\nu_{p}(b'-b)+2}-\deg{\rm S}(N)\sum_{a\in\mathbb{Z}/p^{r}\mathbb{Z}}1\\
 & =-\deg{\rm S}(N)\sum_{\substack{b'\in\mathbb{Z}/p^{r-1}\mathbb{Z}\\
b'\ne0
}
}p^{2\nu_{p}(b')+2}-\deg{\rm S}(N)\cdot p^{r}\\
 & =-\deg{\rm S}(N)\cdot p^{2}\sum_{\substack{b'\in\mathbb{Z}/p^{r-1}\mathbb{Z}\\
b'\ne0
}
}p^{2\nu_{p}(b')}-\deg{\rm S}(N)\cdot p^{r}\\
 & =-\deg{\rm S}(N)\cdot p^{2}\left(p^{2(r-1)-1}-p^{r-2}\right)-\deg{\rm S}(N)\cdot p^{r}\\
 & =-\deg{\rm S}(N)\cdot p^{2r-1}.\tag*{\qedhere}
\end{alignat*}
\end{proof}

We will now describe the intersection matrix $M$ by specifying four blocks which make up $M$. We label the columns (and by symmetry the rows) of $M$ in the following order:{\footnotesize{}
\[
\Lambda_{(1,0)}, \Lambda_{(1,-1)}, \dots, \Lambda_{(1,-a)}, \dots, \Lambda_{(1,-(p^{r}-1))},  \Lambda_{(0,1)}, \Lambda_{(-p,1)}, \dots, \Lambda_{(-pb,1)},\dots, \Lambda_{(-p(p^{r-1}-1),1)}
\]
}so that the $(i,j)$ entry of $M$ is equal to the intersection number between the $i$th row label and $j$th column label. Since $\deg{\rm S}(N)$ is a common factor among each entry of $M$, we will describe the matrix $\frac{1}{\deg{\rm S}(N)}M$ to simplify exposition.

Let $M_{11}$ (resp. $M_{22}$) denote the submatrix of $\frac{1}{\deg{\rm S}(N)}M$ corresponding to the column and row labels of the form $\Lambda_{(1,-a)}$ (resp. $\Lambda_{(-pb,1)}$). We let $M_{12}$ and $M_{21}$ denote the remaining two submatrices of $\frac{1}{\deg{\rm S}(N)}M$ so that
\[
M=\deg{\rm S}(N)\left(\begin{array}{cc}
M_{11} & M_{12}\\
M_{21} & M_{22}
\end{array}\right).
\]
We also let $M(p^{r})$ denote the matrix $M_{11}$ to highlight the dependence on $p^{r}$. By convention, we define $M(p^{0})$ to be the $1\times1$ matrix consisting of the entry $-\frac{1}{p}$. The matrices $M_{11}$ and $M_{22}$ take on a special form.

\begin{defn}
An $n\times n$ \textbf{circulant} matrix $C$ is of the form
\[
C=\left(\begin{array}{ccccc}
c_{0} & c_{n-1} & \cdots & c_{2} & c_{1}\\
c_{1} & c_{0} & \cdots & c_{3} & c_{2}\\
\vdots & \vdots & \ddots & \vdots & \vdots\\
c_{n-2} & c_{n-3} & \cdots & c_{0} & c_{n-1}\\
c_{n-1} & c_{n-2} & \cdots & c_{1} & c_{0}
\end{array}\right)
\]
where each column is equal to the previous column shifted downward by 1, looping around as appropriate.
\end{defn}

Refer to Appendix \ref{apx: circulant matrices} for a discussion on circulant matrices, including an explicit description of the eigenvalues, eigenvectors, and inverse in terms of the entries of the matrix and roots of unity.

\begin{prop}\label{Description of M}\leavevmode
\begin{enumerate}
    \item[a.] The entries of $M_{12}$ and $M_{21}$ are all equal to 1.  
    \item[b.] The matrix $M_{11}=M(p^{r})$ is a $p^{r}\times p^{r}$ circulant matrix whose first entry of the first column is $-p^{2r-1}$. For $2\leq a\leq p^{r}$, the $a^{{\rm th}}$ entry in the first column is equal to $p^{2\nu_{p}(a-1)}$.
    \item[c.] The matrix $M_{22}$ is a $p^{r-1}\times p^{r-1}$ circulant matrix equal to $p^{2}M(p^{r-1})$.
\end{enumerate}
\end{prop}
\begin{proof}
Proposition \ref{int nums for distinct comps} immediately tells us the entries of $M_{12}$ and $M_{21}$ are all 1. Next we describe the first column of $M_{11}$. The first entry is equal to the local self-intersection number of $\Lambda_{(1,0)}$, which is $-p^{2r-1}$ by Proposition \ref{self-int num}. Using Proposition \ref{int nums for distinct comps}, for $2\leq a\leq p^{r}$, the $a^{{\rm th}}$ entry of the first column is equal to the local intersection number 
\[
\frac{1}{\deg{\rm S}(N)}\Lambda_{(1,-(a-1))}.\Lambda_{(1,0)}=p^{2\nu_{p}(a-1)}.
\]

Now we show $M_{11}$ is circulant. Recall the $j$th column of $M_{11}$ corresponds to the label $\Lambda_{(1,-(j-1))}$ for $1\leq j\leq p^{r}$. The $k$th entry in the $j$th column is equal to 
\[
\frac{1}{\deg{\rm S}(N)}\Lambda_{(1,-(k-1))}.\Lambda_{(1,-(j-1))}=\begin{cases}
p^{2\nu_{p}((k-1)-(j-1))} & \mbox{if }j\ne k\\
-p^{2r-1} & \mbox{if }j=k
\end{cases}.
\]
Note for any nonzero $x,y\in\mathbb{Z}$, if $x\equiv y$ (mod $p^{r}$), then $\nu_{p}(x)=\nu_{p}(y)$. Therefore the quantity 
\[
\nu_{p}((k-1)-(j-1))
\]
for $k\ne j$ remains unchanged if we take $(k-1)-(j-1)$ modulo $p^{r}$. We conclude 
\[
\Lambda_{(1,-(k-1))}.\Lambda_{(1,-(j-1))}=\Lambda_{(1,0)}.\Lambda_{(1,-(k-1)+(j-1))}
\]
which says the $j$th column is equal to the first column with every entry shifted downward by $j-1$, looping around as appropriate. Hence $M_{11}$ is a circulant matrix.

Lastly, we show $M_{22}=p^{2}M(p^{r-1})$. Suppose $r>1$. By Proposition \ref{int nums for distinct comps}, the $(i,j)$ entry of $M_{22}$ for $i\ne j$ is equal to
\[
p^{2\nu_{p}((j-1)-(i-1))+2}=p^{2}p^{2\nu_{p}((j-1)-(i-1))}
\]
which is equal to the $(i,j)$ entry of $M(p^{r-1})$ multiplied by $p^{2}$. When $i=j$, the $(i,j)$ entry of $M_{22}$ is $-p^{2r-1}$ while the $(i,j)$ entry of $p^{2}M(p^{r-1})$ is 
\[
-p^{2}p^{2(r-1)-1}=-p^{2}p^{2r-3}=-p^{2r-1}.
\]
Thus $M_{22}=p^{2}M(p^{r-1})$. 

When $r=1$, $M_{22}$ is a $1\times1$ matrix consisting of the entry $(-p)$. By convention $M(p^{0})=(-\frac{1}{p})$ so $M_{22}=p^{2}M(p^{r-1})$ in the case $r=1$.
\end{proof}

\begin{ex*}
The intersection matrix for $\mathfrak{X}(N\cdot5)$ is the $6\times6$ matrix 
\[
\deg{\rm S}(N)\left(\begin{array}{cccccc}
\textbf{-5} & \textbf{1} & \textbf{1} & \textbf{1} & \textbf{1} & 1\\
\textbf{1} & \textbf{-5} & \textbf{1} & \textbf{1} & \textbf{1} & 1\\
\textbf{1} & \textbf{1} & \textbf{-5} & \textbf{1} & \textbf{1} & 1\\
\textbf{1} & \textbf{1} & \textbf{1} & \textbf{-5} & \textbf{1} & 1\\
\textbf{1} & \textbf{1} & \textbf{1} & \textbf{1} & \textbf{-5} & 1\\
1 & 1 & 1 & 1 & 1 & \textit{-5}
\end{array}\right)
\]
where the entries in $\textbf{bold}$ comprise $M_{11}=M(p^{r})$ and the entries in \textit{italic} comprise $M_{22}$. The intersection matrix for $\mathfrak{X}(N3^{2})$ is the $12\times12$ matrix 

\[
\deg{\rm S}(N)\left(\begin{array}{cccccccccccc}
\mathbf{-3^{3}} & \textbf{1} & \textbf{1} & \mathbf{3^{2}} & \textbf{1} & \textbf{1} & \mathbf{3^{2}} & \textbf{1} & \textbf{1} & 1 & 1 & 1\\
\textbf{1} & \mathbf{-3^{3}} & \textbf{1} & \textbf{1} & \mathbf{3^{2}} & \textbf{1} & \textbf{1} & \mathbf{3^{2}} & \textbf{1} & 1 & 1 & 1\\
\textbf{1} & \textbf{1} & \textbf{$-3^{3}$} & \textbf{1} & \textbf{1} & \mathbf{3^{2}} & \textbf{1} & \textbf{1} & \mathbf{3^{2}} & 1 & 1 & 1\\
\mathbf{3^{2}} & \textbf{1} & \textbf{1} & \mathbf{-3^{3}} & \textbf{1} & \textbf{1} & \mathbf{3^{2}} & \textbf{1} & \textbf{1} & 1 & 1 & 1\\
\textbf{1} & \mathbf{3^{2}} & \textbf{1} & \textbf{1} & \mathbf{-3^{3}} & \textbf{1} & \textbf{1} & \mathbf{3^{2}} & \textbf{1} & 1 & 1 & 1\\
\textbf{1} & \textbf{1} & \mathbf{3^{2}} & \textbf{1} & \textbf{1} & \mathbf{-3^{3}} & \textbf{1} & \textbf{1} & \mathbf{3^{2}} & 1 & 1 & 1\\
\mathbf{3^{2}} & \textbf{1} & \textbf{1} & \mathbf{3^{2}} & \textbf{1} & \textbf{1} & \mathbf{-3^{3}} & \textbf{1} & \textbf{1} & 1 & 1 & 1\\
\textbf{1} & \mathbf{3^{2}} & \textbf{1} & \textbf{1} & \mathbf{3^{2}} & \textbf{1} & \textbf{1} & \mathbf{-3^{3}} & \textbf{1} & 1 & 1 & 1\\
\textbf{1} & \textbf{1} & \mathbf{3^{2}} & \textbf{1} & \textbf{1} & \mathbf{3^{2}} & \textbf{1} & \textbf{1} & \mathbf{-3^{3}} & 1 & 1 & 1\\
1 & 1 & 1 & 1 & 1 & 1 & 1 & 1 & 1 & \mathit{-3^3} & \mathit{3^2} & \mathit{3^2}\\
1 & 1 & 1 & 1 & 1 & 1 & 1 & 1 & 1 & \mathit{3^{2}} & \mathit{-3^{3}} & \mathit{3^{2}}\\
1 & 1 & 1 & 1 & 1 & 1 & 1 & 1 & 1 & \mathit{3^{2}} & \mathit{3^{2}} & \mathit{-3^{3}}
\end{array}\right)
\]
\end{ex*}

\subsection{\label{sub:Inverting M(p^r)} \texorpdfstring{Inverting $M(p^{r})$}{}}

Recall our goal is to invert the matrix $T$ obtained by removing the first row and column of $M$. For a general matrix $A$, let $A_{\hat{1,}\hat{1}}$ denote the matrix obtained by removing the first row and column of $A$ and let $\mathbf{1}_{n\times m}$ denote the $n\times m$ matrix whose entries are all 1. Using Proposition \ref{Description of M}, we have the following description of $T$: 
\[
T=\deg{\rm S}(N)\left(\begin{array}{cc}
M(p^{r})_{\hat{1},\hat{1}} & \mathbf{1}_{p^{r}-1\times p^{r-1}}\\
\mathbf{1}_{p^{r-1}\times p^{r}-1} & p^{2}M(p^{r-1})
\end{array}\right).
\]
We will use the following identity to invert $T$.

\begin{prop}[Woodbury Matrix Identity]\label{WMI}
Let $A$ be an $n\times n$ invertible matrix, $C$ an invertible $k\times k$ matrix where $k\leq n$, $U$ an $n\times k$ matrix, and $V$ a $k\times n$ matrix. Then
\[
(A+UCV)^{-1}=A^{-1}-A^{-1}U(C^{-1}+VA^{-1}U)^{-1}VA^{-1}.
\]
\end{prop}
Refer to Appendix \ref{apx: Inv via WMI} for a discussion on using the Woodbury Matrix Identity to compute the inverse of $A+UCV$ in the situation that $A$ is a block diagonal matrix with two blocks, $C$ is the $2\times2$ identity, and both $U$ and $V$ consists of 0's and 1's which we specify later. Note that the inverse of a block diagonal matrix is obtained by inverting each block. Thus computing $A^{-1}$ amounts to computing $M(p^{r})^{-1}$ and $M(p^{r})_{\hat{1},\hat{1}}^{-1}$ which we will do in this subsection.

Before we compute the eigenvalues of $M(p^{r})$, we will need the following technical lemma.

\begin{lem}\label{Sum of RoU Coprime to p}
Let $N\geq1$ and $J\geq1$ be integers. We have 
\[
\sum_{\substack{u=1\\
p\nmid u
}
}^{p^{N}-1}\zeta_{p^{N}}^{-uJ}=\begin{cases}
0 & \mbox{if }p^{N-1}\nmid J\\
-p^{N-1} & \mbox{if }p^{N-1}\mid J\mbox{ and }p^{N}\nmid J\\
p^{N}-p^{N-1} & \mbox{if }p^{N}\mid J
\end{cases}.
\]
\end{lem}
\begin{proof}
By the geometric partial sum formula, we have

\begin{align*}
\sum_{u=1}^{p^{N}-1}\zeta_{p^{N}}^{-uJ} & = \begin{cases}
-1+\frac{1-\zeta_{p^{N}}^{-Jp^{N}}}{1-\zeta_{p^{N}}^{-J}} & \mbox{if }p^{N}\nmid J\\
p^{N}-1 & \mbox{if }p^{N}\mid J
\end{cases}\\
 & = \begin{cases}
-1 & \mbox{if }p^{N}\nmid J\\
p^{N}-1 & \mbox{if }p^{N}\mid J
\end{cases}.
\end{align*}
Next we consider the sum 
\[
\sum_{\substack{u=1\\
p\mid u
}
}^{p^{N}-1}\zeta_{p^{N}}^{-uJ}=\sum_{u=1}^{p^{N-1}-1}\zeta_{p^{N-1}}^{-uJ}=\begin{cases}
-1 & \mbox{if }p^{N-1}\nmid J\\
p^{N-1}-1 & \mbox{if }p^{N-1}\mid J
\end{cases}.
\]
Therefore
\[
\sum_{\substack{u=1\\
p\nmid u
}
}^{p^{N}-1}\zeta_{p^{N}}^{-uJ}=\sum_{u=1}^{p^{N}-1}\zeta_{p^{N}}^{-uJ}-\sum_{\substack{u=1\\
p\mid u
}
}^{p^{N}-1}\zeta_{p^{N}}^{-uJ}=\begin{cases}
0 & \mbox{if }p^{N-1}\nmid J\\
-p^{N-1} & \mbox{if }p^{N-1}\mid J\mbox{ and }p^{N}\nmid J\\
p^{N}-p^{N-1} & \mbox{if }p^{N}\mid J
\end{cases}.\qedhere
\]
\end{proof}

\begin{lem}\label{Eigenvalues of M(p^r)}
The eigenvalues $\lambda_{j}$ of $M(p^{r})$ are $\lambda_{1}=-p^{r-1}$ and 
\[
\lambda_{j}=-p^{2r-2-\nu_{p}(j-1)}(p+1)
\]
for $2\leq j\leq p^{r}$.
\end{lem}
\begin{proof}
By Lemma \ref{Evals of Circulant}, the eigenvalues of an $n\times n$ circulant matrix whose first column has entries $c_{0},c_{1},\dots,c_{n-1}$ are given by 
\[
\lambda_{j}=\sum_{k=0}^{n-1}c_{k}\zeta_{n}^{j(n-k)}.
\]
For $M(p^{r})$, the eigenvalues are therefore given by 
\[
\lambda_{j}=-p^{2r-1}+\sum_{k=1}^{p^{r}-1}p^{2\nu_{p}(k)}\zeta_{p^{r}}^{(j-1)(p^{r}-k)}=-p^{2r-1}+\sum_{k=1}^{p^{r}-1}p^{2\nu_{p}(k)}\zeta_{p^{r}}^{-k(j-1)}
\]
for $j=1,\dots,p^{r}$. 

When $j=1$, we have 
\[
\lambda_{1}=-p^{2r-1}+\sum_{k=1}^{p^{r}-1}p^{2\nu_{p}(k)}=-p^{2r-1}+(p^{2r-1}-p^{r-1})=-p^{r-1}
\]
where we have used Lemma \ref{sum of p^(vp)} to calculate the sum.

Assume $j>1$. When $r=1$, we have 
\[
\lambda_{j}=-p+\sum_{k=1}^{p-1}p^{2\nu_{p}(k)}\zeta_{p}^{-k(j-1)}.
\]
Since $\nu_{p}(k)=0$ for $1\leq k\leq p-1$, we have \[
\lambda_{j}=-p+\sum_{k=1}^{p-1}\left(\zeta_{p}^{-(j-1)}\right)^{k}
\]
As $2\leq j\leq p$, we have $p\nmid(j-1)$ so $\zeta_{p}^{-(j-1)}\ne1$. Using the geometric series partial sum formula, we get 
\[
\lambda_{j}=-p-1+\frac{1-\zeta_{p}^{-(j-1)p}}{1-\zeta_{p}^{-(j-1)}}=-p-1.
\]

Lastly we handle the $r>1$ case. We compute the sum appearing in the expression for $\lambda_{j}$ by breaking it up according to the value of $\nu_{p}(k)$. We have
\[
\sum_{k=1}^{p^{r}-1}p^{2\nu_{p}(k)}\zeta_{p^{r}}^{-k(j-1)}=\sum_{m=0}^{r-1}\left(p^{2m}\sum_{\substack{k=1\\
\nu_{p}(k)=m
}
}^{p^{r}-1}\zeta_{p^{r}}^{-k(j-1)}\right)
\]
We can rewrite the index in each sum as $k=up^{\nu_{p}(k)}$ where $p\nmid u$ and $1\leq u\leq p^{r-\nu_{p}(k)}-1$. Re-indexing, with $u\geq1$, we get 
\[
=\sum_{m=0}^{r-1}\left(p^{2m}\sum_{\substack{u=1\\
p\nmid u
}
}^{p^{r-m}-1}\zeta_{p^{r}}^{-up^{m}(j-1)}\right)
\]
Since $\zeta_{p^{r}}^{p^{m}}=\zeta_{p^{r-m}}$ for $0\leq m<r$, we have
\begin{equation}
=\sum_{m=0}^{r-1}\left(p^{2m}\sum_{\substack{u=1\\
p\nmid u
}
}^{p^{r-m}-1}\zeta_{p^{r-m}}^{-u(j-1)}\right)\label{eq:Eigenvalue Sum}
\end{equation}

Let $\ell=\nu_{p}(j-1)$. We will simplify (\ref{eq:Eigenvalue Sum}) using Lemma \ref{Sum of RoU Coprime to p} with $J=j-1$, and $N=1,\dots,r$. For ease of exposition, we split into two different cases depending on $\ell$ and will consequently obtain our desired expression for $\lambda_{j}$. 

\textbf{Case 1}: Suppose $\ell=0$. By Lemma \ref{Sum of RoU Coprime to p}, each sum in (\ref{eq:Eigenvalue Sum}) is zero except for the last sum corresponding to $m=r-1$. We get 
\[
\lambda_{j}=-p^{2r-1}+p^{2(r-1)}\sum_{\substack{u=1\\
p\nmid u
}
}^{p-1}\zeta_{p}^{-u(j-1)}=-p^{2r-1}+p^{2(r-1)}(-1)=-p^{2r-1}-p^{2(r-1)}.
\]

\textbf{Case 2}: Suppose $1\leq\ell\leq r-1$. Since $p^{\ell+1}\nmid(j-1)$, we have 
\[
\sum_{\substack{u=1\\
p\nmid u
}
}^{p^{N}-1}\zeta_{p^{N}}^{-u(j-1)}=0
\]
for all $N\geq\ell+2$. Expression (\ref{eq:Eigenvalue Sum}) becomes
\begin{align*}	
    & =	\sum_{m=1}^{\ell+1}\left(p^{2(r-m)}\sum_{\substack{u=1\\p\nmid u}}^{p^{m}-1}\zeta_{p^{m}}^{-u(j-1)}\right)\\
	& =	p^{2(r-(\ell+1))}(-p^{\ell})+\sum_{m=1}^{\ell}p^{2(r-m)}(p^{m}-p^{m-1})\\
	& =	-p^{2(r-1)-\ell}+(p-1)p^{2r}\sum_{m=2}^{\ell+1}p^{-m}\\
	& =	-p^{2(r-1)-\ell}+(p-1)p^{2r}\cdot\frac{p^{-(\ell+1)}(p^{\ell}-1)}{p-1}\\
	& =	-p^{2(r-1)-\ell}+p^{2r-(\ell+1)}(p^{\ell}-1).
\end{align*}
We conclude
\begin{align*}
\lambda_{j} & = (-p^{2r-1})-p^{2(r-1)-\ell}+p^{2r-(\ell+1)}(p^{\ell}-1)\\
 & = -p^{2r-1}-p^{2r-\ell-2}+p^{2r+1}-p^{2r-\ell-1}\\
 & = -p^{2r-2-\ell}(p+1)\\
 & = -p^{2r-2-\nu_{p}(j-1)}(p+1)
\end{align*}
as desired. 
\end{proof}

\begin{cor}
$M(p^{r})$ is invertible.
\end{cor}
\begin{proof}
By Lemma \ref{Eigenvalues of M(p^r)}, all the eigenvalues of $M(p^{r})$ are nonzero hence $M(p^{r})$ is invertible.
\end{proof}
Now that we know the eigenvalues of $M(p^{r})$, we can use Proposition \ref{inv of circ mat} to compute the inverse of $M(p^{r})$.

\begin{prop}\label{Inverse of M(p^r)}
Let $b_{i,j}$ denote the $(i,j)$-entry of $M(p^{r})^{-1}$. We have 
\[
b_{i,j}=\begin{cases}
-p^{1-2r}-\frac{p-1}{p+1}rp^{1-2r} & \mbox{if }i=j\\
\\
-p^{1-2r}-\frac{p^{-3r+2}}{p+1}\cdot(-p^{r-1}+\nu_{p}(i-j)p^{r-1}(p-1)) & \mbox{otherwise}
\end{cases}.
\]
\end{prop}

\begin{proof}
By Lemma \ref{inv of circ mat}, the $(i,j)$-entry of $M(p^{r})^{-1}$ is equal to
\[
\frac{1}{p^{r}}\sum_{k=1}^{p^{r}}\lambda_{k}^{-1}\zeta_{p^{r}}^{(k-1)(i-j)}
\]
where $\lambda_{k}$ are the eigenvalues of $M(p^{r})$ as in Lemma \ref{Eigenvalues of M(p^r)}. Continuing,
\begin{align*}
 & = \frac{1}{p^{r}}\left(-p^{1-r}+\sum_{k=2}^{p^{r}}\frac{-1}{p+1}p^{-2r+2+\nu_{p}(k-1)}\zeta_{p^{r}}^{(k-1)(i-j)}\right)\\
 & = -p^{1-2r}-\frac{p^{-3r+2}}{p+1}\sum_{k=2}^{p^{r}}p^{\nu_{p}(k-1)}\zeta_{p^{r}}^{(k-1)(i-j)}\\
 & = -p^{1-2r}-\frac{p^{-3r+2}}{p+1}\sum_{k=1}^{p^{r}-1}p^{\nu_{p}(k)}\zeta_{p^{r}}^{k(i-j)}.
\end{align*}
We split into two cases, breaking down the sum in a similar manner as in the proof of Lemma \ref{Eigenvalues of M(p^r)}. We have 
\[
\sum_{k=1}^{p^{r}-1}p^{\nu_{p}(k)}\zeta_{p^{r}}^{k(i-j)}=\sum_{m=0}^{r-1}\left(p^{m}\sum_{\substack{u=1\\
p\nmid u
}
}^{p^{r-m}-1}\zeta_{p^{r-m}}^{u(i-j)}\right)
\] 

\textbf{Case 1}: Suppose $i=j$. Then
\begin{align*}
\sum_{k=1}^{p^{r}-1}p^{\nu_{p}(k)}\zeta_{p^{r}}^{k(i-j)} & = \sum^{r-1}_{m=0}p^m\sum_{\substack{u=1\\
p\nmid u
}
}^{p^{r-m}-1}1 \\
& = \sum_{m=0}^{r-1}p^{r-1}(p-1)\\
& = rp^{r-1}(p-1)
\end{align*}
Thus the $(i,i)$-entry of $M(p^{r})^{-1}$ is 
\[
-p^{1-2r}-\frac{p^{-3r+2}}{p+1}rp^{r-1}(p-1)=-p^{1-2r}-\frac{p-1}{p+1}\cdot r\cdot p^{1-2r}.
\]

\textbf{Case 2}: Suppose $i\ne j$ and let $\ell=\nu_{p}(i-j)$. Then this situation resembles that of Equation (\ref{eq:Eigenvalue Sum}) which we have already computed.
\begin{align*}
\sum_{k=1}^{p^{r}-1}p^{\nu_{p}(k)}\zeta_{p^{r}}^{k(i-j)} & = \sum_{m=1}^{\ell+1}p^{r-m}\left(\sum_{\substack{u=1\\
p\nmid u
}
}^{p^{m}-1}\zeta_{p^{m}}^{u(i-j)}\right)\\
& = \sum_{m=1}^{\ell+1}p^{r-m}p^{m-1}(p-1)\\
& = -p^{r-1}+\ell p^{r-1}(p-1)
\end{align*}
Thus the $(i,j)$-entry of $M(p^{r})^{-1}$ is 
\[
-p^{1-2r}-\frac{p^{-3r+2}}{p+1}\left(-p^{r-1}+\ell p^{r-1}(p-1)\right)\]
\[
=-p^{1-2r}-\frac{p^{-3r+2}}{p+1}\left(-p^{r-1}+\nu_{p}(i-j)p^{r-1}(p-1)\right)
\]
as desired.
\end{proof}

\subsection{\label{sub:Inverting T}\texorpdfstring{Inverting $M(p^{r})_{\hat{1},\hat{1}}$}{}}

Having calculated the entries of $M(p^{r})^{-1}$, we can calculate the entries of $M(p^{r})_{\hat{1},\hat{1}}^{-1}$ using the following result. We will provide a sketch of the proof. A full proof can be found in \cite[Theorem 2.2]{[JCP16]}.

\begin{prop}\label{JCP16}
Let $A$ be an invertible $n\times n$ matrix and let $A^{-1}=(m_{ij})$. Let $s,t\in\left\{ 1,\dots,n\right\} $ and let $A_{\hat{s},\hat{t}}$ denote the matrix obtain from $A$ by removing the $s^{{\rm th}}$ row and $t^{{\rm th}}$ column. Then the $(i,j)$-entry of $A_{\hat{s},\hat{t}}^{-1}=(a_{ij})$ is given by 
\[
a_{ij}=m_{ij}-\frac{m_{is}m_{tj}}{m_{ts}}
\]
for $i,j=1,\dots,n$ with $i\ne t$ and $j\ne s$.
\end{prop}
\begin{proof}
Write $A=(w_{ij})$. Let $u$ denote the $s^{{\rm th}}$ column of $A^{-1}$ after removing the $t^{{\rm th}}$ component and let $v$ denote the $s^{{\rm th}}$ row of $A$ after removing the $t^{{\rm th}}$ component. Then one can verify 
\[
(A_{\hat{s},\hat{t}})^{-1}=\left(I_{n-1}-uv^{T}\right)^{-1}(A^{-1})_{\hat{t},\hat{s}}
\]
where $I_{n-1}$ is the $(n-1)\times(n-1)$ identity matrix. Using the Sherman-Morrison formula, which is a special case of Proposition \ref{WMI}, to calculate $(I_{n-1}-uv^{T})^{-1}$, we get 
\begin{align*}
(A_{\hat{s},\hat{t}})^{-1} & = \left(I_{n-1}+\frac{uv^{T}}{1-v^{T}u}\right)(A^{-1})_{\hat{t},\hat{s}}\\
 & = \left(I_{n-1}+\frac{uv^{T}}{w_{st}m_{ts}}\right)(A^{-1})_{\hat{t},\hat{s}.}
\end{align*}
Therefore the $(i,j)$ entry of $(A_{\hat{s},\hat{t}})^{-1}$ is
\begin{alignat*}{1}
a_{ij} & =m_{ij}+\frac{m_{is}}{w_{st}m_{ts}}\sum_{k\ne q}w_{sk}m_{kj}\\
 & =m_{ij}-\frac{m_{is}m_{tj}}{m_{ts}}\tag*{\qedhere}
\end{alignat*}
\end{proof}

\begin{prop}\label{Inverse of M(p^r)11}
Let $(a_{i,j})=M(p^{r})_{\hat{1},\hat{1}}^{-1}$ and let $\ell_{i}=\nu_{p}(i)$. We have 
\[
a_{i,j}=\begin{cases}
{\displaystyle -p^{1-2r}-\frac{p-1}{p+1}rp^{1-2r}+\frac{p^{1-2r}(\ell_{i}p-\ell_{i}+p)^{2}}{(p+1)(pr+p-r+1)}} & \mbox{if }i=j\\
\\
{\displaystyle -\frac{p^{1-2r}(\ell_{i-j}p+p-\ell_{i-j})}{p+1}+\frac{p^{1-2r}(\ell_{i}p+p-\ell_{i})(\ell_{j}p+p-\ell_{j})}{(p+1)(pr+p-r+1)}} & \mbox{otherwise}
\end{cases}
\]
\end{prop}
\begin{proof}
We will apply Proposition \ref{Inverse of M(p^r)} and Proposition \ref{JCP16} to compute $a_{i,j}$. When $i=j$, we have
\begin{align*}
a_{i,i} & = m_{i,i}-\frac{m_{i,1}m_{1,i}}{m_{1,1}}\\
 & = -p^{1-2r}-\frac{p-1}{p+1}rp^{1-2r}-\frac{\left(-p^{1-2r}-\frac{p^{-3r+2}}{p+1}\cdot(-p^{r-1}+\ell_{i+1}p^{r-1}(p-1))\right)^{2}}{-p^{1-2r}-\frac{p-1}{p+1}rp^{1-2r}}\\
 & = -p^{1-2r}-\frac{p-1}{p+1}rp^{1-2r}+\frac{p^{1-2r}(\ell_{i+1}p-\ell_{i+1}+p)^{2}}{(p+1)(pr+p-r+1)}.
\end{align*}
When $i\ne j$, we have
\begin{align*}
a_{i,j} & = m_{i,j}-\frac{m_{i,1}m_{1,j}}{m_{1,1}}\\
 & = -p^{1-2r}-\frac{p^{-3r+2}}{p+1}\cdot(-p^{r-1}+\ell_{i-j}p^{r-1}(p-1))\\ 
 &   -\frac{\left(-p^{1-2r}-\frac{p^{-3r+2}}{p+1}\cdot(-p^{r-1}+\ell_{i+1}p^{r-1}(p-1))\right)\left(-p^{1-2r}-\frac{p^{-3r+2}}{p+1}\cdot(-p^{r-1}+\ell_{j+1}p^{r-1}(p-1))\right)}{\left(-p^{1-2r}-\frac{p-1}{p+1}rp^{1-2r}\right)}\\
 & = -\frac{p^{1-2r}(\ell_{i-j}p-\ell_{i-j}+p)}{p+1}+\frac{\frac{p^{1-2r}(\ell_{i+1}p-\ell_{i+1}+p)}{p+1}\cdot\frac{p^{1-2r}(\ell_{j+1}p-\ell_{j+1}+p)}{p+1}}{\frac{p^{1-2r}(pr+p-r+1)}{p+1}}\\
 & = -\frac{p^{1-2r}(\ell_{i-j}p+p-\ell_{i-j})}{p+1}+\frac{p^{1-2r}(\ell_{i+1}p+p-\ell_{i+1})(\ell_{j+1}p+p-\ell_{j+1})}{(p+1)(pr+p-r+1)}.
\end{align*}
Note that the indices $i,j$ in Proposition \ref{JCP16} range $2\leq i,j\leq p^{r}$ in the situation $s=t=1$; we will shift our index down by 1 so that $1\leq i,j\leq p^{r}-1$, giving our desired expression for $a_{i,j}$ in the statement of the proposition.
\end{proof}

\subsection{\texorpdfstring{Inverting $T$}{}}
Recall in Section \ref{sub:Inverting M(p^r)} we wrote $T=A+N$ where
\[
A=\left(\begin{array}{cc}
M(p^{r})_{\hat{1},\hat{1}} & 0\\
0 & p^{2}M(p^{r-1})
\end{array}\right)\mbox{ and }N=\left(\begin{array}{cc}
0 & \mathbf{1}_{p^{r}-1\times p^{r-1}}\\
\mathbf{1}_{p^{r-1}\times p^{r}-1} & 0
\end{array}\right).
\]
The matrix $N$ is rank 2 and can be written as $N=UI_{2}V$ where $I_{2}$ is the $2\times2$ identity matrix, $U$ is the $(p^{r}-1+p^{r-1})\times2$ matrix whose first and last column are the same as those of $N$,
and $V$ is the $2\times(p^{r}-1+p^{r-1})$ matrix 
\[
V=\left(\begin{array}{cccccccc}
1 & 1 & \cdots & 1 & 0 & 0 & \cdots & 0\\
0 & 0 & \cdots & 0 & 1 & 1 & \cdots & 1
\end{array}\right)
\]
where the first $p^{r}-1$ entries of the first row of $V$ are all 1 with the remaining $p^{r-1}$ entries are all 0 and the first $p^{r}-1$ entries of the second row are 0 while the remaining $p^{r-1}$ entries are all 1. Note that $V=U^{T}$.

In the more general situation where $A$ is an invertible block diagonal matrix with
\[
A^{-1}=\left(\begin{array}{cc}
\begin{array}{ccc}
a_{11} & \cdots & a_{1n}\\
\vdots & \ddots & \vdots\\
a_{n1} & \cdots & a_{nn}
\end{array} & \mathbf{0}\\
\mathbf{0} & \begin{array}{ccc}
b_{11} & \cdots & b_{1m}\\
\vdots & \ddots & \vdots\\
b_{m1} & \cdots & b_{mm}
\end{array}
\end{array}\right)
\]
an explicit formula for the entries of $(A+UV)^{-1}$ is provided in Proposition \ref{Gen Formula for T^-1}. We state it here for convenience:

\begin{prop}\label{General Formula for T^-1}
Let $T=A+UV$ where $A,U,$ and $V$ are the given matrices above. Let $c^{i,j}$ denote the $(i,j)$-entry of $T^{-1}$. We have 
\[
c^{i,j}=\begin{cases}
{\displaystyle a_{i,j}+\frac{\beta}{1-\alpha\beta}\left(\sum_{k=1}^{n}a_{i,k}\right)\left(\sum_{k=1}^{n}a_{k,j}\right)} & \mbox{if }1\leq i,j\leq n\\
{\displaystyle \frac{-1}{1-\alpha\beta}\left(\sum_{k=1}^{n}a_{i,k}\right)\left(\sum_{k=1}^{m}b_{k,j}\right)} & \mbox{if }1\leq i\leq n\mbox{ and }n<j\leq n+m\\
{\displaystyle \frac{-1}{1-\alpha\beta}\left(\sum_{k=1}^{m}b_{i,k}\right)\left(\sum_{k=1}^{n}a_{k,j}\right)} & \mbox{if }1\leq j\leq n\mbox{ and }n<i\leq n+m\\
{\displaystyle b_{i,j}+\frac{\alpha}{1-\alpha\beta}\left(\sum_{k=1}^{m}b_{i,k}\right)\left(\sum_{k=1}^{m}b_{k,j}\right)} & \mbox{if }n<i,j\leq n+m
\end{cases}
\]
where $\alpha=\sum a_{i,j}$ is the sum of all entries in the first block in $A^{-1}$ and $\beta=\sum b_{i,j}$ is the sum of all entries in the second block in $A^{-1}$.
\end{prop}

To obtain a closed formula for the entries $c^{i,j}$, we will therefore need to calculate the row and column sums of $M(p^{r})_{\hat{1},\hat{1}}^{-1}$ and $p^{-2}M(p^{r-1})^{-1}$. Note that both these matrices are symmetric so it suffices to compute, say, the row sums. The following lemma will be used when computing these sums.

\begin{lem}\label{Sum of P-adic Vals}
We have 
\[
\sum_{m=1}^{p^{r}-1}\nu_{p}(m)=\frac{p^{r}-pr+r-1}{p-1}.
\]
\end{lem}
\begin{proof}
The number of positive integers $<p^{r}$ with valuation $\nu_{p}(m)=\ell$ is precisely $\varphi(p^{r-\ell})=p^{r-\ell-1}(p-1)$. Hence 
\[
\sum_{m=1}^{p^{r}-1}\nu_{p}(m)=\sum_{\ell=0}^{r-1}\ell p^{r-\ell-1}(p-1)=p^{r-1}(p-1)\sum_{\ell=0}^{r-1}\ell p^{-\ell}.
\]

Using the identity 
\[
\sum_{\ell=0}^{n}\ell X^{\ell}=\frac{X(nX^{n+1}-(n+1)X^{n}+1)}{(X-1)^{2}}
\]
with $X=p^{-1}$ and $n=r-1$, we obtain
\begin{alignat*}{1}
p^{r-1}(p-1)\sum_{\ell=0}^{r-1}\ell p^{-\ell} & =p^{r-1}(p-1)\cdot\frac{p^{-1}((r-1)p^{-r}-rp^{-(r-1)}+1)}{(p^{-1}-1)^{2}}\\
 & =\frac{p^{r}-pr+r-1}{p-1}.\tag*{\qedhere}
 \end{alignat*}
\end{proof}

\begin{cor}\label{Sum of Padic Vals i not j}
Let $1\leq i\leq p^{r}-1$. We have
\[
\sum_{\substack{j=1\\
j\ne i
}
}^{p^{r}-1}\nu_{p}(j-i)=-\nu_{p}(i)+\frac{p^{r}-pr+r-1}{p-1}.
\]
\end{cor}
\begin{proof}
We have 
\[
\sum_{\substack{j=1\\
j\ne i
}
}^{p^{r}-1}\nu_{p}(j-i)=\sum_{\substack{j=1-i\\
j\ne0
}
}^{p^{r}-1-i}\nu_{p}(j)=\sum_{j=1-i}^{-1}\nu_{p}(j)+\sum_{j=1}^{p^{r-1}-1-i}\nu_{p}(j).
\]
Since $\nu_{p}(j)<r$, we have $\nu_{p}(j)=\nu_{p}(j+p^{r})$. Continuing,
\begin{align*}
 & = \sum_{j=1-i}^{-1}\nu_{p}(j+p^{r})+\sum_{j=1}^{p^{r-1}-i-1}\nu_{p}(j)\\
 & = \sum_{j=p^{r}-i+1}^{p^{r}-1}\nu_{p}(j)+\sum_{j=1}^{p^{r}-i-1}\nu_{p}(j)\\
 & = -\nu_{p}(p^{r}-i)+\sum_{j=1}^{p^{r}-1}\nu_{p}(j)\\
 & = -\nu_{p}(i)+\sum_{j=1}^{p^{r}-1}\nu_{p}(j).
\end{align*}
By Lemma \ref{Sum of P-adic Vals}, we finally get 
\[
=-\nu_{p}(i)+\frac{p^{r}-pr+r-1}{p-1}.\qedhere
\]
\end{proof}

Write $M(p^{r})_{\hat{1},\hat{1}}^{-1}=(a_{ij})$ and $p^{-2}M(p^{r-1})^{-1}=(b_{ij})$. We first compute the row and column sums of $p^{-2}M(p^{r-1})^{-1}$ and the quantity $\beta$, the sum of all the entries of $(b_{ij})$. 

\begin{lem}\label{Row Sums for b_ij}
Consider the matrix $p^{-2}M(p^{r-1})^{-1}=(b_{ij})$. We have 
\[
\sum_{m=1}^{p^{r-1}}b_{i,m}=-p^{-r}\mbox{ and }\beta=-p^{-1}
\]
for all $1\leq i\leq p^{r-1}$.
\end{lem}
\begin{proof}
Since $M(p^{r-1})$ is circulant, the inverse $M(p^{r-1})^{-1}$ is circulant by Corollary \ref{Inv of Circ is Circ}. Therefore all the row sums are the same. Furthermore, by Lemma \ref{Eigenvalues of M(p^r)}, $\lambda_{1}=-p^{r-2}$ is an eigenvalue of $M(p^{r-1})$ with corresponding eigenvector $v_{1}=(1,1,\dots,1)$ (see Lemma \ref{Evals of Circulant}). 

Note that the entries of $M(p^{r-1})v_{1}=-p^{r-2}v_{1}$ are precisely the row sums of $M(p^{r-1})$. Hence the row sums of $M(p^{r-1})$ are all $-p^{r-2}$. Observe that 
\[
M(p^{r-1})^{-1}v_{1}=-p^{2-r}v_{1}
\]
so the row sums of $M(p^{r-1})^{-1}$ are all $-p^{2-r}$. Consequently the row sums of $p^{-2}M(p^{r-1})^{-1}$ are all $-p^{-r}$. The matrix $p^{-2}M(p^{r-1})^{-1}$ has $p^{r-1}$ rows so 
\[
\beta=(p^{r-1})(-p^{-r})=-p^{-1}
\]
as desired. 
\end{proof}

Next we compute the row sums for $M(p^{r})_{\hat{1},\hat{1}}^{-1}$ which is substantially more tedious than Lemma \ref{Row Sums for b_ij}, noting that $M(p^{r})_{\hat{1},\hat{1}}$ fails to be circulant in general.

\begin{lem}\label{Row Sums for a_ij}
Consider the matrix $M(p^{r})_{\hat{1},\hat{1}}^{-1}=(a_{ij})$ and let $\ell_{i}=\nu_{p}(i)$. Fix a row $i\geq1$. We have 
\[
\sum_{j=1}^{p^{r}-1}a_{ij}=\frac{p^{1-r}(-pr+r-1)}{(p-1)r+p+1}+\left(\frac{p^{1-r}(p-1)}{(p-1)r+p+1}\right)\ell_{i}
\]
and
\[
\alpha=-\frac{(p-1)pr}{(p-1)r+p+1}.
\]
\end{lem}
\begin{proof}
Using Proposition \ref{Inverse of M(p^r)11}, we compute
\begin{align}
\sum_{j=1}^{p^{r}-1}a_{ij}            &=a_{ii}+\sum_{\substack{j=1\\j\ne i}
}^{p^{r}-1}a_{ij}\nonumber \\
& =-p^{1-2r}-\frac{p-1}{p+1}rp^{1-2r}+\frac{p^{1-2r}(\nu_{p}(i)p-\nu_{p}(i)+p)^{2}}{(p+1)(pr+p-r+1)}\nonumber \\
&+\sum_{\substack{j=1\\j\ne i}
}^{p^{r}-1}\left(-\frac{p^{1-2r}(\nu_{p}(i-j)p+p-\nu_{p}(i-j))}{p+1}+\frac{p^{1-2r}(\nu_{p}(i)p+p-\nu_{p}(i))(\nu_{p}(j)p+p-\nu_{p}(j))}{(p+1)(pr+p-r+1)}\right)\nonumber \\
& =-p^{1-2r}-\frac{p-1}{p+1}rp^{1-2r}+\frac{p^{1-2r}(\nu_{p}(i)p-\nu_{p}(i)+p)^{2}}{(p+1)(pr+p-r+1)}\label{eq:Three Expressions}\\
& +\sum_{\substack{j=1\\j\ne i}
}^{p^{r}-1}\left(-\frac{p^{1-2r}(\nu_{p}(i-j)p+p-\nu_{p}(i-j))}{p+1}\right)\label{eq: Sum_1}\\
& +\sum_{\substack{j=1\\j\ne i}
}^{p^{r}-1}\left(\frac{p^{1-2r}(\nu_{p}(i)p+p-\nu_{p}(i))(\nu_{p}(j)p+p-\nu_{p}(j))}{(p+1)(pr+p-r+1)}\right).\label{eq: Sum_2}
\end{align}
We will now compute the sums above, starting with (\ref{eq: Sum_1}).
\small{
\begin{alignat*}{1}
\sum_{\substack{j=1\\
j\ne i
}
}^{p^{r}-1}\left(-\frac{p^{1-2r}(\nu_{p}(i-j)p+p-\nu_{p}(i-j))}{p+1}\right) & =-\frac{p^{1-2r}}{p+1}\sum_{\substack{j=1\\
j\ne i
}
}^{p^{r}-1}(\nu_{p}(i-j)p+p-\nu_{p}(i-j))\\
 & =-\frac{p^{1-2r}}{p+1}\left[\sum_{\substack{j=1\\
j\ne i
}
}^{p^{r}-1}p+(p-1)\sum_{\substack{j=1\\
j\ne i
}
}^{p^{r}-1}\nu_{p}(j-i)\right]\\
 & =-\frac{p^{1-2r}}{p+1}\left[p(p^{r}-2)+(p-1)\sum_{\substack{j=1\\
j\ne i
}
}^{p^{r}-1}\nu_{p}(j-i)\right].
\end{alignat*}
}
\normalsize
Using Corollary \ref{Sum of Padic Vals i not j}, we have
\begin{alignat*}{1}
 & =-\frac{p^{1-2r}}{p+1}\left[p(p^{r}-2)+(p-1)\left(-\nu_{p}(i)+\frac{p^{r}-pr+r-1}{p-1}\right)\right]\\
 & =-\frac{p^{1-2r}}{p+1}\left[p(p^{r}-2)-(p-1)\nu_{p}(i)+p^{r}-pr+r-1\right].
\end{alignat*}

Now we compute (\ref{eq: Sum_2}). Note that the third expression in (\ref{eq:Three Expressions}) is what would be the $j=i$ term in the sum. We can combine them and compute instead

\begin{alignat*}{1}
 & \sum_{j=1}^{p^{r}-1}\frac{p^{1-2r}(\nu_{p}(i)p+p-\nu_{p}(i))(\nu_{p}(j)p+p-\nu_{p}(j))}{(p+1)(pr+p-r+1)}\\
 & =\frac{p^{1-2r}(\nu_{p}(i)p+p-\nu_{p}(i))}{(p+1)(pr+p-r+1)}\sum_{j=1}^{p^{r}-1}(\nu_{p}(j)p+p-\nu_{p}(j))\\
 & =\frac{p^{1-2r}(\nu_{p}(i)p+p-\nu_{p}(i))}{(p+1)(pr+p-r+1)}\sum_{j=1}^{p^{r}-1}(p+\nu_{p}(j)(p-1))\\
 & =\frac{p^{2-2r}(\nu_{p}(i)p+p-\nu_{p}(i))}{(p+1)(pr+p-r+1)}(p^{r}-1)+\frac{p^{1-2r}(p-1)(\nu_{p}(i)p+p-\nu_{p}(i))}{(p+1)(pr+p-r+1)}\sum_{j=1}^{p^{r}-1}\nu_{p}(j)\\
 & =\frac{p^{2-2r}(\nu_{p}(i)p+p-\nu_{p}(i))}{(p+1)(pr+p-r+1)}(p^{r}-1)+\frac{p^{1-2r}(p-1)(\nu_{p}(i)p+p-\nu_{p}(i))}{(p+1)(pr+p-r+1)}\left(\frac{p^{r}-pr+r-1}{p-1}\right)\\
 & =(\nu_{p}(i)p+p-\nu_{p}(i))\left(\frac{p^{2-2r}(p^{r}-1)}{(p+1)(pr+p-r+1)}+\frac{p^{1-2r}(p^{r}-pr+r-1)}{(p+1)(pr+p-r+1)}\right)\\
 & =(\nu_{p}(i)p+p-\nu_{p}(i))p^{1-2r}\left(\frac{p^{r}}{(p-1)r+p+1}-\frac{1}{p+1}\right).
\end{alignat*}
Combining everything together, the $i$th row sum is

\footnotesize{
\begin{alignat*}{1}
\sum_{j=1}^{p^{r}-1}a_{ij} & =-p^{1-2r}-\frac{p-1}{p+1}rp^{1-2r}-\frac{p^{1-2r}}{p+1}\left[p(p^{r}-2)-(p-1)\nu_{p}(i)+p^{r}-pr+r-1\right]\\
 & +(\nu_{p}(i)p+p-\nu_{p}(i))p^{1-2r}\left(\frac{p^{r}}{(p-1)r+p+1}-\frac{1}{p+1}\right)\\
 & =-p^{1-2r}-\frac{p-1}{p+1}rp^{1-2r}-\frac{p^{1-2r}}{p+1}p(p^{r}-2)+\frac{p^{1-2r}}{p+1}(p-1)\nu_{p}(i)-\frac{p^{1-2r}}{p+1}(p^{r}-pr+r-1)\\
 & +((p-1)\nu_{p}(i)+p)\left(\frac{p^{r}p^{1-2r}}{(p-1)r+p+1}-\frac{p^{1-2r}}{p+1}\right)\\
 & =-p^{1-2r}-\frac{p-1}{p+1}rp^{1-2r}-\frac{p^{1-2r}}{p+1}p(p^{r}-2)+\frac{p^{1-2r}}{p+1}(p-1)\nu_{p}(i)-\frac{p^{1-2r}}{p+1}(p^{r}-pr+r-1)\\
 & +(p-1)\left(\frac{p^{1-r}}{(p-1)r+p+1}-\frac{p^{1-2r}}{p+1}\right)\nu_{p}(i)+\left(\frac{p^{1-r}}{(p-1)r+p+1}-\frac{p^{1-2r}}{p+1}\right)p\\
 & =-p^{1-2r}-\frac{p-1}{p+1}rp^{1-2r}-\frac{p^{1-2r}}{p+1}p(p^{r}-2)+\frac{p^{1-2r}}{p+1}(p-1)\nu_{p}(i)-\frac{p^{1-2r}}{p+1}(p^{r}-pr+r-1)\\
 & +\left(\frac{p^{1-r}(p-1)}{(p-1)r+p+1}-\frac{p^{1-2r}(p-1)}{p+1}\right)\nu_{p}(i)+\left(\frac{p^{2-r}}{(p-1)r+p+1}-\frac{p^{2-2r}}{p+1}\right)\\
 & =-p^{1-2r}-\frac{p-1}{p+1}rp^{1-2r}-\frac{p^{1-2r}}{p+1}p(p^{r}-2)-\frac{p^{1-2r}}{p+1}(p^{r}-pr+r-1)\\
 & +\left(\frac{p^{2-r}}{(p-1)r+p+1}-\frac{p^{2-2r}}{p+1}\right)\\
 & =\frac{p^{1-2r}(p-(p+1)p^{r})}{p+1}+\left(\frac{p^{2-r}}{(p-1)r+p+1}-\frac{p^{2-2r}}{p+1}\right)+\left(\frac{(p-1)p^{1-r}}{(p-1)r+p+1}\right)\nu_{p}(i)\\
 & =\frac{p^{1-r}(-pr+r-1)}{(p-1)r+p+1}+\left(\frac{p^{1-r}(p-1)}{(p-1)r+p+1}\right)\nu_{p}(i).
\end{alignat*}
}
\normalsize
Lastly, we compute $\alpha$:
\begin{alignat*}{1}
\alpha & =\sum_{i=1}^{p^{r}-1}\left(\frac{p^{1-r}(-pr+r-1)}{(p-1)r+p+1}+\left(\frac{p^{1-r}(p-1)}{(p-1)r+p+1}\right)\nu_{p}(i)\right)\\
 & =\sum_{i=1}^{p^{r}-1}\frac{p^{1-r}(-pr+r-1)}{(p-1)r+p+1}+\frac{p^{1-r}(p-1)}{(p-1)r+p+1}\sum_{i=1}^{p^{r}-1}\nu_{p}(i)\\
 & =\frac{p^{1-r}(-pr+r-1)}{(p-1)r+p+1}(p^{r}-1)+\frac{p^{1-r}(p-1)}{(p-1)r+p+1}\left(\frac{p^{r}-pr+r-1}{p-1}\right)\\
 & =-\frac{(p-1)pr}{(p-1)r+p+1}.\tag*{\qedhere}
\end{alignat*}
\end{proof}

We will now provide an explicit description of the entries of $T^{-1}$. 

\begin{thm}\label{Entries of Inv Int Matrix}
Let $c^{i,j}$ denote the $(i,j)$-entry of $T^{-1}$. We have {\footnotesize{}
\[
\deg{\rm S}(N)c^{ij}=\begin{cases}
{\displaystyle -\frac{2p^{1-2r}(pr-r+1)}{p+1}+\frac{2p^{1-2r}(p-1)}{p+1}\nu_{p}(i)} & \mbox{if }1\leq i,j\leq p^{r}-1\mbox{ and }i=j\\
\\
{\displaystyle -p^{1-2r}\frac{(pr-r+1)}{p+1}-\frac{p^{1-2r}(p-1)}{p+1}\nu_{p}\left(\frac{1}{j}-\frac{1}{i}\right)} & \mbox{if }1\leq i,j\leq p^{r}-1\mbox{ and }i\ne j\\
\\
{\displaystyle -\frac{p^{1-2r}(pr-r+1)}{p+1}+\frac{p^{1-2r}(p-1)}{p+1}\nu_{p}(i)} & \mbox{if }1\leq i\leq p^{r}-1\mbox{ and }p^{r}\leq j\leq p^{r}-1+p^{r-1}\\
\\
{\displaystyle -\frac{p^{1-2r}(pr-r+1)}{p+1}+\frac{p^{1-2r}(p-1)}{p+1}\nu_{p}(j)} & \mbox{if }1\leq j\leq p^{r}-1\mbox{ and }p^{r}\leq i\leq p^{r}-1+p^{r-1}\\
\\
{\displaystyle -\frac{2p^{1-2r}(pr-r+1)}{p+1}} & \mbox{if }p^{r}\leq i,j\leq p^{r}-1+p^{r-1}\mbox{ and }i=j\\
\\
{\displaystyle -\frac{p^{1-2r}(pr+p-r)}{p+1}-\frac{p^{1-2r}(p-1)}{p+1}\nu_{p}(i-j)} & \mbox{if }p^{r}\leq i,j\leq p^{r}-1+p^{r-1}\mbox{ and }i\ne j
\end{cases}.
\]
}{\footnotesize \par}
\end{thm}

\begin{proof}
We will break into four cases, using Proposition \ref{General Formula for T^-1} to calculate the $c^{i,j}$ along with our results in Lemma \ref{Row Sums for b_ij} and Lemma \ref{Row Sums for a_ij}. 

\textbf{Case 1}. Suppose $1\leq i,j\leq p^{r}-1$. Then 
\[
c^{ij}=a_{i,j}+\frac{\beta}{1-\alpha\beta}\left(\sum_{k=1}^{p^{r}-1}a_{i,k}\right)\left(\sum_{k=1}^{p^{r}-1}a_{k,j}\right).
\]
If $i=j$, then $c^{ii}$ is equal to {\small{}
\begin{alignat*}{1}
 & =-p^{1-2r}-\frac{p-1}{p+1}rp^{1-2r}+\frac{p^{1-2r}(\nu_{p}(i)p-\nu_{p}(i)+p)^{2}}{(p+1)(pr+p-r+1)}+\frac{\beta}{1-\alpha\beta}\left(\frac{p^{1-r}(-pr+r-1)}{(p-1)r+p+1}+\left(\frac{p^{1-r}(p-1)}{(p-1)r+p+1}\right)\nu_{p}(i)\right)^{2}\\
 & =-p^{1-2r}-\frac{p-1}{p+1}rp^{1-2r}+\frac{p^{1-2r}(\nu_{p}(i)p-\nu_{p}(i)+p)^{2}}{(p+1)(pr+p-r+1)}-\frac{pr+p-r+1}{p(p+1)}\left(\frac{p^{2-2r}((p-1)\nu_{p}(i)+(1-p)r-1)^{2}}{(pr+p-r+1)^{2}}\right)\\
 & =-p^{1-2r}-\frac{p-1}{p+1}rp^{1-2r}+\frac{p^{1-2r}(\nu_{p}(i)p-\nu_{p}(i)+p)^{2}}{(p+1)(pr+p-r+1)}-\frac{p^{1-2r}(pr+p-r+1)}{(p+1)}\left(\frac{((p-1)\nu_{p}(i)+(1-p)r-1)^{2}}{(pr+p-r+1)^{2}}\right)\\
 & =-p^{1-2r}-\frac{p-1}{p+1}rp^{1-2r}+\frac{p^{1-2r}(\nu_{p}(i)p-\nu_{p}(i)+p)^{2}}{(p+1)(pr+p-r+1)}-\frac{p^{1-2r}((p-1)\nu_{p}(i)+(1-p)r-1)^{2}}{(p+1)(pr+p-r+1)}\\
 & =-p^{1-2r}-\frac{p-1}{p+1}rp^{1-2r}+\frac{p^{1-2r}(\nu_{p}(i)p-\nu_{p}(i)+p)^{2}-p^{1-2r}((p-1)\nu_{p}(i)+(1-p)r-1)^{2}}{(p+1)(pr+p-r+1)}\\
 & =-p^{1-2r}-\frac{p-1}{p+1}rp^{1-2r}+\frac{(p-1)(2\nu_{p}(i)-r+1)p^{1-2r}}{p+1}\\
 & =\frac{2p^{1-2r}((p-1)\nu_{p}(i)-pr+r-1)}{p+1}\\
 & =-\frac{2p^{1-2r}(pr-r+1)}{p+1}+\frac{2p^{1-2r}(p-1)}{p+1}\nu_{p}(i).
\end{alignat*}
}{\small \par}
If $i\ne j$, then $c^{ij}$ is equal to {\small{}
\begin{alignat*}{1}
 & =-\frac{p^{1-2r}(\nu_{p}(i-j)p+p-\nu_{p}(i-j))}{p+1}+\frac{p^{1-2r}(\nu_{p}(i)p+p-\nu_{p}(i))(\nu_{p}(j)p+p-\nu_{p}(j))}{(p+1)(pr+p-r+1)}+\\
 & +\frac{\beta}{1-\alpha\beta}\left(\frac{p^{1-r}(-pr+r-1)}{(p-1)r+p+1}+\left(\frac{p^{1-r}(p-1)}{(p-1)r+p+1}\right)\nu_{p}(i)\right)\left(\frac{p^{1-r}(-pr+r-1)}{(p-1)r+p+1}+\frac{p^{1-r}(p-1)}{(p-1)r+p+1}\nu_{p}(j)\right)\\
 & =-\frac{p^{1-2r}(\nu_{p}(i-j)p+p-\nu_{p}(i-j))}{p+1}+\frac{p^{1-2r}(\nu_{p}(i)p+p-\nu_{p}(i))(\nu_{p}(j)p+p-\nu_{p}(j))}{(p+1)(pr+p-r+1)}\\
 & -\frac{pr+p-r+1}{p(p+1)}\left(\frac{p^{1-r}(-pr+r-1)}{(p-1)r+p+1}+\left(\frac{p^{1-r}(p-1)}{(p-1)r+p+1}\right)\nu_{p}(i)\right)\left(\frac{p^{1-r}(-pr+r-1)}{(p-1)r+p+1}+\frac{p^{1-r}(p-1)}{(p-1)r+p+1}\nu_{p}(j)\right)\\
 & =-\frac{p^{1-2r}(\nu_{p}(i-j)p+p-\nu_{p}(i-j))}{p+1}+\frac{p^{1-2r}(\nu_{p}(i)p+p-\nu_{p}(i))(\nu_{p}(j)p+p-\nu_{p}(j))}{(p+1)(pr+p-r+1)}\\
 & -\frac{pr+p-r+1}{p(p+1)}\left(\frac{p^{2-2r}(pr-p\nu_{p}(i)-r+\nu_{p}(i)+1)(pr-p\nu_{p}(j)-r+\nu_{p}(j)+1)}{(pr+p-r+1)^{2}}\right)\\
 & =-\frac{p^{1-2r}(\nu_{p}(i-j)p+p-\nu_{p}(i-j))}{p+1}+\frac{p^{1-2r}(\nu_{p}(i)p+p-\nu_{p}(i))(\nu_{p}(j)p+p-\nu_{p}(j))}{(p+1)(pr+p-r+1)}\\
 & -\frac{p^{1-2r}(pr-p\nu_{p}(i)-r+\nu_{p}(i)+1)(pr-p\nu_{p}(j)-r+\nu_{p}(j)+1)}{(p+1)(pr+p-r+1)}\\
 & =-\frac{p^{1-2r}(\nu_{p}(i-j)p+p-\nu_{p}(i-j))}{p+1}-\frac{(p-1)p^{1-2r}(pr+p-r+1)(r-\nu_{p}(i)-\nu_{p}(j)-1)}{(p+1)(pr+p-r+1)}\\
 & =-\frac{p^{1-2r}(\nu_{p}(i-j)p+p-\nu_{p}(i-j))}{p+1}-\frac{(p-1)p^{1-2r}(r-\nu_{p}(i)-\nu_{p}(j)-1)}{p+1}\\
 & =-p^{1-2r}\frac{\nu_{p}(i-j)p+p-\nu_{p}(i-j)+(p-1)(r-\nu_{p}(i)-\nu_{p}(j)-1)}{p+1}\\
 & =-p^{1-2r}\frac{(pr-r+1)}{p+1}-\frac{p^{1-2r}(p-1)}{p+1}\left(\nu_{p}(i-j)-\nu_{p}(i)-\nu_{p}(j)\right)\\
 & =-p^{1-2r}\frac{(pr-r+1)}{p+1}-\frac{p^{1-2r}(p-1)}{p+1}\nu_{p}\left(\frac{1}{j}-\frac{1}{i}\right).
\end{alignat*}
}{\small \par}

\normalsize
\textbf{Case 2}. Suppose $1\leq i\leq p^{r}-1$ and $p^{r}\leq j\leq p^{r}-1+p^{r-1}$.
Then
\begin{align*}
c^{ij} & = \frac{-1}{1-\alpha\beta}\left(\sum_{k=1}^{n}a_{i,k}\right)\left(\sum_{k=1}^{m}b_{k,j}\right)\\
 & = \frac{-1}{1-\alpha\beta}\left(\frac{p^{1-r}(-pr+r-1)}{(p-1)r+p+1}+\left(\frac{p^{1-r}(p-1)}{(p-1)r+p+1}\right)\nu_{p}(i)\right)\left(-p^{-r}\right)\\
 & = -\left(\frac{pr+p-r+1}{p+1}\right)\left(-\frac{p^{1-2r}(-pr+r-1)}{(p-1)r+p+1}-\left(\frac{p^{1-2r}(p-1)}{(p-1)r+p+1}\right)\nu_{p}(i)\right)\\
 & = -\frac{p^{1-2r}(pr-r+1)}{p+1}+\frac{p^{1-2r}(p-1)}{p+1}\nu_{p}(i).
\end{align*}

\textbf{Case 3}. Suppose $p^{r}\leq i\leq p^{r}-1+p^{r-1}$ and $1\leq j\leq p^{r}-1$.
Then
\begin{align*}
c^{ij} & = \frac{-1}{1-\alpha\beta}\left(\sum_{k=1}^{m}b_{i,k}\right)\left(\sum_{k=1}^{n}a_{k,j}\right)\\
 & = -\left(\frac{pr+p-r+1}{p+1}\right)(-p^{-r})\left(\frac{p^{1-r}(-pr+r-1)}{(p-1)r+p+1}+\left(\frac{p^{1-r}(p-1)}{(p-1)r+p+1}\right)\nu_{p}(j)\right)\\
 & = -\frac{p^{1-2r}(pr-r+1)}{p+1}+\frac{p^{1-2r}(p-1)}{p+1}\nu_{p}(j).
\end{align*}

\textbf{Case 4}. Suppose $p^{r}\leq i,j\leq p^{r}-1+p^{r-1}$. If
$i=j$, then
\begin{align*}
c^{ii} & = b_{i,i}+\frac{\alpha}{1-\alpha\beta}\left(\sum_{k=1}^{m}b_{i,k}\right)\left(\sum_{k=1}^{m}b_{k,i}\right)\\
 & = \frac{1}{p^{2}}\left(-p^{3-2r}-\frac{p-1}{p+1}(r-1)p^{3-2r}\right)+\frac{\alpha}{1-\alpha\beta}(-p^{-r})^{2}\\
 & = \frac{2p^{1-2r}(-pr+r-1)}{p+1}=-\frac{2p^{1-2r}(pr-r+1)}{p+1}
\end{align*}
If $i\ne j$, then
\begin{alignat*}{1}
c^{ij} & =\frac{1}{p^{2}}\left(-p^{3-2r}-\frac{p^{5-3r}}{p+1}(-p^{r-2}+\nu_{p}(i-j)p^{r-2}(p-1))\right)+\frac{\alpha}{1-\alpha\beta}(-p^{-r})^{2}\\
 & =-\frac{((p-1)\nu_{p}(i-j)+p)p^{1-2r}}{p+1}-\frac{(p-1)pr}{p+1}(-p^{-r})^{2}\\
 & =-\frac{((p-1)\nu_{p}(i-j)+p)p^{1-2r}+(p-1)p^{1-2r}r}{p+1}\\
 & =-\frac{p^{1-2r}(pr+p-r)}{p+1}-\frac{p^{1-2r}(p-1)}{p+1}\nu_{p}(i-j).\tag*{\qedhere}
\end{alignat*}
\end{proof}

The following corollary will be useful when we find an upper bound for the exponent in Theorem \ref{Upper Bound}.
\begin{cor}\label{Entries cij Are Negative}
Let $c^{i,j}$ denote the $(i,j)$-entry of $T^{-1}$. Then each $c^{i,j}$ is negative. 
\end{cor}
\begin{proof}
Based on our result in Theorem \ref{Entries of Inv Int Matrix}, we will show the third case is negative \emph{i.e. }we will show 
\begin{equation}
-\frac{p^{1-2r}(pr-r+1)}{p+1}+\frac{p^{1-2r}(p-1)}{p+1}\nu_{p}(i)\label{eq: negentry}
\end{equation}
is negative for $1\leq i\leq p^{r}-1$. The other cases are either clearly negative or are essentially the same as this case. 

Since the largest value $\nu_{p}(i)$ attains is $r-1$, the largest value expression (\ref{eq: negentry}) attains is 
\[
-\frac{p^{1-2r}(pr-r+1)}{p+1}+\frac{p^{1-2r}(p-1)}{p+1}(r-1)=-\frac{p^{2-2r}}{p+1}
\]
which is always negative. 
\end{proof}


\section{\label{sec:Computing-deg}Computing the degree of the modular sheaf}

Throughout this chapter, we will keep the notation of §\ref{sec:Intersection-Matrix}. Unless otherwise stated, we let $R=\mathbb{Z}_{p}[\zeta_{Np^{r}}]$. We will compute the degree of $\underline{\omega}^{\otimes2}$ restricted to an irreducible component $\Lambda$ of $\bar{\mathfrak{X}}$ and ultimately compute an upper bound for the exponent $e$. 

\subsection{Decomposing the Modular Sheaf}
We will make use of the Kodaira-Spencer isomorphism, as stated in Theorem \ref{K-S Iso}, which we restate here for convenience.

\begin{thm}
Let $R$ be a noetherian, regular, excellent $\mathbb{Z}[\zeta_{N}]$-algebra containing $1/N$. The Kodaira-Spencer isomorphism  $\underline{\omega}_{\mathfrak{Y}(N)}^{\otimes2}\simeq\Omega_{\mathfrak{Y}(N)/R}^{1}$
on $\mathfrak{Y}(N)$ extends to an isomorphism on $\mathfrak{X}(N)$
\[
\underline{\omega}_{\mathfrak{X}(N)}^{\otimes2}\simeq\Omega_{\mathfrak{X}(N)/R}^{1}(\mathfrak{C}(N)).
\]
\end{thm}

We will need the following definition, which we take from \cite[{}6.4.18]{[Liu02]}.

\begin{defn}\label{def: RDS}
Let $f:X\rightarrow Y$ be a proper morphism of relative dimension $\leq r$. A\textbf{ relative ($r$-)dualizing sheaf }for $f:X\rightarrow Y$ is a quasi-coherent sheaf $\Omega_{f}$ on $X$, endowed with a homomorphism of ${\cal O}_{X}$-modules 
\[
{\rm tr}_{f}:R^{r}f_{*}\Omega_{f}\rightarrow{\cal O}_{Y}
\]
such that for any quasi-coherent sheaf ${\cal F}$ on $X$, the natural bilinear map
\[
f_{*}{\cal H}om_{{\cal O}_{X}}({\cal F},\Omega_{f})\times R^{r}f_{*}{\cal F}\rightarrow R^{r}f_{*}\Omega_{f}\overset{{\rm tr}_{f}}{\longrightarrow}{\cal O}_{Y}
\]
induces an isomorphism
\[
f_{*}{\cal H}om_{{\cal O}_{X}}({\cal F},\Omega_{f})\simeq{\cal H}om_{{\cal O}_{Y}}(R^{r}f_{*}{\cal F},{\cal O}_{Y}).
\]
\end{defn}
By \cite[{}6.4.19]{[Liu02]}, uniqueness of $\Omega_{f}$ is automatic once we have existence. If $Y$ is locally noetherian and $f:X\rightarrow Y$ is a projective morphism with fibers of dimension $\leq r$, then as remarked in \cite[{}6.4.30]{[Liu02]}, the relative $r$-dualizing sheaf exists.  Furthermore, by \cite[Theorem 6.4.32]{[Liu02]},  the relative dualizing sheaf is isomorphic to the canonical sheaf $\Omega_{X/Y}$ (see \cite[{}6.4.7]{[Liu02]}) whenever $f$ is a flat projective l.c.i. and $Y$ is locally noetherian. When $f$ is smooth, the canonical sheaf coincides with the sheaf of Kahler
differentials $\Omega_{X/Y}^{1}$.
We will use the following result, which is in \cite[Theorem 6.4.9]{[Liu02]} known as the adjunction formula, to eventually relate $\underline{\omega}_{\mathfrak{X}(Np^{r})}^{\otimes2}$ with the relative dualizing sheaf of $\mathfrak{X}(Np^{r})$ and of
$\mathfrak{X}(N)$.

\begin{thm}\label{Adj Formula}
Let $f:X\rightarrow Y$ and $g:Y\rightarrow Z$ be quasi-projective l.c.i.s. We have a canonical isomorphism of canonical sheaves
\[
\Omega_{X/Z}\simeq\Omega_{X/Y}\otimes_{{\cal O}_{X}}f^{*}\Omega_{Y/Z}.
\]
\end{thm}

We cannot directly apply the Kodaira-Spencer isomorphism to our modular curve $\mathfrak{X}$ since $p$ and consequently the level $Np^{r}$, is not invertible in $R$. Instead, we will apply it to the modular curve $\mathfrak{X}(N)$ over $R$ since the level $N$ is invertible in $\mathbb{Z}_{p}\subset R$. Consider 
\begin{equation}
\mathfrak{X}\overset{{\rm pr}}{\longrightarrow}\mathfrak{X}(N)_{/R}\overset{g}{\longrightarrow}{\rm Spec}(R)\label{eq: pr map}
\end{equation}
where ${\rm pr}$ is the projection map and $g$ is the structural morphism. For convenience, we let $\mathfrak{X}(N)$ denote the base change $\mathfrak{X}(N)_{/R}$. 

According to \cite[{}6.3.18]{[Liu02]}, if $X\rightarrow Y$ is a morphism of finite type of regular locally noetherian schemes, then $X\rightarrow Y$ is an l.c.i. Therefore the maps ${\rm pr}$ and $g$ are l.c.i.s. By Theorem \ref{X(N) is arith surf}, $g$ and $g\circ{\rm pr}$ are projective. Hence by \cite[{}3.3.32(e)]{[KM85]}, ${\rm pr}$ is projective. Applying the adjunction formula to (\ref{eq: pr map}), we have 
\[
\Omega_{\mathfrak{X}/R}\simeq\Omega_{\mathfrak{X}/\mathfrak{X}(N)}\otimes_{{\cal O}_{\mathfrak{X}}}{\rm pr}^{*}\Omega_{\mathfrak{X}(N)/R}.
\]
Combining this with the Kodaira-Spencer isomorphism applied to $\mathfrak{X}(N)$, we get
\begin{equation}
\Omega_{\mathfrak{X}/R}\simeq\Omega_{\mathfrak{X}/\mathfrak{X}(N)}\otimes_{{\cal O}_{\mathfrak{X}(Np^{r})}}{\rm pr}^{*}\underline{\omega}_{\mathfrak{X}(N)}^{\otimes2}(-\mathfrak{C}(N)).\label{eq:Adj + KS}
\end{equation}
We will later show in Lemma \ref{Pullback by pr} 
\[
{\rm pr}^{*}\underline{\omega}_{\mathfrak{X}(N)}^{\otimes2}(-\mathfrak{C}(N))\simeq\underline{\omega}_{\mathfrak{X}}^{\otimes2}(-p^{r}\mathfrak{C}(Np^{r}))
\]
which is where the sheaf $\underline{\omega}_{\mathfrak{X}}^{\otimes2}$ appears in (\ref{eq:Adj + KS}). This will then allow us to identify $\underline{\omega}_{\mathfrak{X}}^{\otimes2}$ with $\Omega_{\mathfrak{X}/R}(\mathfrak{C}(Np^{r}))$, the relative dualizing sheaf twisted by the cuspidal divisor. Thus, computing $\deg(\underline{\omega}_{\mathfrak{X}}^{\otimes2}|_{\Lambda})$ amounts to computing $\deg(\Omega_{\mathfrak{X}}(\mathfrak{C}(Np^{r}))|_{\Lambda})$. Our first step will be to investigate $\Omega_{\mathfrak{X}/\mathfrak{X}(N)}$. 

\subsection{\label{sub: The RDS}\texorpdfstring{The relative dualizing sheaf $\Omega_{\mathfrak{X}/\mathfrak{X}(N)}$}{}}

In this section, our goal will be to better understand the relative dualizing sheaf $\Omega_{\mathfrak{X}/\mathfrak{X}(N)}$. Since $\Omega_{\mathfrak{X}/\mathfrak{X}(N)}$ is invertible, we have $\Omega_{\mathfrak{X}/\mathfrak{X}(N)}\simeq{\cal O}_{\mathfrak{X}}({\cal R})$ for some divisor ${\cal R}$ of $\mathfrak{X}$. As we will see, ${\cal R}$ will be the divisor associated to the \emph{different} of the morphism ${\rm pr}:\mathfrak{X}\rightarrow\mathfrak{X}(N)_{/R}.$

We begin by discussing the trace map which generalizes the usual notion over a finite extension of fields. Let $A\rightarrow B$ be a finite, flat map of noetherian rings. According to \cite[Tag0BSY]{[Stacks]}, $B$ is a finite locally free $A$-module and so we can consider the trace ${\rm Tr}_{B/A}(b)$ of the $A$-linear map $B\rightarrow B$ given by multiplication by $b$. This gives us an $A$-linear map ${\rm Tr}_{B/A}:B\rightarrow A$. The following definition is from \cite[Tag0BW0]{[Stacks]}. 

\begin{defn}
Let $A\rightarrow B$ be a ring map and let $K={\rm Frac}(A)$, the total ring of fractions of $A$ (see \cite[02C5]{[Stacks]}, note when $A$ is a domain, ${\rm Frac}(A)$ coincides with the field of fractions), and $L=B\otimes_{A}K$. We say the \textbf{Dedekind different is defined }if $A$ is noetherian, $A\rightarrow B$ is finite and maps any non-zerodivisor of $A$ to a non-zerodivisor of $B$, and $K\rightarrow L$ is \'{e}tale. In this situation, $K\rightarrow L$ is finite flat. Let 
\[
{\cal L}_{B/A}=\left\{ x\in L:{\rm Tr}_{L/K}(bx)\in A\mbox{ for all }b\in B\right\} .
\]
We define the \textbf{Dedekind different }of $A\rightarrow B$ to be the inverse of ${\cal L}_{B/A}$: 
\[
\mathfrak{D}_{B/A}={\cal L}_{B/A}^{-1}=\left\{ x\in L:x{\cal L}_{B/A}\subset B\right\} 
\]
viewed as a sub $B$-module of $L$.
\end{defn}

\begin{rmk}
Let $A$ be a Dedekind domain, $K={\rm Frac}(A)$, $L$ a finite separable extension of $K$, and $B$ the integral closure of $A$ in $L$. In this situation, \cite[§4.3]{[Ser79]} defines the different $\mathfrak{D}_{B/A}$ in the same manner as we have done. Since $A$ is normal and noetherian, by \cite[Tag032L]{[Stacks]}, $A\rightarrow B$ is finite. Furthermore, $L={\rm Frac}(B)$ and $L=B\otimes_{A}K$ so indeed, the Dedekind different is defined for $A\rightarrow B$.  We record a few useful facts for calculating the different in this situation.
\end{rmk}

\begin{prop}\label{AKLB Different Props}
Let $A$ be a Dedekind domain, $K={\rm Frac}(A)$, $L$ a finite separable extension of $K$, and $B$ the integral closure of $A$ in $L$. 
\begin{enumerate}
    \item[a.]  Let $\mathfrak{P}$ be a non-zero prime of $B$ such that the corresponding residue extension is separable and let $e_{\mathfrak{P}}$ denote the ramification index of $\mathfrak{P}$. Then the exponent of $\mathfrak{P}$ in the different $\mathfrak{D}_{B/A}$ is greater than or equal to $e_{\mathfrak{P}}-1$ with equality precisely when $\mathfrak{P}$ is tamely ramified. 
    \item[b.] Suppose for each prime $\mathfrak{P}$ of $B$, the corresponding residue extension is separable. The annihilator of the $B$-module $\Omega_{B/A}^{1}$ of Kahler differentials is equal to $\mathfrak{D}_{B/A}$.
\end{enumerate}
\end{prop}
\begin{proof}
(a) is \cite[III, §6, Prop 13]{[Ser79]} while (b) is \cite[III, §7, Prop. 14]{[Ser79]}.
\end{proof}

\begin{lem}\label{Localize of Different}
Suppose the Dedekind different is defined for $A\rightarrow B$. Let $S\subset A$ be a multiplicatively closed subset such that the Dedekind different is defined for $S^{-1}A\rightarrow S^{-1}B$. Then $S^{-1}\mathfrak{D}_{B/A}=\mathfrak{D}_{S^{-1}B/S^{-1}A}$ as $S^{-1}B$-modules. 
\end{lem}
\begin{proof}
First we show $S^{-1}{\cal L}_{B/A}={\cal L}_{S^{-1}B/S^{-1}A}$. By definition, 
\[
{\cal L}_{S^{-1}B/S^{-1}A}=\left\{ \frac{x}{s}\in S^{-1}L:{\rm Tr}_{S^{-1}L/K}\left(\frac{b}{s'}\frac{x}{s}\right)\in S^{-1}A\mbox{ for all }\frac{b}{s'}\in S^{-1}B\right\} .
\]
Now 
\[
{\rm Tr}_{S^{-1}L/K}\left(\frac{bx}{ss'}\right)\in S^{-1}A\mbox{ for all }\frac{b}{s'}\in S^{-1}B
\]
\[
\iff\frac{1}{ss'}{\rm Tr}_{S^{-1}L/K}(bx)\in S^{-1}A\mbox{ for all }\frac{b}{s'}\in S^{-1}B
\]
\[
\iff{\rm Tr}_{L/K}(bx)\in A\mbox{ for all }b\in B.
\]
Thus ${\cal L}_{S^{-1}B/S^{-1}A}$ can be identified with 
\[
=\left\{ \frac{x}{s}\in S^{-1}L:{\rm Tr}_{L/K}\left(bx\right)\in A\mbox{ for all }b\in B\right\} =S^{-1}{\cal L}_{B/A}.
\]
Lastly, we show $S^{-1}({\cal L}_{B/A}^{-1})=(S^{-1}{\cal L}_{B/A})^{-1}$. Let $\frac{x}{s}\in S^{-1}({\cal L}_{B/A}^{-1})$ so $x{\cal L}_{B/A}\subseteq B$. Then $\frac{x}{s}S^{-1}{\cal L}_{B/A}\subseteq S^{-1}B$ hence $S^{-1}({\cal L}_{B/A}^{-1})\subseteq(S^{-1}{\cal L}_{B/A})^{-1}$. 

For the other inclusion, we first note that ${\cal L}_{B/A}$ is finitely generated since $B$ is noetherian. Let $x_{1},\dots,x_{n}\in{\cal L}_{B/A}$ denote the generators of ${\cal L}_{B/A}$ over $B$. Let $x\in(S^{-1}{\cal L}_{B/A})^{-1}$ so $xS^{-1}{\cal L}_{B/A}\subseteq S^{-1}B$. Then for each $i=1,\dots,n$ there exists $s_{i}\in S$ such that $s_{i}xx_{i}\in B$. Let $s=\prod_{i=1}^{n}s_{i}$. Then $sxx_{i}\in B$ so $sx{\cal L}_{B/A}\subseteq B$. Therefore $x{\cal L}_{B/A}\subset S^{-1}B$ so $x\in S^{-1}({\cal L}_{B/A}^{-1})$. 

In conclusion, 
\[
S^{-1}\mathfrak{D}_{B/A}=S^{-1}({\cal L}_{B/A}^{-1})=(S^{-1}{\cal L}_{B/A})^{-1}=({\cal L}_{S^{-1}B/S^{-1}A})^{-1}=\mathfrak{D}_{S^{-1}B/S^{-1}A}.\qedhere
\]
\end{proof}

Let $f:Y\rightarrow X$ be a proper\footnote{More generally, one can define the different of a locally quasi-finite morphism of locally noetherian schemes, as in \cite[Tag0BTC]{[Stacks]}.} morphism of locally noetherian schemes. According to \cite[Tag0BVG]{[Stacks]}, the relative dualizing sheaf $\Omega_{f}$ is the unique coherent ${\cal O}_{Y}$-module such that for every pair of affine opens ${\rm Spec}(B)\subset Y$ and ${\rm Spec}(A)\subset X$ with $f({\rm Spec}(B))\subset{\rm Spec}(A)$, we have a canonical isomorphism 
\[
H^{0}({\rm Spec}(B),\Omega_{f})\simeq{\rm Hom}_{A}(B,A).
\]
If we further assume $f$ is flat, then by \cite[Tag0BVJ]{[Stacks]}, there exists a global section $\tau_{Y/X}\in H^{0}(Y,\Omega_{f})$ such that whenever $A\rightarrow B$ is finite, $\tau_{Y/X}|_{{\rm Spec}(B)}$ is identified with ${\rm Tr}_{B/A}$ under the isomorphism. 

\begin{defn}
Let $f:Y\rightarrow X$ be a flat, proper morphism of noetherian schemes. The \textbf{different} $\mathfrak{D}_{f}$ is the annihilator of the cokernel
\[
{\rm Coker}({\cal O}_{Y}\overset{\tau_{Y/X}}{\longrightarrow}\Omega_{f})
\]
which is a coherent ideal sheaf $\mathfrak{D}_{f}\subset{\cal O}_{Y}$.
\end{defn}

By \cite[Tag0BW5]{[Stacks]}, we have $\mathfrak{D}_{f}|_{{\rm Spec}(B)}=(\mathfrak{D}_{B/A})^{\sim}$ where $(\mathfrak{D}_{B/A})^{\sim}$ is the quasi-coherent sheaf induced by the $B$-module $\mathfrak{D}_{B/A}$.

\begin{lem}\label{Stalks of Different}
Let $f:Y\rightarrow X$ be a proper
morphism of noetherian schemes. Let ${\rm Spec}(B)\subset Y$ and ${\rm Spec}(A)\subset X$ be affine open subschemes such that $f({\rm Spec}(B))\subset{\rm Spec}(A)$. Let $x\in{\rm Spec}(B)$ and suppose the Dedekind different is defined for ${\cal O}_{X,f(x)}\rightarrow{\cal O}_{Y,x}$. Then $\mathfrak{D}_{f,x}\simeq\mathfrak{D}_{{\cal O}_{Y,x}/{\cal O}_{X,f(x)}}$.
\end{lem}
\begin{proof}
Let $\mathfrak{p}$ be the prime of $B$ corresponding to $x$ and $\mathfrak{q}$ be the prime of $A$ corresponding to $f(x)$. We have 
\[
\mathfrak{D}_{f,x}\simeq(\mathfrak{D}_{B/A})_{\mathfrak{p}}\simeq\mathfrak{D}_{B_{\mathfrak{p}}/A_{\mathfrak{q}}}=\mathfrak{D}_{{\cal O}_{Y,x}/{\cal O}_{X,f(x)}}.\qedhere
\]
\end{proof}

Let ${\cal R}\subseteq X$ denote the closed subscheme associated to $\mathfrak{D}_{f}$ and let $\left[{\cal R}\right]$ denote the Weil divisor associated to ${\cal R}$. The following is \cite[Tag0BWA]{[Stacks]}.

\begin{prop}\label{RDS and Different Divisor}
Let $f:Y\rightarrow X$ be a proper morphism of noetherian schemes. If $\Omega_{f}$ is invertible and $f$ is \'{e}tale at the associated points of $Y$, then ${\cal R}$ is an effective Cartier divisor and $\Omega_{f}\simeq{\cal O}_{Y}({\cal R})$.
\end{prop}

Explicitly, the Weil divisor associated to ${\cal R}$ is 
\[
\left[{\cal R}\right]=\sum_{\substack{x\in Y\\
\dim{\cal O}_{Y,x}=1
}
}{\rm mult}_{x}({\cal R})\cdot\overline{\left\{ x\right\} }.
\]
By definition, ${\rm mult}_{x}({\cal R})={\rm length}_{{\cal O}_{Y,x}}({\cal O}_{Y,x}/{\cal R}_{x}{\cal O}_{Y,x})$ and ${\cal R}_{x}=\mathfrak{D}_{f,x}$. We will now show that the usual projection morphism ${\rm pr}:\mathfrak{X}\rightarrow\mathfrak{X}(N)_{/R}$ satisfies the necessary properties to use the above results on the different. By \cite[{}5.5.1(1)]{[KM85]}, the map ${\rm pr}$ is finite and flat so the different $\mathfrak{D}_{{\rm pr}}$ is defined. 

Suppose $x\in\mathfrak{X}$ is a codimension 1 point. Let 
\[
{\rm pr}_{x}:{\cal O}_{\mathfrak{X}(N),{\rm pr}(x)}\rightarrow{\cal O}_{\mathfrak{X},x}
\]
denote the induced map on stalks. Since $x$ is of codimension 1 and ${\rm pr}$ is flat, we have ${\rm pr}(x)$ is also of codimension 1. Consequently, both ${\cal O}_{\mathfrak{X},x}$ and ${\cal O}_{\mathfrak{X}(N),{\rm pr}(x)}$
are DVRs.

\begin{prop}\label{Props of prx}
Let $x\in\mathfrak{X}$ be a codimension 1 point. The induced map on stalks ${\rm pr}_{x}$ is a finite ring map. Furthermore, the induced map on fraction fields 
\[
K={\rm Frac}({\cal O}_{\mathfrak{X}(N),{\rm pr}(x)})\rightarrow L={\rm Frac}({\cal O}_{\mathfrak{X}(Np^{r}),x})
\]
is finite separable and ${\cal O}_{\mathfrak{X},x}$ is the integral closure of ${\cal O}_{\mathfrak{X}(N),{\rm pr}(x)}$ in $L$.
\end{prop}
\begin{proof}
Let ${\rm Spec}(B)\subset\mathfrak{X}(N)$ be an affine open containing ${\rm pr}(x)$ and ${\rm Spec}(A)\subset\mathfrak{X}$ be an affine open containing $x$ such that ${\rm pr}({\rm Spec}(A))\subset{\rm Spec}(B)$. Since ${\rm pr}$ is finite, it is also integral so the induced map $A\rightarrow B$ is integral. By \cite[Tag034K]{[Stacks]}, the induced map on localizations remains integral so ${\rm pr}_{x}$ is integral. Since $\mathfrak{X}(N)$ is an integral scheme, we have inclusions
\[
A\hookrightarrow{\cal O}_{\mathfrak{X}(N),{\rm pr}(x)}\hookrightarrow{\rm Frac}(A).
\]
Therefore $K={\rm Frac}(A)$. Similarly, we conclude $L={\rm Frac}(B).$ Both $\mathfrak{X}(N)$ and $\mathfrak{X}$ are normal schemes so ${\cal O}_{\mathfrak{X}(N),{\rm pr}(x)}$ (resp. ${\cal O}_{\mathfrak{X},x}$) is integrally closed in $K$ (resp. $L$). Having established 
\[
{\rm pr}_{x}:{\cal O}_{\mathfrak{X}(N),{\rm pr}(x)}\rightarrow{\cal O}_{\mathfrak{X},x}
\]
is integral, the integral closure of ${\cal O}_{\mathfrak{X}(N),{\rm pr}(x)}$ in $L$ is precisely ${\cal O}_{\mathfrak{X},x}$. 

Since ${\rm pr}$ is a finite map between two integral schemes, the extension of function fields 
\[
K(\mathfrak{X}(N))\rightarrow K(\mathfrak{X})
\]
is a finite extension of characteristic zero fields, and hence separable. Note that $K(\mathfrak{X}(N))=K$ and $K(\mathfrak{X})=L$ so $L/K$ is a finite separable extension. By \cite[I, §IV, Prop. 8]{[Ser79]}, we can conclude ${\rm pr}_{x}$ is finite.
\end{proof}

Since $\mathfrak{X}$ is integral, it's only associated point is its unique generic point. In Proposition \ref{Props of prx}, we deduced the map on function fields $K(\mathfrak{X})\rightarrow K(\mathfrak{X}(N))$ is finite separable so ${\rm pr}$ is \'{e}tale at the generic point of $\mathfrak{X}$. By Proposition \ref{RDS and Different Divisor}, we have
\begin{equation}
\Omega_{\mathfrak{X}/\mathfrak{X}(N)}\simeq{\cal O}_{\mathfrak{X}}([{\cal R}])={\cal O}_{\mathfrak{X}}\left(\sum_{\substack{x\in\mathfrak{X}\\
\dim{\cal O}_{\mathfrak{X},x}=1
}
}{\rm mult}_{x}({\cal R})\cdot\overline{\left\{ x\right\} }\right)\label{eq: RDS iso and diff}
\end{equation}
Since ${\cal O}_{\mathfrak{X},x}$ is a DVR in this case, ${\rm mult}_{x}({\cal R})$ is equal to the valuation of ${\cal R}_{x}=\mathfrak{D}_{{\rm pr},x}$ in ${\cal O}_{\mathfrak{X},x}$. By Corollary \ref{Stalks of Different} and Proposition \ref{Props of prx}, $\mathfrak{D}_{{\rm pr},x}=\mathfrak{D}_{{\cal O}_{\mathfrak{X},x}/{\cal O}_{\mathfrak{X}(N),{\rm pr}(x)}}$. Also, by Proposition \ref{AKLB Different Props}(b), we can identify $\mathfrak{D}_{{\rm pr},x}$ with 
\[
{\rm Ann}_{{\cal O}_{\mathfrak{X},x}}(\Omega_{{\cal O}_{\mathfrak{X},x}/{\cal O}_{\mathfrak{X}(N),{\rm pr}(x)}}^{1}).
\]

Let $d_{x}$ denote the valuation of the different $\mathfrak{D}_{{\rm pr},x}$ in the DVR ${\cal O}_{\mathfrak{X},x}$. Recall by \cite[{}8.3.4]{[Liu02]} the codimension 1 points of $\mathfrak{X}$ are precisely the closed points of the generic fiber, and the generic points of the special fiber. We will compute $d_{x}$ when $x$ is a closed point of the generic fiber. If $x$ is a generic point of the special fiber, then we show all $d_{x}$ contributions are the same which taken together contribute nothing to the different of ${\rm pr}$. However, one can compute the $d_{x}$ explicitly using strict Henselizations. 

\subsection{\texorpdfstring{Computing $d_{x}$ for the closed points of the generic fiber}{}}

The generic fiber of $\mathfrak{X}$ is open so it suffices to compute $d_{x}$ over $\mathbb{Q}_{p}(\zeta_{Np^{r}})$ i.e. for the map 
\[
{\rm pr}_{x}:{\cal O}_{\mathfrak{X}(N)_{/\mathbb{Q}_{p}(\zeta_{Np^{r}})},{\rm pr}(x)}\rightarrow{\cal O}_{\mathfrak{X},x}.
\]
First we show the value of $d_{x}$ does not change after base changing $\mathbb{Q}_{p}(\zeta_{Np^{r}})$. In particular, we can compute the value of $d_{x}$ over $\mathbb{C}$ and use the classical theory of modular curves as compact Riemann surfaces. 

\begin{lem}
Let $\pi:X\rightarrow Y$ be a finite type morphism of normal curves over a field $K$ and let $L$ be a field extension of $K$. Let $p:X_{L}\rightarrow X$ denote the usual projection morphism from base change and let $x\in X_{L}$. Then we have $d_{x}=d_{p(x)}$. 
\end{lem}
\begin{proof}
We can equate $d_{x}$ with the valuation of the annihilator ideal of $\Omega_{X_{L}/Y_{L},x}^{1}$ in ${\cal O}_{X_{L},x}$. Since Kahler differentials are compatible with base change, we have 
\[
\Omega_{X_{L}/Y_{L},x}^{1}\simeq(p^{*}\Omega_{X/Y}^{1})_{x}\simeq\Omega_{X/Y,p(x)}^{1}\otimes_{{\cal O}_{X,p(x)}}{\cal O}_{X_{L},x}.
\]
Since $K\rightarrow L$ is flat, the map ${\cal O}_{X,p(x)}\rightarrow{\cal O}_{X_{L},x}$ is flat. Furthermore, $\Omega_{X/Y,p(x)}^{1}$ is finite over ${\cal O}_{X,p(x)}$ since $\pi$ is finite type. Thus, by Lemma \ref{Ann Base Change}, we have
\begin{align*}
{\rm Ann}_{{\cal O}_{X_{L}},x}(\Omega_{X_{L}/Y_{L},x}^{1}) & = {\rm Ann}_{{\cal O}_{X_{L}},x}(\Omega_{X/Y,p(x)}^{1}\otimes_{{\cal O}_{X,p(x)}}{\cal O}_{X_{L},x})\\
 & = {\rm Ann}_{{\cal O}_{X},p(x)}(\Omega_{X/Y,x}^{1}){\cal O}_{X_{L},x}.
\end{align*}
Hence $d_{x}=d_{p(x)}$. 
\end{proof}

Thus we can compute $d_{x}$ in the situation that our modular curves are over $\mathbb{C}$. Let $X(M)$ denote the modular curve over $\mathbb{C}$ of level $\Gamma(M)$.  In this situation, ${\cal O}_{X(Np^{r}),x}$ is tamely ramified over ${\cal O}_{X(N),{\rm pr}(x)}$ so by Proposition \ref{AKLB Different Props}(a), we have $d_{x}=e_{x}-1$ where $e_{x}$ is the ramification index of $x$. We will investigate the ramification index of all points in $X(Np^{r})$ under the usual projection map ${\rm pr}:X(Np^{r})\rightarrow X(N)$. 

Let $\Gamma\leq{\rm SL}_{2}(\mathbb{Z})$ be a congruence subgroup. We also let ${\cal H}$ denote the (complex) upper half plane and ${\cal H}^{*}={\cal H}\cup\mathbb{P}^{1}(\mathbb{Q})$. As we will see shortly, investigating the ramification of a point $z\in{\cal H}^{*}$ amounts to understanding the stabilizer group 
\[
\Gamma_{z}=\left\{ \gamma\in\Gamma:\gamma\cdot z=z\right\} .
\]
The points of $z\in{\cal H}^{*}$ can be classified depending on which elements of $\Gamma$ fix $z$. Note that 
\[
\left(\begin{array}{cc}
a & b\\
c & d
\end{array}\right)\cdot z=z\iff cz^{2}+(d-a)z-b=0
\]
for matrices not equal to $\pm I$. Hence the fixed points of a given $\gamma\in\Gamma$ are either conjugate complex numbers, a single real number, or two distinct real numbers. Based on this fact, the following definition from \cite[§1.3, p.27]{[Miy76]}.

\begin{defn}
We say $z\in{\cal H}^{*}$ is a/an
\begin{itemize}
    \item \textbf{elliptic point of $\Gamma$} if there exists $\gamma\in\Gamma_{z}$ such that $\left|{\rm tr}(\gamma)\right|<2$, or equivalently $\gamma$ has two distinct fixed points $z$ and $\bar{z}$. 
    \item \textbf{cusp} \textbf{of $\Gamma$} if there exists $\gamma\in\Gamma_{z}$ such that $\left|{\rm tr}(\gamma)\right|=2$, or equivalently $\gamma$ has a unique real fixed point $z$.
    \item \textbf{hyperbolic point of $\Gamma$} if there exists $\gamma\in\Gamma_{z}$ such that $\left|{\rm tr}(\gamma)\right|>2$, or equivalently $\gamma$ has two distinct real fixed points.
    \item \textbf{ordinary point of $\Gamma$ }if $\Gamma$ does not fix $z$, excluding $\pm I$.
\end{itemize}
\end{defn}

Since the matrices we consider are in ${\rm SL}_{2}(\mathbb{Z})$, we will not have any hyperbolic points appearing.  We denote $\bar{\Gamma}=\Gamma/\Gamma\cap\left\{ \pm I\right\} $. The following is \cite[Theorem 1.5.4]{[Miy76]} which describes the
stabilizer groups $\Gamma_{z}$.

\begin{thm}\leavevmode
\begin{enumerate}
    \item[a.] If $z\in{\cal H}$ is an elliptic point of $\Gamma$, then $\Gamma_{z}$ is a finite cyclic group.
    \item[b.]  If $z\in\mathbb{Q}\cup\{\infty\}$ is a cusp of $\Gamma$, then $\bar{\Gamma}_{z}\simeq\mathbb{Z}$.
\end{enumerate}
\end{thm}

\noindent The following proposition, which is \cite[Prop. 1.37]{[Shi71]}, relates the ramification index with the index of stabilizer groups.

\begin{prop}\label{Ram Index Stabilizer}
Let $\Gamma'\leq\Gamma$ be a finite index subgroup and consider the projection $\pi:\Gamma'\backslash{\cal H}^{*}\rightarrow\Gamma\backslash{\cal H}^{*}$. The ramification index $e_{z}$ of a point $z\in\Gamma'\backslash{\cal H}^{*}$ under $\pi$ is equal to 
\[
e_{z}=\left[\bar{\Gamma}_{z}:\bar{\Gamma'}_{z}\right].
\]
\end{prop}

We are now ready to compute the ramification index under our map ${\rm pr}$.

\begin{prop}\label{Ramification of pr over C}
Let ${\rm pr}:X(Np^{r})\rightarrow X(N)$ denote the usual projection map and let $z\in X(Np^{r})$. We have
\[
e_{z}=\begin{cases}
p^{r} & \mbox{if }z\mbox{ is a cusp}\\
1 & \mbox{otherwise}
\end{cases}.
\]
\end{prop}
\begin{proof}
By \cite[Prop. 1.39]{[Shi71]}, the congruence subgroup $\Gamma(M)$ has no elliptic elements for any $M>1$. Consequently, the ordinary points are precisely the points of $\Gamma(M)\backslash{\cal H}$. We will now apply Proposition \ref{Ram Index Stabilizer} in the case $\Gamma'=\Gamma(Np^{r})$ and $\Gamma=\Gamma(N)$. 

Note that the image of a cusp under ${\rm pr}$ is again a cusp and similarly for ordinary points. If $z\in X(Np^{r})$ is ordinary, then its stabilizer group is trivial so $e_{z}=1$. If $z$ is a cusp, then there exists $\gamma\in{\rm SL}_{2}(\mathbb{Z})$ such that $\gamma\cdot z=\infty$. Therefore
\[
\overline{\Gamma(N)}_{z}/\overline{\Gamma(Np^{r})}_{z}\simeq\gamma\overline{\Gamma(N)}_{z}\gamma^{-1}/\gamma\overline{\Gamma(Np^{r})}_{z}\gamma^{-1}=\overline{\Gamma(N)}_{\infty}/\overline{\Gamma(Np^{r})}_{\infty}
\]
so it suffices to compute $e_{\infty}$. 

Let $A(d)=\left(\begin{array}{cc}
1 & d\\
0 & 1
\end{array}\right).$ According to \cite[bottom p. 22]{[Shi71]}, we have for $M>1$ that $\overline{\Gamma(M)}_{\infty}=\left\langle A(d)\right\rangle .$ Furthermore, $-I\in\Gamma(M)$ if and only if $-1\equiv1$ modulo $M$, or equivalently $M=2$. Since $M\geq3$ in our situation, we always have $\overline{\Gamma(M)}_{\infty}=\Gamma(M)_{\infty}$. Note that for any $m\geq0$, we have $A(d)^{m}=A(md)$. Therefore the order of $A(d)$ in  $\overline{\Gamma(N)}_{\infty}/\overline{\Gamma(Np^{r})}_{\infty}$ is equal to $p^{r}$ so 
\[
e_{\infty}=[\overline{\Gamma(N)}_{\infty}:\overline{\Gamma(Np^{r})}_{\infty}]=p^{r}.\qedhere
\]
\end{proof}

\subsection{\label{sub:Computing deg(omega|lambda)}\texorpdfstring{Relating the modular sheaf with $\Omega_{\mathfrak{X}/R}$}{}}

Let $x\in\mathfrak{X}$ be a generic point of the special fiber. Recall the value of $d_{x}$ is equal to the valuation of the different ideal corresponding to the induced map on stalks
\[
{\rm pr}_{x}:{\cal O}_{\mathfrak{X}(N),{\rm pr}(x)}\rightarrow{\cal O}_{\mathfrak{X},x}.
\]
We will use the discussion in Section \ref{sub:cusp=000026specialfiber} to provide an explicit description of ${\rm pr}_{x}$. Recall in Theorem \ref{Iso Ig and ExIg} we have a commutative diagram 
\[
\xymatrix{{\rm Ig}(p^{r},N)\ar@{->}[r]^{\simeq}\ar@{->}[d]^{\rho} & {\rm ExIg}(p^{r},i,N)\ar@{->}[d]^{\rho'}\\
\mathfrak{X}(N)_{\mathbb{F}_{q}}^{(\sigma^{-i})}\ar@{->}[r]^{F^{i}} & \mathfrak{X}(N)_{\mathbb{F}_{q}}
}
\]
where $\rho$ and $\rho'$ are the usual projection maps. By Theorem \ref{Special Fiber Comps}, the restriction of ${\rm pr}$ to any irreducible component of $\bar{\mathfrak{X}}$ is the map $\rho'$.
We get a commutative diagram: 

\begin{gather}
\begin{gathered}\xymatrix{{\cal O}_{\mathfrak{X}(N)_{/\mathbb{F}_{q}},{\rm pr}(x)}\ar@{->}[r]^{F^{r}}\ar@{->}[d]^{{\rm pr}_{x}} & {\cal O}_{\mathfrak{X}(N)_{/\mathbb{F}_{q}}^{(p^{-r})},{\rm pr}(x)}\ar@{->}[d]^{{\rm \rho}_{x}}\\
{\cal O}_{{\rm ExIg}(p^{r},r,N),x}\ar@{->}[r]^{\simeq} & {\cal O}_{{\rm Ig}(p^{r},N),x}
}
\end{gathered}
\label{eq: Diag of Res Fields-1}
\end{gather}
In particular, the map ${\rm pr}_{x}$ is the same for each generic point $x$ of the special fiber. Hence, the value of $d_{x}$ is independent of $x$ in this situation.

\begin{lem}\label{RDS of X(Np^r)/X(N)}
We have 
\[
\Omega_{\mathfrak{X}/\mathfrak{X}(N)}\simeq{\cal O}_{\mathfrak{X}}((p^{r}-1)\mathfrak{C}(Np^{r})).
\]
\end{lem}
\begin{proof}
Recall from equation (\ref{eq: RDS iso and diff}) and the paragraphs proceeding it, we have 
\[
\Omega_{\mathfrak{X}/\mathfrak{X}(N)}\simeq{\cal O}_{\mathfrak{X}}\left(\sum_{x\in\mathfrak{X}{}_{/\mathbb{Q}_{p}(\zeta_{Np^{r}})}}d_{x}\cdot\overline{\left\{ x\right\} }+\sum_{x\in\mathfrak{X}_{/\mathbb{F}_{q}}}d_{x}\cdot\overline{\left\{ x\right\} }\right)
\]
where the first sum is over closed points of the generic fiber and the second sum is over the generic points of the irreducible components of the special fiber. By Proposition \ref{Ramification of pr over C}, we have 
\[
\Omega_{\mathfrak{X}/\mathfrak{X}(N)}\simeq{\cal O}_{\mathfrak{X}}\left(\sum_{x\in\mathfrak{C}(Np^{r})}(p^{r}-1)\cdot\overline{\left\{ x\right\} }+\sum_{x\in\mathfrak{X}_{/\mathbb{F}_{q}}}d_{x}\cdot\overline{\left\{ x\right\} }\right)
\]
where the first sum is over all the cusps in the generic fiber. As discussed above, the $d_{x}$ values appearing in the second sum are independent of $x$. Using the fact that the special fiber is reduced, we have 
\[
\Omega_{\mathfrak{X}/\mathfrak{X}(N)}\simeq{\cal O}_{\mathfrak{X}}((p^{r}-1)\mathfrak{C}(Np^{r})+d_{x}\mathfrak{X}_{/\mathbb{F}_{q}}).
\]
Note that $\mathfrak{X}_{/\mathbb{F}_{q}}$, when viewed as a divisor, is principal (see Proposition \ref{Special Fiber as Div}). Thus 
\[
\Omega_{\mathfrak{X}/\mathfrak{X}(N)}\simeq{\cal O}_{\mathfrak{X}}((p^{r}-1)\mathfrak{C}(Np^{r}))
\]
as desired.
\end{proof}

Recall the map from (\ref{eq: pr map}):
\[
\mathfrak{X}\overset{{\rm pr}}{\longrightarrow}\mathfrak{X}(N)_{/R}\overset{g}{\longrightarrow}{\rm Spec}(R)
\]
Again, for convenience, we let $\mathfrak{X}(N)$ denote the base change $\mathfrak{X}(N)_{/R}$. To compute $\deg(\underline{\omega}_{\mathfrak{X}}^{\otimes2}|_{\Lambda})$, our strategy is to first prove that $\underline{\omega}_{\mathfrak{X}}^{\otimes2}\simeq\Omega_{\mathfrak{X}/R}(\mathfrak{C}(Np^{r}))$ using the isomorphism in (\ref{eq:Adj + KS}):
\[
\Omega_{\mathfrak{X}/R}\simeq\Omega_{\mathfrak{X}/\mathfrak{X}(N)}\otimes_{{\cal O}_{\mathfrak{X}}}{\rm pr}^{*}\underline{\omega}_{\mathfrak{X}(N)}^{\otimes2}(-\mathfrak{C}(N)).
\]
Our next step is to provide a better description of ${\rm pr}^{*}\underline{\omega}_{\mathfrak{X}(N)}^{\otimes2}(-\mathfrak{C}(N))$. 

\begin{lem}\label{Pullback by pr}
We have an isomorphism of sheaves on $\mathfrak{X}(Np^{r})$
\[
{\rm pr}^{*}\underline{\omega}_{\mathfrak{X}(N)}^{\otimes2}(-\mathfrak{C}(N))\simeq\underline{\omega}_{\mathfrak{X}}^{\otimes2}\otimes_{{\cal O}_{\mathfrak{X}}}{\cal O}_{\mathfrak{X}}(-p^{r}\mathfrak{C}(Np^{r})).
\]
\end{lem}
\begin{proof}
We have 
\[
{\rm pr}^{*}\underline{\omega}_{\mathfrak{X}(N)}^{\otimes2}(-\mathfrak{C}(N))\simeq{\rm pr}^{*}\underline{\omega}_{\mathfrak{X}(N)}^{\otimes2}\otimes_{{\cal O}_{\mathfrak{X}(Np^{r})}}{\rm pr}^{*}{\cal O}_{\mathfrak{X}(N)}(-\mathfrak{C}(N)).
\]
By Proposition \ref{Base change} and Proposition \ref{Pullback of Modular Sheaf}, we have 
\[
{\rm pr}^{*}\underline{\omega}_{\mathfrak{X}(N)}^{\otimes2}=\underline{\omega}_{\mathfrak{X}(Np^{r})}^{\otimes2}.
\]
Since ${\rm pr}$ maps $\mathfrak{C}(Np^{r})$ onto $\mathfrak{C}(N)$ and is ramified at the cusps of degree $p^{r}$ by\emph{ }Proposition \ref{Ramification of pr over C}, we have 
\[
{\rm pr}^{*}{\cal O}_{\mathfrak{X}(N)}(-\mathfrak{C}(N))\simeq{\cal O}_{\mathfrak{X}(Np^{r})}(-p^{r}\mathfrak{C}(Np^{r})).
\]
Thus 
\[
{\rm pr}^{*}\underline{\omega}_{\mathfrak{X}(N)}^{\otimes2}(-\mathfrak{C}(N))\simeq\underline{\omega}_{\mathfrak{X}}^{\otimes2}\otimes_{{\cal O}_{\mathfrak{X}}}{\cal O}_{\mathfrak{X}}(-p^{r}\mathfrak{C}(Np^{r}))
\]
as desired. 
\end{proof}

\begin{thm}\label{iso of modular sheaf}
We have $\underline{\omega}_{\mathfrak{X}}^{\otimes2}\simeq\Omega_{\mathfrak{X}/R}(\mathfrak{C}(Np^{r})).$ 
\end{thm}
\begin{proof}
Recall the isomorphism in (\ref{eq:Adj + KS}) states
\[
\Omega_{\mathfrak{X}/R}\simeq\Omega_{\mathfrak{X}/\mathfrak{X}(N)}\otimes_{{\cal O}_{\mathfrak{X}}}{\rm pr}^{*}\underline{\omega}_{\mathfrak{X}(N)}^{\otimes2}(-\mathfrak{C}(N)).
\]
By Lemma \ref{RDS of X(Np^r)/X(N)} and Lemma \ref{Pullback by pr}, we have 
\[
\Omega_{\mathfrak{X}/R}\simeq{\cal O}_{\mathfrak{X}}((p^{r}-1)\mathfrak{C}(Np^{r}))\otimes\underline{\omega}_{\mathfrak{X}}^{\otimes2}\otimes{\cal O}_{\mathfrak{X}}(-p^{r}\mathfrak{C}(Np^{r}))\simeq\underline{\omega}_{\mathfrak{X}}^{\otimes2}(-\mathfrak{C}(Np^{r})).\qedhere
\]
\end{proof}

\begin{cor}\label{Main Formula for Deg Omega}
Let $\Lambda$ be an irreducible component of the special fiber of $\mathfrak{X}$. We have 
\[
\deg(\underline{\omega}_{\mathfrak{X}}^{\otimes2}|_{\Lambda})=\deg(\Omega_{\mathfrak{X}/R}|_{\Lambda})+\deg(\mathfrak{C}(Np^{r})|_{\Lambda}).
\]
\end{cor}

The following will allow us to calculate $\deg(\mathfrak{C}(Np^{r})|_{\Lambda})$.

\begin{lem}\label{Deg of Cusps Lem}
$\deg(\mathfrak{C}(Np^{r})|_{\Lambda})$ is equal to the number of cusps $\overline{\left\{ x\right\} }$ which intersect $\Lambda$.
\end{lem}
\begin{proof}
Viewed as a divisor, $\mathfrak{C}(Np^{r})$ is equal to the closure of all the cusps of the generic fiber of $\mathfrak{X}$ by Proposition \ref{Description of C(N) as Divisor}. By Theorem \ref{Props of Int Num}(d),
\[
\deg(\mathfrak{C}(Np^{r})|_{\Lambda})=\mathfrak{C}(Np^{r}).\Lambda=\sum_{x}\overline{\left\{ x\right\} }.\Lambda
\]
where the sum is indexed over $x\in\mathfrak{C}(\mathfrak{X}{}_{/\mathbb{Q}_{p}(\zeta_{Np^{r}})})$ \emph{i.e. }over the cusps of the generic fiber. 

All the cusps are rational, so by Corollary \ref{Int of Ratl Pt}, $\overline{\left\{ x\right\} }$ intersects at a single irreducible component of the special fiber. Therefore
\[
\overline{\left\{ x\right\} }.\Lambda=\begin{cases}
1 & \mbox{if }\overline{\left\{ x\right\} }\cap\Lambda\ne\emptyset\\
0 & \mbox{otherwise}
\end{cases}
\]
so $\deg(\mathfrak{C}(Np^{r})|_{\Lambda})$ is equal to the number of cusps which intersect $\Lambda$ as desired.
\end{proof}

For convenience, we let $C(Np^{r})=\mathfrak{C}(\mathfrak{X}{}_{/\mathbb{Q}_{p}(\zeta_{Np^{r}}})$, the set of all cusps of the generic fiber. Using Lemma \ref{Deg of Cusps Lem} and Proposition \ref{num of irr comps int cusp}, we obtain the following:

\begin{prop}\label{Deg of Cusp Restricted}
We have 
\[
\deg(\mathfrak{C}(Np^{r})|_{\Lambda})=\varphi(p^{r})\#C(N).
\]
\end{prop}

\noindent To compute $\deg(\Omega_{\mathfrak{X}(Np^{r})/R}|_{\Lambda})$, we will use the following result in \cite[Theorem 9.1.37]{[Liu02]}.

\begin{thm}\label{Adj for Irr Comps}
Let $X\rightarrow S$ be a regular fibered surface, $s\in S$ a closed point, and $E\in{\rm Div}_{s}(X)$ such that $0<E\leq X_{s}$. Then we have 
\[
\Omega_{E/k(s)}\simeq({\cal O}_{X}(E)\otimes\Omega_{X/S})|_{E}.
\]
\end{thm}

\begin{cor}\label{Deg RDS Restricted}
We have 
\[
\deg(\Omega_{\mathfrak{X}/R}|_{\Lambda})=p^{r}\varphi(p^{r})\frac{\#{\rm SL}_{2}(\mathbb{Z}/N\mathbb{Z})}{24}-\varphi(p^{r})\#C(N)+\deg{\rm S}(N)\cdot p^{2r-1}
\]
\end{cor}
\begin{proof}
Applying Theorem \ref{Adj for Irr Comps} to $\mathfrak{X}\rightarrow{\rm Spec}(R)$ and $E=\Lambda$, we get 
\[
\Omega_{\Lambda/\mathbb{F}_{q}}\simeq({\cal O}_{\mathfrak{X}}(\Lambda)\otimes\Omega_{\mathfrak{X}/R})|_{\Lambda}.
\]
Therefore
\begin{align*}
\deg(\Omega_{\Lambda/\mathbb{F}_{q}}) & = \deg(({\cal O}_{\mathfrak{X}}(\Lambda)\otimes\Omega_{\mathfrak{X}/R})|_{\Lambda})\\
 & = \deg({\cal O}_{\mathfrak{X}}(\Lambda)|_{\Lambda})+\deg(\Omega_{\mathfrak{X}/R}|_{\Lambda})\\
 & = \Lambda.\Lambda+\deg(\Omega_{\mathfrak{X}/R}|_{\Lambda}).
\end{align*}
Thus 
\[
\deg(\Omega_{\mathfrak{X}/R}|_{\Lambda})=\deg(\Omega_{\Lambda/\mathbb{F}_{q}})-\Lambda.\Lambda.
\]

By \cite[Corollary 7.3.31(a)]{[Liu02]}, 
\[
\deg(\Omega_{\Lambda/\mathbb{F}_{q}})=2\rho_{a}(\Lambda)-2
\]
where $\rho_{a}(\Lambda)$ is the arithmetic genus of $\Lambda$. By \cite[{}12.9.4]{[KM85]}
along with \cite[Corollary 10.13.12]{[KM85]}, we have
\begin{align*}
2\rho_{a}(\Lambda) & = 2\rho_{a}({\rm Ig}(p^{r},N))\\
 & = p^{r}\varphi(p^{r})\frac{\#{\rm SL}_{2}(\mathbb{Z}/N\mathbb{Z})}{24}+2-\varphi(p^{r})\#C(N)
\end{align*}
We have already shown $\Lambda.\Lambda=-\deg{\rm S}(N)\cdot p^{2r-1}$
in Proposition \ref{self-int num}. Thus
\begin{align*}
\deg(\Omega_{\mathfrak{X}/R}|_{\Lambda}) & = (2\rho_{a}(\Lambda)-2)-\Lambda.\Lambda\\
 & = p^{r}\varphi(p^{r})\frac{\#{\rm SL}_{2}(\mathbb{Z}/N\mathbb{Z})}{24}-\varphi(p^{r})\#C(N)+\deg{\rm S}(N)\cdot p^{2r-1}.
\end{align*}
as desired.
\end{proof}

Recall that the space of cusp forms, by definition, are the global sections of $\underline{\omega}_{\mathfrak{X}}^{\otimes2}(-\mathfrak{C}(Np^{r}))$. 

\begin{thm}\label{Deg of Res Modular Sheaf }
Let $k\geq1$ be an integer. We
have 
\[
\deg(\underline{\omega}_{\mathfrak{X}}^{\otimes2k}|_{\Lambda})=k\cdot\#{\rm SL}_{2}(\mathbb{Z}/N\mathbb{Z})\left[\frac{1}{12}(p-1)p^{2r-1}\right]
\]
and 
\[
\deg(\underline{\omega}_{\mathfrak{X}}^{\otimes2k}(-\mathfrak{C}(Np^{r}))|_{\Lambda})=k\cdot\#{\rm SL}_{2}(\mathbb{Z}/N\mathbb{Z})(p-1)\left[\frac{p^{2r-1}}{12}-\frac{p^{r-1}}{2N}\right].
\]
\end{thm}
\begin{proof}
We first focus on the case $k=1$. By Corollary \ref{Main Formula for Deg Omega}, Proposition \ref{Deg of Cusp Restricted}, and Corollary \ref{Deg RDS Restricted}, we have
\begin{align*}
\deg(\underline{\omega}_{\mathfrak{X}}^{\otimes2}|_{\Lambda}) & = p^{r}\varphi(p^{r})\frac{\#{\rm SL}_{2}(\mathbb{Z}/N\mathbb{Z})}{24}-\varphi(p^{r})\#C(N)+\deg{\rm S}(N)\cdot p^{2r-1}+\varphi(p^{r})\#C(N)\\
 & = p^{r}\varphi(p^{r})\frac{\#{\rm SL}_{2}(\mathbb{Z}/N\mathbb{Z})}{24}+\deg{\rm S}(N)\cdot p^{2r-1}
\end{align*}

\noindent By \cite[{}12.4.5]{[KM85]}, $\deg{\rm S}(N)=\frac{(p-1)\#{\rm SL}_{2}(\mathbb{Z}/N\mathbb{Z})}{24}$. Continuing, we have
\begin{alignat*}{1}
\deg(\underline{\omega}_{\mathfrak{X}}^{\otimes2}|_{\Lambda}) & =p^{r}\varphi(p^{r})\frac{\#{\rm SL}_{2}(\mathbb{Z}/N\mathbb{Z})}{24}+\frac{(p-1)\#{\rm SL}_{2}(\mathbb{Z}/N\mathbb{Z})}{24}\cdot p^{2r-1}\\
 & =\#{\rm SL}_{2}(\mathbb{Z}/N\mathbb{Z})\left[\frac{p^{r}\varphi(p^{r})}{24}+\frac{p^{2r-1}(p-1)}{24}\right]\\
 & =\#{\rm SL}_{2}(\mathbb{Z}/N\mathbb{Z})\left[\frac{1}{12}(p-1)p^{2r-1}\right].
\end{alignat*}
By Proposition \ref{Deg of Cusp Restricted} along with Lemma \ref{Num of Cusps}, we have 
\begin{alignat*}{1}
\deg(\underline{\omega}_{\mathfrak{X}}^{\otimes2}(-\mathfrak{C}(Np^{r}))|_{\Lambda}) & =\deg(\underline{\omega}_{\mathfrak{X}}^{\otimes2})-\deg(\mathfrak{C}(Np^{r})|_{\Lambda})\\
 & =\deg(\underline{\omega}_{\mathfrak{X}}^{\otimes2})-\varphi(p^{r})\#C(N)\\
 & =\#{\rm SL}_{2}(\mathbb{Z}/N\mathbb{Z})\left[\frac{1}{12}(p-1)p^{2r-1}\right]-\left[\varphi(p^{r})\frac{{\rm \#SL}_{2}(\mathbb{Z}/N\mathbb{Z})}{2N}\right]\\
 & =\#{\rm SL}_{2}(\mathbb{Z}/N\mathbb{Z})(p-1)\left[\frac{p^{2r-1}}{12}-\frac{p^{r-1}}{2N}\right].
\end{alignat*}

In general, taking tensor powers commutes with pullback of sheaves (see \cite[Tag01CD]{[Stacks]}). Hence $\underline{\omega}_{\mathfrak{X}}^{\otimes2k}|_{\Lambda}\simeq\left(\underline{\omega}_{\mathfrak{X}}^{\otimes2}|_{\Lambda}\right)^{\otimes k}$ and consequently
\[
\deg(\underline{\omega}_{\mathfrak{X}}^{\otimes2k}|_{\Lambda})=\deg\left(\left(\underline{\omega}_{\mathfrak{X}}^{\otimes2}|_{\Lambda}\right)^{\otimes k}\right)=k\cdot\deg(\underline{\omega}_{\mathfrak{X}}^{\otimes2}|_{\Lambda})
\]
which gives our desired result. 
\end{proof}

\subsection{An Upper Bound}

Recall at the end of Remark \ref{rmk on exp} we arrived at the following expression for the exponent
\[
e=\max_{\substack{1\le i\leq n\\
f\in H^{0}(\mathfrak{X},\underline{\omega}^{\otimes2})
}
}\left\{ \nu_{\Lambda_{0}}(f)-\nu_{\Lambda_{i}}(f)\right\} .
\]
At the end of Section \ref{sub:description_of_e} we established 
\[
\nu_{\Lambda}(f)-\nu_{\Lambda_{(1,0)}}(f)=\sum_{\Lambda'\ne\Lambda_{(1,0)}}\left(\deg(\underline{\omega}^{\otimes k}|_{\Lambda'})-H_{f}.\Lambda'\right)c^{\Lambda,\Lambda'}
\]
where the sum is over all irreducible components of the special fiber excluding $\Lambda_{(1,0)}$. We will now provide an upper bound for $e$. First we need to compute the sums $\sum c^{\Lambda,\Lambda'}$, where $c^{\Lambda,\Lambda'}$ is the entry of $T^{-1}$ corresponding to row label $\Lambda$ and column label $\Lambda'$ (see Section \ref{sub:Inverting M(p^r)}). 

\begin{lem}\label{P-adic Sum of Inv}
For $1\leq a'\leq p^{r}-1$, we have 
\[
\sum_{\substack{a=1\\
a\ne a'
}
}^{p^{r}-1}\nu_{p}\left(\frac{1}{a}-\frac{1}{a'}\right)=-(p^{r}-2)\nu_{p}(a').
\]
\end{lem}
\begin{proof}
We have 
\[
\sum_{\substack{a=1\\
a\ne a'
}
}^{p^{r}-1}\nu_{p}\left(\frac{1}{a}-\frac{1}{a'}\right)=\sum_{\substack{a=1\\
a\ne a'
}
}^{p^{r}-1}\nu_{p}\left(\frac{a'-a}{aa'}\right)=\sum_{\substack{a=1\\
a\ne a'
}
}^{p^{r}-1}\nu_{p}(a'-a)-\sum_{\substack{a=1\\
a\ne a'
}
}^{p^{r}-1}\nu_{p}\left(a\right)-\sum_{\substack{a=1\\
a\ne a'
}
}^{p^{r}-1}\nu_{p}\left(a'\right)
\]
\begin{equation}
=\sum_{\substack{a=1\\
a\ne a'
}
}^{p^{r}-1}\nu_{p}(a-a')-\left(\sum_{a=1}^{p^{r}-1}\nu_{p}(a)-\nu_{p}(a')\right)-(p^{r}-2)\nu_{p}(a')\label{eq:Eqn w/ Valuations}
\end{equation}
Note that
\begin{align*}
\sum_{\substack{a=1\\
a\ne a'
}
}^{p^{r}-1}\nu_{p}(a-a') & =  \sum_{a=1}^{a'-1}\nu_{p}(a-a')+\sum_{a=a'+1}^{p^{r}-1}\nu_{p}(a-a')\\
 & = \sum_{a=1}^{a'-1}\nu_{p}(a)+\sum_{a=1}^{p^{r}-1-a'}\nu_{p}(a)\\
 & = \sum_{a=-a'+1}^{-1}\nu_{p}(a)+\sum_{a=1}^{p^{r}-1-a'}\nu_{p}(a)\\
 & = \sum_{a=p^{r}-a'+1}^{p^{r}-1}\nu_{p}(a)+\sum_{a=1}^{p^{r}-a'-1}\nu_{p}(a)\\
 & = -\nu_{p}(a')+\sum_{a=1}^{p^{r}-1}\nu_{p}(a).
\end{align*}
 Continuing equation (\ref{eq:Eqn w/ Valuations}), we have
\begin{alignat*}{1}
 & =-\nu_{p}(a')+\sum_{a=1}^{p^{r}-1}\nu_{p}(a)-\left(\sum_{a=1}^{p^{r}-1}\nu_{p}(a)-\nu_{p}(a')\right)-(p^{r}-2)\nu_{p}(a')\\
 & =-(p^{r}-2)\nu_{p}(a').\tag*{\qedhere}
\end{alignat*}
\end{proof}

\begin{prop}\label{Row Sums of cij}
We have 
\[
\deg{\rm S}(N)\sum_{\Lambda'\ne\Lambda_{(1,0)}}c^{\Lambda,\Lambda'}=\begin{cases}
{\displaystyle -p^{-r}(pr-r+1)+p^{-r}(p-1)\nu_{p}(a)} & \mbox{ if }\Lambda=\Lambda_{(1,-a)}\\
\\
{\displaystyle -p^{-r}(pr-r+1)} & \mbox{ if }\Lambda=\Lambda_{(-pb,1)}
\end{cases}
\]
\end{prop}
\begin{proof}
We split into two cases, depending on $\Lambda$. 

\textbf{Case 1}: Suppose $\Lambda=\Lambda_{(1,-a')}$. Using Theorem \ref{Entries of Inv Int Matrix},
we compute
\begin{align*}
\sum_{\Lambda\ne\Lambda_{(1,0)}}c^{\Lambda_{(1,a')},\Lambda} & = \sum_{a=1}^{p^{r}-1}c^{\Lambda_{(1,-a')},\Lambda_{(1,-a)},}+\sum_{b=0}^{p^{r-1}-1}c^{\Lambda_{(1,-a')},\Lambda_{(-pb,1)}}\\
& = c^{\Lambda_{(1,a')},\Lambda_{(1,a')}}+\sum_{\substack{a=1\\
a\ne a'
}
}^{p^{r}-1}\left(-p^{1-2r}\frac{(pr-r+1)}{p+1}-\frac{p^{1-2r}(p-1)}{p+1}\nu_{p}\left(\frac{1}{a}-\frac{1}{a'}\right)\right)\\
 &   +\sum_{b=0}^{p^{r-1}-1}\left(-\frac{p^{1-2r}(pr-r+1)}{p+1}+\frac{p^{1-2r}(p-1)}{p+1}\nu_{p}(a')\right)\\
 & = c^{\Lambda_{(1,a')},\Lambda_{(1,a')}}+(p^{r}-2)\left(-p^{1-2r}\frac{(pr-r+1)}{p+1}\right)-\frac{p^{1-2r}(p-1)}{p+1}\sum_{\substack{a=1\\
a\ne a'
}
}^{p^{r}-1}\nu_{p}\left(\frac{1}{a}-\frac{1}{a'}\right)\\
 &   -(p^{r-1})\frac{p^{1-2r}(pr-r+1)}{p+1}+p^{r-1}\frac{p^{1-2r}(p-1)}{p+1}\nu_{p}(a')
\end{align*}
 Using Lemma \ref{P-adic Sum of Inv} on the sum $\sum v_{p}(\nicefrac{1}{a}-\nicefrac{1}{a'})$,
we have{\small{}
\begin{alignat*}{1}
 & =c^{\Lambda_{(1,a')},\Lambda_{(1,a')}}+(p^{r}-2)\left(-p^{1-2r}\frac{(pr-r+1)}{p+1}\right)-\frac{p^{1-2r}(p-1)}{p+1}\left(-(p^{r}-2)\nu_{p}(a')\right)\\
 & -(p^{r-1})\frac{p^{1-2r}(pr-r+1)}{p+1}+p^{r-1}\frac{p^{1-2r}(p-1)}{p+1}\nu_{p}(a')\\
 & =c^{\Lambda_{(1,a')},\Lambda_{(1,a')}}+(p^{r}-2+p^{r-1})\left(-p^{1-2r}\frac{(pr-r+1)}{p+1}\right)-\frac{p^{1-2r}(p-1)}{p+1}\left(-(p^{r}-2)\nu_{p}(a')\right)+p^{r-1}\frac{p^{1-2r}(p-1)}{p+1}\nu_{p}(a')\\
 & =c^{\Lambda_{(1,a')},\Lambda_{(1,a')}}+(p^{r}-2+p^{r-1})\left(-p^{1-2r}\frac{(pr-r+1)}{p+1}\right)+(p^{r}+p^{r-1}-2)\frac{p^{1-2r}(p-1)}{p+1}\nu_{p}(a')\\
 & ={\displaystyle -\frac{2p^{1-2r}(pr-r+1)}{p+1}+\frac{2p^{1-2r}(p-1)}{p+1}\nu_{p}(a')}+(p^{r}-2+p^{r-1})\left(-p^{1-2r}\frac{(pr-r+1)}{p+1}\right)\\
 & +(p^{r}+p^{r-1}-2)\frac{p^{1-2r}(p-1)}{p+1}\nu_{p}(a')\\
 & =-p^{-r}(pr-r+1)+p^{-r}(p-1)\nu_{p}(a').
\end{alignat*}
} 

\textbf{Case 2}: Suppose $\Lambda'=\Lambda_{(-pb',1)}$. Then
\begin{align*}
\sum_{\Lambda\ne\Lambda_{(1,0)}}c^{\Lambda_{(-pb',1)},\Lambda} & = \sum_{a=1}^{p^{r}-1}c^{\Lambda_{(-pb',1)},\Lambda_{(1,-a)},}+\sum_{b=0}^{p^{r-1}-1}c^{\Lambda_{(-pb',1)},\Lambda_{(-pb,1)}}\\
 & = \sum_{a=1}^{p^{r}-1}\left(-\frac{p^{1-2r}(pr-r+1)}{p+1}+\frac{p^{1-2r}(p-1)}{p+1}\nu_{p}(a)\right)\\
 &  -\frac{2p^{1-2r}(pr-r+1)}{p+1}+\sum_{\substack{b=0\\
b\ne b'
}
}^{p^{r-1}-1}\left(-\frac{p^{1-2r}(pr+p-r)}{p+1}-\frac{p^{1-2r}(p-1)}{p+1}\nu_{p}(b'-b)\right)\\
 & = -(p^{r}-1)\frac{p^{1-2r}(pr-r+1)}{p+1}+\frac{p^{1-2r}(p-1)}{p+1}\sum_{a=1}^{p^{r}-1}\nu_{p}(a)\\
 &  -\frac{2p^{1-2r}(pr-r+1)}{p+1}-(p^{r-1}-1)\frac{p^{1-2r}(pr+p-r)}{p+1}-\frac{p^{1-2r}(p-1)}{p+1}\sum_{\substack{b=0\\
b\ne b'
}
}^{p^{r-1}-1}\nu_{p}(b'-b)
\end{align*}
Using Lemma \ref{Sum of P-adic Vals}, we get
\begin{alignat*}{1}
 & =-(p^{r}-1)\frac{p^{1-2r}(pr-r+1)}{p+1}+\frac{p^{1-2r}(p-1)}{p+1}\left(\frac{p^{r}-pr+r-1}{p-1}\right)\\
 & -\frac{2p^{1-2r}(pr-r+1)}{p+1}-(p^{r-1}-1)\frac{p^{1-2r}(pr+p-r)}{p+1}-\frac{p^{1-2r}(p-1)}{p+1}\left(-\nu_{p}(b')+\frac{p^{r-1}-p(r-1)+r-2}{p-1}+\nu_{p}(b')\right)\\
 & =-\frac{(p-1)rp^{1-r}}{p+1}-\frac{2p^{1-2r}(pr-r+1)}{p+1}-(p^{r-1}-1)\frac{p^{1-2r}(pr+p-r)}{p+1}-\frac{p^{1-2r}(p-1)}{p+1}\left(\frac{p^{r-1}-p(r-1)+r-2}{p-1}\right)\\
 & =\frac{p^{-2r}(-p^{r+1}-rp^{r+2}+rp^{r}-p^{2}(r-1)+p(r-2))}{p+1}-\frac{p^{1-2r}(p-1)}{p+1}\left(\frac{p^{r-1}-p(r-1)+r-2}{p-1}\right)\\
 & =-p^{-r}(pr-r+1).\tag*{\qedhere}
\end{alignat*}
\end{proof}

\begin{thm}\label{Upper Bound}
Let $k\geq1$, $N\geq3$, and $r\geq1$ be integers and $p\geq2$ be a prime such that $p\nmid N$. Let $e$ be the exponent of $\pi$ in the annihilator of $M_{2k}(\Gamma(Np^{r}),\mathbb{Z}_{p}[\zeta_{Np^{r}}])/H^{0}(\mathfrak{X},\underline{\omega}^{\otimes2k})$. Then 
\[
e\leq2kp^{r-1}(pr-r+1).
\]
\end{thm}
\begin{proof}
For any $f\in M_{2k}(Np^{r},\mathbb{Z}_{p}[\zeta_{Np^{r}}])$, we always have $\nu_{\Lambda_{(1,0)}}(f)\geq0$. Thus
\begin{align*}
-\nu_{\Lambda}(f) & \leq \nu_{\Lambda_{(1,0)}}(f)-\nu_{\Lambda}(f)\\
 & = \sum_{\Lambda'\ne\Lambda_{(1,0)}}\left(H_{f}.\Lambda'-\deg(\underline{\omega}^{\otimes2k}|_{\Lambda'})\right)c^{\Lambda,\Lambda'}\\
 & = \sum_{\Lambda'\ne\Lambda_{(1,0)}}(H_{f}.\Lambda')c^{\Lambda,\Lambda'}-\sum_{\Lambda'\ne\Lambda_{(1,0)}}\deg(\underline{\omega}^{\otimes2k}|_{\Lambda'})c^{\Lambda,\Lambda'}.
\end{align*}

By Corollary \ref{Entries cij Are Negative}, $c^{\Lambda,\Lambda'}<0$. Also note that $H_{f}$ is an effective horizontal divisor since $f$ has no poles while $\Lambda'$ is an effective vertical divisor. Since $H_{f}$ and $\Lambda'$ do not have any common components, the intersection number $H_{f}.\Lambda'$ is positive. Thus $\sum_{\Lambda'\ne\Lambda_{(1,0)}}(H_{f}.\Lambda')c^{\Lambda,\Lambda'}$ will always be negative. Continuing, we have 
\[
\leq-\sum_{\Lambda'\ne\Lambda_{(1,0)}}\deg(\underline{\omega}^{\otimes2k}|_{\Lambda'})c^{\Lambda,\Lambda'}
\]
Using Theorem \ref{Deg of Res Modular Sheaf }, the above expression becomes 
\[
-\frac{k(p-1)p^{2r-1}\#{\rm SL}_{2}(\mathbb{Z}/N\mathbb{Z})}{12}\frac{24}{(p-1)\#{\rm SL}_{2}(\mathbb{Z}/N\mathbb{Z})}\deg{\rm S}(N)\sum_{\Lambda'\ne\Lambda_{(1,0)}}c^{\Lambda,\Lambda'}
\]
\[
=-kp^{2r-1}\deg{\rm S}(N)\sum_{\Lambda'\ne\Lambda_{(1,0)}}c^{\Lambda,\Lambda'}.
\]
Using Proposition \ref{Row Sums of cij}, we have
\[
=\begin{cases}
{\displaystyle 2kp^{r-1}(pr-r+1)-2kp^{r-1}(p-1)\nu_{p}(a)} & \mbox{ if }\Lambda=\Lambda_{(1,-a)}\\
\\
{\displaystyle 2kp^{r-1}(pr-r+1)} & \mbox{ if }\Lambda=\Lambda_{(-pb,1)}
\end{cases}
\]
which is maximized whenever $a$ is coprime to $p$ and attains a value of $2kp^{r-1}(pr-r+1)$.
\end{proof}

By replacing $\underline{\omega}^{\otimes2}$ with $\underline{\omega}^{\otimes2}(-\mathfrak{C}(Np^{r})$ in the proof of Theorem \ref{Upper Bound}, we obtain an upper bound for the exponent in the situation of cusp forms. 

\begin{cor}\label{Upper Bound for Cusp Forms}
Let $k\geq1,N\geq3$, and $r\geq1$ be integers and $p\geq2$ a prime such that $p\nmid N$. The exponent $e$ of $\pi$ in the annihilator of $S_{2k}(\Gamma(Np^{r}),\mathbb{Z}_{p}[\zeta_{Np^{r}}])/H^{0}(\mathfrak{X},\underline{\omega}^{\otimes2k}(-\mathfrak{C}(Np^{r}))$ is bounded above by 
\[
2kp^{r-1}(pr-r+1)-\frac{12k}{Np}(pr-r+1).
\]
\end{cor}
\begin{proof}
Using Theorem \ref{Deg of Res Modular Sheaf } to compute $\deg(\underline{\omega}^{\otimes2k}(-\mathfrak{C}(Np^{r}))|_{\Lambda})$, we get an upper bound of 
\[
-k\#{\rm SL}_{2}(\mathbb{Z}/N\mathbb{Z})(p-1)\left[\frac{p^{2r-1}}{12}-\frac{p^{r-1}}{2N}\right]\frac{24}{(p-1)\#{\rm SL}_{2}(\mathbb{Z}/N\mathbb{Z})}\deg{\rm S}(N)\sum_{\Lambda'\ne\Lambda_{(1,0)}}c^{\Lambda,\Lambda'}
\]

\begin{alignat*}{1}
 & =k\left[\frac{12p^{r-1}}{N}-2p^{2r-1}\right]\deg{\rm S}(N)\sum_{\Lambda'\ne\Lambda_{(1,0)}}c^{\Lambda,\Lambda'}\\
 & \leq-k\left[\frac{12p^{r-1}}{N}-2p^{2r-1}\right]p^{-r}(pr-r+1)\\
 & =2kp^{r-1}(pr-r+1)-\frac{12k}{Np}(pr-r+1).\tag*{\qedhere}
\end{alignat*}
\end{proof}

\begin{rmk}\label{Compare Edi Bound}
We compare our result in Corollary \ref{Upper Bound for Cusp Forms} with Edixhoven's method in \cite{[Edi06]}. He considers the situation of weight $2$ and level $\Gamma_{0}(N)$ cusp forms\emph{ }where ${\rm ord}_{p}(N)=1$, and bounds $e$ via the inequality 
\[
e<\frac{\deg(\Omega_{\mathfrak{X}_{0}(N)}|_{\Lambda})}{\deg{\rm S}(N)}.
\]
We will show a similar inequality holds our situation. Let $\Lambda_{0},\Lambda_{1},\dots,\Lambda_{n}$ denote the irreducible components of $\bar{\mathfrak{X}}$ where $\Lambda_{0}$ contains the cusp $\infty$. Let $f\in M_{2}(Np^{r},\mathbb{Z}_{p}[\zeta_{Np^{r}}])$ be a non-zero cusp form such that $\nu_{\Lambda_{0}}(f)=0$ and, without loss of generality, let $-m:=\nu_{\Lambda_{1}}(f)<0$ be the minimum among the values in $\left\{ \nu_{\Lambda_{i}}(f)\right\} _{i=1}^{n}$. By scaling, we can assume $\nu_{\Lambda_{0}}(f)=m$, $\nu_{\Lambda_{1}}(f)=0$, and $\nu_{\Lambda_{i}}(f)\geq0$ for $i=2,\dots,n$ so that $f$ has non-negative valuation along every irreducible component \emph{i.e.} $f\in H^{0}(\mathfrak{X},\underline{\omega}^{\otimes2})$. We can write the divisor associated to $f$ as 
\[
{\rm div}(f)=H_{f}+m\Lambda_{0}+\sum_{i=2}^{n}\nu_{\Lambda_{i}}(f)\Lambda_{i}
\]
where $H_{f}$ is the horizontal part of the divisor ${\rm div}(f)$. By Theorem \ref{iso of modular sheaf}, the sheaf of cusp forms $\underline{\omega}^{\otimes2}(-\mathfrak{C}(Np^{r}))$ is isomorphic to the relative dualizing sheaf $\Omega_{\mathfrak{X}/R}$. Since $\mathfrak{X}$ is an arithmetic surface, we can use intersection theory (see Theorem \ref{Props of Int Num}), along with Proposition
\ref{twist by div(f)}, to compute:{\small{}
\begin{align*}
\deg\Omega_{\mathfrak{X}/R}|_{\Lambda_{1}}(-m{\rm S}(N)) & = \deg(\Omega_{\mathfrak{X}/R}|_{\Lambda_{1}})-m\deg{\rm S}(N)\\
 & = \deg({\rm div}(f)|_{\Lambda_{1}})-m\deg{\rm S}(N)\\
 & = {\rm div}(f).\Lambda_{1}-m\deg{\rm S}(N)\\
 & = H_{f}.\Lambda_{1}+m\Lambda_{0}.\Lambda_{1}+\sum_{i=2}^{n}\nu_{\Lambda_{i}}(f)\Lambda_{i}.\Lambda_{1}-m\deg{\rm S}(N)\\
 & = H_{f}.\Lambda_{1}+m(\Lambda_{0}.\Lambda_{1}-\deg{\rm S}(N))+\sum_{i=2}^{n}\nu_{\Lambda_{i}}(f)\Lambda_{i}.\Lambda_{1}.
\end{align*}
}Since $f$ has no poles, $H_{f}$ is effective so $H_{f}.\Lambda_{1}\geq0$. By assumption, $\nu_{\Lambda_{i}}(f)\geq0$ for $i=2,\dots,n$; we also have $\Lambda_{i}.\Lambda_{1}\geq0$ since $\Lambda_{i}$ and $\Lambda_{1}$ have no common components. Lastly, $\Lambda_{0}$ and $\Lambda_{1}$ intersect precisely at the supersingular points, so $\Lambda_{0}.\Lambda_{1}\geq\deg{\rm S}(N)$. Therefore we conclude
\[
0\leq\deg\Omega_{\mathfrak{X}/R}|_{\Lambda_{1}}(-m{\rm S}(N))=\deg(\Omega_{\mathfrak{X}/R}|_{\Lambda_{1}})-m\deg{\rm S}(N).
\]

Using Corollary \ref{Deg RDS Restricted}, we have: {\small{}
\begin{align*}
m & < \frac{\deg(\Omega_{\mathfrak{X}/R}|_{\Lambda_{1}})}{\deg{\rm S}(N)}\\
 & = \frac{p^{r}\varphi(p^{r})\frac{\#{\rm SL}_{2}(\mathbb{Z}/N\mathbb{Z})}{24}-\varphi(p^{r})\#C(N)+\deg{\rm S}(N)\cdot p^{2r-1}}{\deg{\rm S}(N)}\\
 & = \frac{24p^{r}\varphi(p^{r})\#{\rm SL}_{2}(\mathbb{Z}/N\mathbb{Z})}{24(p-1)\#{\rm SL}_{2}(\mathbb{Z}/N\mathbb{Z})}-\frac{24\varphi(p^{r})\#{\rm SL}_{2}(\mathbb{Z}/N\mathbb{Z})}{2N(p-1)\#{\rm SL}_{2}(\mathbb{Z}/N\mathbb{Z})}+p^{2r-1}\\
 & = \frac{p^{r}p^{r-1}(p-1)}{p-1}-\frac{12p^{r-1}(p-1)}{N(p-1)}+p^{2r-1}\\
 & = 2p^{2r-1}-\frac{12p^{r-1}}{N}
\end{align*}
}which provides an upper bound for $e$. When $r=1$, this bound agrees with the bound in Corollary \ref{Upper Bound for Cusp Forms}, but in general is larger by a factor of $p^{r}/(pr-r+1)$. 
\end{rmk}

\section{\label{sec:A-Lower-Bound}A Lower Bound}
In this section, we will use Klein forms to build an explicit modular form in $M_{2}(\Gamma_{1}(p^{r}),\mathbb{Z}_{p}[\zeta_{p}])$. By viewing these modular forms at level $\Gamma(Np^{r})$ with coefficients in $\mathbb{Z}_{p}[\zeta_{Np^{r}}]$, we will obtain a lower bound for the exponent.

\begin{defn}\label{Klein form}
A \textbf{nearly holomorphic modular form }is a modular form which is allowed to be meromorphic at the cusps. 
\end{defn}

\noindent Fix $(r_{1},r_{2})\in\mathbb{Q}^{2}-\mathbb{Z}^{2}$. Let $\tau\in{\cal H}$, $q=e^{2\pi i\tau}$, and $q_{z}=e^{2\pi iz}$ where $z=r_{1}\tau+r_{2}$. Define the \textbf{Klein form} 
\begin{equation}
\kappa_{(r_{1},r_{2})}(\tau)=e^{\pi ir_{2}(r_{1}-1)}q^{\frac{1}{2}r_{1}(r_{1}-1)}(1-q_{z})\prod_{n=1}^{\infty}(1-q^{n}q_{z})(1-q^{n}q_{z}^{-1})(1-q^{n})^{-2}.\label{eq: Klein q-exp}
\end{equation}

The following result, which is \cite[Theorem. 2.6]{[EKS11]}, builds upon results in \cite[§2.1, §3.4]{[KL81]} which establish a criterion for when a product of Klein forms is a nearly holomorphic modular form. For $r\in\mathbb{R}$, we let $\left\langle r\right\rangle $ denote the fractional part of $r$. Note $\left\langle r\right\rangle =0$ precisely when $r\in\mathbb{Z}$. 

\begin{thm}\label{Prod of Klein}
For an integer $N\geq2$, let $\left\{ m(t)\right\} _{t=1}^{N-1}$ be a family of integers. Then the product 
\[
\kappa(\tau)=\prod_{t=1}^{N-1}\kappa_{(t/N,0)}(N\tau)^{m(t)}
\]
is a nearly holomorphic modular form for $\Gamma_{1}(N)$ of weight
$k=-\sum_{t=1}^{N-1}m(t)$ if 
\[
\sum_{t=1}^{N-1}m(t)t^{2}\equiv0\;({\rm mod\,}\gcd(2,N)\cdot N).
\]
Furthermore, for $\alpha=\left(\begin{array}{cc}
a & b\\
c & d
\end{array}\right)\in{\rm SL}_{2}(\mathbb{Z})$ we have 
\[
{\rm ord}_{q}\left(\kappa(\tau)|[\alpha]_{k}\right)=\frac{\gcd(c,N)^{2}}{2N}\sum_{t=1}^{N-1}m(t)\left\langle \frac{at}{\gcd(c,N)}\right\rangle \left(\left\langle \frac{at}{\gcd(c,N)}\right\rangle -1\right).
\]
\end{thm}

Using this result when $N=p^{r}>3$, we will choose a family of integers $\left\{ m(t)\right\} _{t=1}^{N-1}$ such that $k=2$ and ${\rm ord}_{q}\left(\kappa(\tau)|[\alpha]_{2}\right)\geq0$ for all $\alpha\in{\rm SL}_{2}(\mathbb{Z})$. This guarantees that $\kappa(\tau)$ is a weight 2 (holomorphic) modular form of level $\Gamma_{1}(p^{r})$. As $\Gamma(Np^{r})\leq\Gamma(p^{r})\leq\Gamma_{1}(p^{r})$, this also gives us a modular form of level $\Gamma(Np^{r})$ as originally desired. 

It may be natural to choose $m(t)$ to be zero for most values of $t$ in order to simplify the quadratic condition in Theorem \ref{Prod of Klein} and the expression for the order. Consider the situation when $m(t_{1}),m(t_{2}),$ and $m(t_{3})$ are the only non-zero values for some distinct $t_{1},t_{2},t_{3}\in\left\{ 1,2,\dots,p^{r}-1\right\} $. By Theorem \ref{Prod of Klein}, we seek $m(t_{1}),m(t_{2}),$ and $m(t_{3})$ such that 
\begin{equation}
m(t_{1})+m(t_{2})+m(t_{3})=-2\label{eq:KleinCond1}
\end{equation}
\begin{equation}
m(t_{1})t_{1}^{2}+m(t_{2})t_{2}^{2}+m(t_{3})t_{3}^{2}\equiv0\,({\rm mod\,}2p^{r})\label{eq:KleinCond2}
\end{equation}
and
\begin{equation}
\frac{\gcd(c,p^{r})^{2}}{2p^{r}}\left(\sum_{i=1}^{3}m(t_{i})\left\langle \frac{at_{i}}{\gcd(c,p^{r})}\right\rangle \left(\left\langle \frac{at_{i}}{\gcd(c,p^{r})}\right\rangle -1\right)\right)\geq0\label{eq:KleinCond3}
\end{equation}
for all $a,c\in\mathbb{Z}$. 

We will choose the $m(t_{i})$ satisfying these conditions in three separate cases depending on the level: $p>5$, $p=5$, and $p^{r}>3$ with $r\geq2$. It is not clear how one constructs similar modular forms of the remaining levels $p=2$ and $p=3$. We first consider the case of level $p>5$. To satisfy equation (\ref{eq:KleinCond1}), we can take $m(t_{1})=m(t_{2})=-2$ and $m(t_{3})=2$. The next equation (\ref{eq:KleinCond2}) becomes 
\[
-t_{1}^{2}-t_{2}^{2}+t_{3}^{2}\equiv0\,({\rm mod\,}p)\mbox{ or equivalently }t_{1}^{2}+t_{2}^{2}\equiv t_{3}^{2}\,({\rm mod\,}p)
\]
which can be satisfied if we take $(t_{1},t_{2},t_{3})$ to be a Pythagorean triple. As we now show, taking $t_{1}=3,t_{2}=4,$ and $t_{3}=5$ suffices if $p>5$. For convenience, we let $(x,y)$ denote $\gcd(x,y)$ for $x,y\in\mathbb{Z}$. 

\begin{prop}\label{Prod of Klein Forms is Hol}
The product of Klein forms 
\[
\kappa(\tau)=\kappa_{(3/p,0)}(p\tau)^{-2}\kappa_{(4/p,0)}(p\tau)^{-2}\kappa_{(5/p,0)}(p\tau)^{2}
\]
is a weight 2 modular form of level $\Gamma_{1}(p)$ for $p>5$. 
\end{prop}
\begin{proof}
Define $m:\{1,\dots,p-1\}\rightarrow\mathbb{Z}$ by 
\[
m(t)=\begin{cases}
-2 & \mbox{if }t=3,4\\
2 & \mbox{if }t=5\\
0 & \mbox{otherwise}
\end{cases}
\]
By Theorem \ref{Prod of Klein}, the product $\kappa(\tau)$ defined by $m(t)$ is a weight 2 nearly holomorphic modular form of level $\Gamma_{1}(p)$. It remains to show $\kappa(\tau)$ is holomorphic at the cusps. The order of $\kappa(\tau)$ is given by 
\[
{\rm ord}_{q}(\kappa(\tau)|[\alpha]_{2})
\]
{\small{}
\[
=\frac{(c,p)^{2}}{2p}\left(-2\left\langle \frac{3a}{(c,p)}\right\rangle \left(\left\langle \frac{3a}{(c,p)}\right\rangle -1\right)-2\left\langle \frac{4a}{(c,p)}\right\rangle \left(\left\langle \frac{4a}{(c,p)}\right\rangle -1\right)+2\left\langle \frac{5a}{(c,p)}\right\rangle \left(\left\langle \frac{5a}{(c,p)}\right\rangle -1\right)\right)
\]
}where $\alpha=\left(\begin{array}{cc}
a & b\\
c & d
\end{array}\right)\in{\rm SL}_{2}(\mathbb{Z})$. We will show this expression is always non-negative.

Since $(c,p)^{2}/2p>0$, we can ignore this factor. Dividing through by 2, our goal is to show
\[
-\left\langle \frac{3a}{(c,p)}\right\rangle \left(\left\langle \frac{3a}{(c,p)}\right\rangle -1\right)-\left\langle \frac{4a}{(c,p)}\right\rangle \left(\left\langle \frac{4a}{(c,p)}\right\rangle -1\right)+\left\langle \frac{5a}{(c,p)}\right\rangle \left(\left\langle \frac{5a}{(c,p)}\right\rangle -1\right)
\]
or equivalently 
\[
\left\langle \frac{3a}{(c,p)}\right\rangle \left(1-\left\langle \frac{3a}{(c,p)}\right\rangle \right)+\left\langle \frac{4a}{(c,p)}\right\rangle \left(1-\left\langle \frac{4a}{(c,p)}\right\rangle \right)-\left\langle \frac{5a}{(c,p)}\right\rangle \left(1-\left\langle \frac{5a}{(c,p)}\right\rangle \right)
\]
is non-negative for all $a,c\in\mathbb{Z}$.

Let 
\[
f(x)=\left\langle 3x\right\rangle \left(1-\left\langle 3x\right\rangle \right)+\left\langle 4x\right\rangle \left(1-\left\langle 4x\right\rangle \right)-\left\langle 5x\right\rangle \left(1-\left\langle 5x\right\rangle \right),
\]
noting that $f(a/(c,p))$ is the expression we are showing is non-negative. Since $\left\langle x\right\rangle $ is periodic with period $0\leq x<1$, it suffices to show $f(x)\geq0$ for all $0\leq x<1$. Since
\[
\left\langle -x\right\rangle (1-\left\langle -x\right\rangle )=(1-\left\langle x\right\rangle )(1-(1-\left\langle x\right\rangle ))=\left\langle x\right\rangle (1-\left\langle x\right\rangle )
\]
we can conclude $f(-x)=f(x)$. Therefore it suffices to show $f(x)\geq0$ for all $0\leq x\leq\nicefrac{1}{2}$. 

We accomplish this by finding an explicit piecewise defined expression for $f(x)$ and then computing its derivative. Note that 
\begin{align*}
\left\langle 3x\right\rangle \left(1-\left\langle 3x\right\rangle \right) & = \begin{cases}
3x(1-3x) & \mbox{if }0\leq x<\nicefrac{1}{3}\\
3(x-\nicefrac{1}{3})(1-3(x-\nicefrac{1}{3})) & \mbox{if }\nicefrac{1}{3}\leq x<\nicefrac{2}{3}\\
3(x-\nicefrac{2}{3})(1-3(x-\nicefrac{2}{3})) & \mbox{if }\nicefrac{2}{3}\leq x<1
\end{cases}\\
 & = \begin{cases}
3x(1-3x) & \mbox{if }0\leq x<\nicefrac{1}{3}\\
(3x-1)(2-3x) & \mbox{if }\nicefrac{1}{3}\leq x<\nicefrac{2}{3}\\
3(3x-2)(1-x) & \mbox{if }\nicefrac{2}{3}\leq x<1
\end{cases}
\end{align*}
Similar expressions can be obtained for $\left\langle 4x\right\rangle \left(1-\left\langle 4x\right\rangle \right)$ and $\left\langle 5x\right\rangle \left(1-\left\langle 5x\right\rangle \right)$. Putting these together, we have 
\[
f(x)=\begin{cases}
3x(1-3x)+4x(1-4x)-5x(1-5x) & \mbox{if }0\leq x<\nicefrac{1}{5}\\
3x(1-3x)+4x(1-4x)-(5x-1)(2-5x) & \mbox{if }\nicefrac{1}{5}\leq x<\nicefrac{1}{4}\\
3x(1-3x)+2(4x-1)(1-2x)-(5x-1)(2-5x) & \mbox{if }\nicefrac{1}{4}\leq x<\nicefrac{1}{3}\\
(3x-1)(2-3x)+2(4x-1)(1-2x)-(5x-1)(2-5x) & \mbox{if }\nicefrac{1}{3}\leq x<\nicefrac{2}{5}\\
(3x-1)(2-3x)+2(4x-1)(1-2x)-(5x-2)(3-5x) & \mbox{if }\nicefrac{2}{5}\leq x<\nicefrac{2}{4}
\end{cases}.
\]
Thus we can compute the derivative directly: 
\[
f'(x)=\begin{cases}
2x & \mbox{if }0\leq x<\nicefrac{1}{5}\\
2-8x & \mbox{if }\nicefrac{1}{5}\leq x<\nicefrac{1}{4}\\
0 & \mbox{if }\nicefrac{1}{4}\leq x<\nicefrac{1}{3}\\
6x-2 & \mbox{if }\nicefrac{1}{3}\leq x<\nicefrac{2}{5}\\
2-4x & \mbox{if }\nicefrac{2}{5}\leq x<\nicefrac{2}{4}
\end{cases}
\]
This shows that the pieces of $f(x)$ are either strictly increasing, strictly decreasing, or constant in their appropriate interval. For each interval, the following table shows if $f(x)$ is increasing, decreasing, or constant based on $f'(x)$. The value of $f$ at the rightmost endpoint is also calculated for each interval.

\begin{table}[h]
\centering{}%
\begin{tabular}{|c|c|c|c|c|c|}
\hline 
$a\leq x<b$ & $0\leq x<\nicefrac{1}{5}$ & $\nicefrac{1}{5}\leq x<\nicefrac{1}{4}$ & $\nicefrac{1}{4}\leq x<\nicefrac{1}{3}$ & $\nicefrac{1}{3}\leq x<\nicefrac{2}{5}$ & $\nicefrac{2}{5}\leq x<\nicefrac{1}{2}$\tabularnewline
\hline 
Inc, Dec, Con & increasing & decreasing & constant & increasing & decreasing\tabularnewline
\hline 
$f(b)$ & $\nicefrac{2}{5}$ & $0$ & $0$ & $\nicefrac{2}{5}$ & $0$\tabularnewline
\hline 
\end{tabular}
\end{table}
\noindent Thus we conclude $f(x)\geq0$ for all $x$. Consequently, $\kappa(\tau)$ is holomorphic at each cusp so is indeed a modular form.
\end{proof}

We compute the valuation of its $q$-expansion at each cusp. The $q$-expansion at $\infty$ for $\kappa_{(r_{1},r_{2})}(\tau)$ was given earlier in equation (\ref{Prod of Klein}). In particular, the $q$-expansion of $\kappa_{(\nicefrac{a}{p},0)}(p\tau)$ is 
\[
\kappa_{(\nicefrac{a}{p},0)}(p\tau)=q^{p\frac{1}{2}\frac{a}{p}(\frac{a}{p}-1)}(1-q^{a})\prod_{n=1}^{\infty}(1-q^{pn}q^{a})(1-q^{pn}q^{-a})(1-q^{pn})^{-2}.
\]
For convenience, let $H_a(q):=\prod_{n=1}^{\infty}(1-q^{pn}q^{a})(1-q^{pn}q^{-a})(1-q^{pn})^{-2}$ and note $1/H_a(q)\in 1+\mathbb{Z}[[q]]$.  The $q$-expansion of $\kappa(\tau)$ is

\begin{align*}
\kappa(\tau) & = \kappa_{(\nicefrac{3}{p},0)}(p\tau)^{-2}\kappa_{(\nicefrac{4}{p},0)}(p\tau)^{-2}\kappa_{(\nicefrac{5}{p},0)}(p\tau)^{2}\\
& = \frac{q^{p\frac{1}{2}\frac{5}{p}(\frac{5}{p}-1)2}(1-q^{5})^{2}}{q^{p\frac{1}{2}\frac{3}{p}(\frac{3}{p}-1)2}(1-q^{3})^{2}q^{p\frac{1}{2}\frac{4}{p}(\frac{4}{p}-1)2}(1-q^{4})^{2}}\left(\frac{H_{5}(q)}{H_{3}(q)H_{4}(q)}\right) \\
& = q^{2}\frac{(1-q^{5})^{2}}{(1-q^{3})^{2}(1-q^{4})^{2}}\left(\frac{H_{5}(q)}{H_{3}(q)H_{4}(q)}\right)
\end{align*}
\noindent which has all integral coefficients. 

To compute the $q$-expansion of $\kappa(\tau)$ at the other cusps, we use the following identity for Klein forms, which can be found in \cite[Prop. 2.1]{[EKS11]}. For any $\alpha=\left(\begin{array}{cc}
a & b\\
c & d
\end{array}\right)\in{\rm SL}_{2}(\mathbb{Z})$, we have 
\[
\kappa_{(r_{1},r_{2})}(\alpha\tau)=\kappa\left(\frac{a\tau+b}{c\tau+d}\right)=(c\tau+d)^{-1}\kappa_{\alpha(r_{1},r_{2})}(\tau).
\]
We let $\pi=1-\zeta_{p^{r}}$ which is a uniformizer of $\mathbb{Z}_{p}[\zeta_{Np^{r}}]$. 

\begin{prop}\label{Val of q-exp around 0}
The $\pi$-adic valuation of $\kappa(\tau)$ at its $q$-expansions around the cusp $0$ is $-2p$ i.e. $\nu_{0}(\kappa(\tau))=-2p$. 
\end{prop}
\begin{proof}
Let $\sigma=\left(\begin{array}{cc}
0 & 1\\
-1 & 0
\end{array}\right)\in{\rm SL}_{2}(\mathbb{Z})$; note $\sigma\cdot\infty=0$. The $q$-expansion of $\kappa(\tau)$ at the cusp $0$ is given by 
\[
\kappa(\tau)|[\sigma]_{2}=(-\tau)^{-2}\kappa(\sigma\cdot\tau)=\frac{\kappa_{(\nicefrac{5}{p},0)}(p\cdot\sigma\tau)^{2}}{\tau^{2}\kappa_{(3/p,0)}(p\cdot\sigma\tau)^{2}\kappa_{(4/p,0)}(p\cdot\sigma\tau)^{2}}.
\]
Observe that 
\[
\kappa_{(r_{1},r_{2})}(p\cdot\sigma\tau)=\kappa_{(r_{1},r_{2})}(-p/\tau)=\kappa_{(r_{1},r_{2})}(\sigma\cdot(\tau/p))
\]
\[
=(-\tau/p)^{-1}\kappa_{\sigma(r_{1},r_{2})}(\tau/p)=-(\tau/p)^{-1}\kappa_{(r_{2},-r_{1})}(\tau/p).
\]
Thus
\begin{align*}
\kappa(\tau)|[\sigma]_{2} & = \frac{(\tau/p)^{-2}\kappa_{(0,-\nicefrac{5}{p})}(\tau/p)^{2}}{\tau^{2}(\tau/p)^{-2}\kappa_{(0,-\nicefrac{3}{p})}(\tau/p)^{2}(\tau/p)^{-2}\kappa_{(0,-\nicefrac{4}{p})}(\tau/p)^{2}}\\
 & = \frac{\kappa_{(0,-\nicefrac{5}{p})}(\tau/p)^{2}}{p^{2}\kappa_{(0,-\nicefrac{3}{p})}(\tau/p)^{2}\kappa_{(0,-\nicefrac{4}{p})}(\tau/p)^{2}}.
\end{align*}
Using equation (\ref{eq: Klein q-exp}), the $q$-expansion of $\kappa_{(0,-\nicefrac{a}{p})}(\tau/p)$ is 
\[
\kappa_{(0,-\nicefrac{a}{p})}(\tau/p)=e^{\pi i\frac{a}{p}}(1-\zeta_{p}^{-a})\prod_{n=1}^{\infty}(1-q^{n/p}\zeta_{p}^{-a})(1-q^{n/p}\zeta_{p}^{a})(1-q^{n/p})^{-2}
\]
where $\zeta_{p}=e^{2\pi i/p}$. Therefore 
\begin{align*}
\kappa(\tau)|[\sigma]_{2} & = \frac{e^{\pi i\frac{5}{p}}(1-\zeta_{p}^{-5})^{2}}{p^{2}e^{\pi i\frac{3}{p}}(1-\zeta_{p}^{-3})^{2}e^{\pi i\frac{4}{p}}(1-\zeta_{p}^{-4})}H(q^{1/n})\\
& = \frac{\zeta_{p}^{-1}(1-\zeta_{p}^{-5})^{2}}{p^{2}(1-\zeta_{p}^{-3})^{2}(1-\zeta_{p}^{-4})^{2}}H(q^{1/n})
\end{align*}
for some $H(q^{1/n})\in 1+\mathbb{Z}[[q^{1/n}]]$. Note that $1-\zeta_{p}^{a}$ is a uniformizer for $\mathbb{Z}_{p}[\zeta_{Np}]$
for $a$ coprime to $p$. Furthermore, 
\[
(1-\zeta_{p}^{a})^{(p-1)}\mathbb{Z}_{p}[\zeta_{Np}]=p\mathbb{Z}_{p}[\zeta_{Np}].
\]
Thus the minimal valuation among the coefficients of the $q$-expansion of $\kappa(\tau)|[\sigma]_{2}$ is 
\[
\nu_{\pi}\left(\frac{\zeta_{p}^{-1}(1-\zeta_{p}^{-5})^{2}}{p^{2}(1-\zeta_{p}^{-3})^{2}(1-\zeta_{p}^{-4})^{2}}\right)=\nu_{\pi}\left(\frac{1}{p^{2}(1-\zeta_{p}^{-4})^{2}}\right)=-2\left(\nu_{\pi}(p)+\nu_{\pi}(1-\zeta_{p}^{-4})\right)
\]
\[
=-2\left((p-1)+1\right)=-2p.\qedhere
\]
\end{proof}

Next we will handle the case of $p=5$. 
\begin{prop}\label{Level 5 Prod of Klein}
The product of Klein forms
\[
\kappa(\tau)=\kappa_{(\nicefrac{1}{5},0)}(5\tau)^{4}\kappa_{(\nicefrac{2}{5},0)}(5\tau)^{-2}\kappa_{(\nicefrac{3}{5},0)}(5\tau)^{-4}
\]
is a weight 2 modular form of level $\Gamma_{1}(5)$. Furthermore,
$\nu_{0}(\kappa(\tau))=-2\cdot5$.
\end{prop}
\begin{proof}
Define $m:\{1,2,3,4\}\rightarrow\mathbb{Z}$ by 
\[
m(t)=\begin{cases}
4 & \mbox{if }t=1\\
-2 & \mbox{if }t=2\\
-4 & \mbox{if }t=3\\
0 & \mbox{if }t=4
\end{cases}
\]
Note that
\[
m(1)+m(2)+m(3)=-2
\]
and 
\[
m(1)\cdot1^{2}+m(2)\cdot2^{2}+m(3)\cdot3^{2}=4\cdot1^{2}-2\cdot2^{2}-4\cdot3^{2}
\]
\[
=-40\equiv0\;({\rm mod}\,10).
\]
By Theorem \ref{Prod of Klein}, the product $\kappa(\tau)$ defined by $m(t)$ is a weight 2 nearly holomorphic modular form of level $\Gamma_{1}(5)$. It remains to show $\kappa(\tau)$ is holomorphic at the cusps. Similar to the proof in Proposition \ref{Prod of Klein Forms is Hol}, it suffices to show 
\[
4\left\langle \frac{a}{(c,5)}\right\rangle \left(\left\langle \frac{a}{(c,5)}\right\rangle -1\right)-2\left\langle \frac{2a}{(c,5)}\right\rangle \left(\left\langle \frac{2a}{(c,5)}\right\rangle -1\right)-4\left\langle \frac{3a}{(c,5)}\right\rangle \left(\left\langle \frac{3a}{(c,5)}\right\rangle -1\right)\geq0
\]
for all $a,c\in\mathbb{Z}$. Let 
\[
f(x)=4\left\langle x\right\rangle \left(\left\langle x\right\rangle -1\right)-2\left\langle 2x\right\rangle \left(\left\langle 2x\right\rangle -1\right)-4\left\langle 3x\right\rangle \left(\left\langle 3x\right\rangle -1\right).
\]
As before, it suffices to show $f(x)\geq0$ for $0\leq x\leq\nicefrac{1}{2}$. Since we consider $x$ of the form $a/(c,5)$, it suffices to show $f(x)\geq0$ when $x=0,\nicefrac{1}{5},\nicefrac{2}{5}$ or $\nicefrac{1}{2}$. We compute directly:
\[
f(0)=f(\nicefrac{1}{2})=0
\]
\[
f(\nicefrac{1}{5})=-4\left(\frac{1}{5}\right)\left(\frac{4}{5}\right)+6\left(\frac{2}{5}\right)\left(\frac{3}{5}\right)=\frac{4}{5}
\]
\[
f(\nicefrac{2}{5})=-4\left(\frac{2}{5}\right)\left(\frac{3}{5}\right)+6\left(\frac{4}{5}\right)\left(\frac{1}{5}\right)=0.
\]
Consequently, $\kappa(\tau)$ is holomorphic at each cusp so is indeed a modular form. 

The same proof holds from Proposition \ref{Val of q-exp around 0} to show the $q$-expansion around the cusp 0 of our modular form $\kappa(\tau)$ from Proposition \ref{Level 5 Prod of Klein} has $\pi$-adic valuation equal to $-2\cdot5$. 
\end{proof}

Lastly we handle the case $p^{r}>3$ with $r\geq2$. 

\begin{prop}\label{v_0(kappa) r2}
Suppose $p^{r}>3$ and $r\geq2$. The product of Klein forms 
\[
\kappa(\tau)=\kappa_{(\nicefrac{p^{r}-1}{p^{r}},0)}(p^{r}\tau)^{-2}\kappa_{(\nicefrac{p^{r-1}}{p^{r}},0)}(p^{r}\tau)^{-2}\kappa_{(\nicefrac{1}{p^{r}},0)}(p^{r}\tau)^{2}
\]
is a weight 2 modular form of level $\Gamma_{1}(p^{r})$. Furthermore, $\nu_{0}(\kappa(\tau))=-2p^{r-1}(pr-r+1)$. 
\end{prop}
\begin{proof}
Define $m(t):\{1,\dots,p^{r}-1\}\rightarrow\mathbb{Z}$ by 
\[
m(t)=\begin{cases}
-2 & \mbox{if }t=p^{r-1},p^{r}-1\\
2 & \mbox{if }t=1\\
0 & {\rm otherwise}
\end{cases}
\]
Since $r\geq2$, the values $p^{r-1},p^{r}-1,$ and $1$ are all distinct. We have 
\begin{align*}
\sum_{t=1}^{p^{r}-1}m(t)t^{2} & = -2\cdot(p^{r-1})^{2}-2\cdot(p^{r}-1)^{2}+2\cdot1^{2}\\
 & = -2p^{2r-2}-2(p^{r}-1)^{2}+2\\
 & \equiv -2(p^{r}-1)^{2}+2\mod p^{r}\\
 & \equiv -2+2\mod p^{r}\\
 & \equiv 0\mod2p^{r}
\end{align*}
noting that $r\geq2$ so $2r-2\geq r$. Furthermore, $-\sum_{t=1}^{p^{r-1}}m(t)=2$ so $\kappa(\tau)$ is a weight 2 level $\Gamma_{1}(p^{r})$ nearly holomorphic modular form. It remains to show $\kappa(\tau)$ has non-negative order at each cusp. 

The order of $\kappa(\tau)$ is given by
\[
{\rm ord}_{q}(\kappa(\tau)|[\alpha]_{2})
\]
{\small{}
\[
=\frac{(c,p^{r})^{2}}{2p^{r}}\left(-2\left\langle \frac{(p^{r}-1)a}{(c,p^{r})}\right\rangle \left(\left\langle \frac{(p^{r}-1)a}{(c,p^{r})}\right\rangle -1\right)-2\left\langle \frac{p^{r-1}a}{(c,p^{r})}\right\rangle \left(\left\langle \frac{p^{r-1}a}{(c,p)}\right\rangle -1\right)+2\left\langle \frac{a}{(c,p^{r})}\right\rangle \left(\left\langle \frac{a}{(c,p^{r})}\right\rangle -1\right)\right)
\]
}where $\alpha=\left(\begin{array}{cc}
a & b\\
c & d
\end{array}\right)\in{\rm SL}_{2}(\mathbb{Z})$. We will show this expression is always non-negative. Similar to the proof of Proposition \ref{Prod of Klein Forms is Hol}, we let
\[
f(x)=\left\langle (p^{r}-1)x\right\rangle \left(1-\left\langle (p^{r}-1)x\right\rangle \right)+\left\langle p^{r-1}x\right\rangle \left(1-\left\langle p^{r-1}x\right\rangle \right)-\left\langle x\right\rangle \left(1-\left\langle x\right\rangle \right).
\]
Writing $(c,p^{r})=p^{d}$ for some $0\leq d\leq r$, it suffices to show $f(a/p^{d})\geq0$ for all $0\leq a/p^{d}\leq\nicefrac{1}{2}$. Note that
\[
\frac{(p^{r}-1)a}{p^{d}}=\frac{p^{r}}{p^{d}}-\frac{a}{p^{d}}.
\]
Since $p^{r}/p^{d}\in\mathbb{Z}$ and $a/p^{d}\leq1$, we have 
\[
\left\langle \frac{(p^{r}-1)a}{(c,p^{r})}\right\rangle =1-\frac{a}{p^{d}}.
\]
Therefore 
\begin{align*}
f(a/p^{d}) & = \left(1-\frac{a}{p^{d}}\right)\left(\frac{a}{p^{d}}\right)+\left\langle \frac{p^{r-1}a}{p^{d}}\right\rangle \left(1-\left\langle \frac{p^{r-1}a}{p^{d}}\right\rangle \right)-\left(\frac{a}{p^{d}}\right)\left(1-\frac{a}{p^{d}}\right)\\
 & = \left\langle \frac{p^{r-1}a}{p^{d}}\right\rangle \left(1-\left\langle \frac{p^{r-1}a}{p^{d}}\right\rangle \right)
\end{align*}
which is always non-negative. Thus $\kappa(\tau)$ is holomorphic at each cusp. 

We compute $\nu_{0}(\kappa(\tau))$ in the same way as in Proposition \ref{Val of q-exp around 0}. Let $\sigma=\left(\begin{array}{cc}
0 & 1\\
-1 & 0
\end{array}\right)$; the $q$-expansion of $\kappa(\tau)$ at the cusp $0$ is given by 
\begin{align*}
\kappa(\tau)|_{2}[\sigma] & = (-\tau)^{-2}\kappa(\sigma\cdot\tau)\\
 & = \frac{\kappa_{(\nicefrac{1}{p^{r}},0)}(p^{r}\cdot\sigma\tau)^{2}}{\kappa_{(\nicefrac{p^{r}-1}{p^{r}},0)}(p^{r}\cdot\sigma\tau)^{2}\kappa_{(\nicefrac{p^{r-1}}{p^{r}},0)}(p^{r}\cdot\sigma\tau)^{2}}.
\end{align*}
Observe that 
\begin{align*}
\kappa_{(a/p^{r},0)}(p^{r}\cdot\sigma\tau) & = \kappa_{(a/p^{r},0)}(-p^{r}/\tau)\\
 & = \kappa_{(a/p^{r},0)}(\sigma\cdot(\tau/p^{r}))\\
 & = (\tau/p^{r})^{-1}\kappa_{\sigma(a/p^{r},0)}(\tau/p^{r})\\
 & = (\tau/p^{r})^{-1}\kappa_{(0,-a/p^{r})}(\tau/p^{r}).
\end{align*}

Similar to our computation in the proof of Proposition \ref{Val of q-exp around 0}, we have 
\begin{align*}
\kappa(\tau)|_{2}[\sigma] & = \frac{\kappa_{(0,-\nicefrac{1}{p^{r}})}(\tau/p^{r})^{2}}{\tau^{2}(\tau/p^{r})^{-2}\kappa_{(0,-\nicefrac{p^{r}-1}{p^{r}})}(\tau/p^{r})^{2}\kappa_{(0,-\nicefrac{p^{r-1}}{p^{r}})}(\tau/p^{r})^{2}}\\
 & = \frac{\kappa_{(0,-\nicefrac{1}{p^{r}})}(\tau/p^{r})^{2}}{p^{2r}\kappa_{(0,-\nicefrac{p^{r}-1}{p^{r}})}(\tau/p^{r})^{2}\kappa_{(0,-\nicefrac{p^{r-1}}{p^{r}})}(\tau/p^{r})^{2}}\\
 & = \frac{e^{-2\pi i/p^{r}}(1-\zeta_{p^{r}}^{-1})^{2}}{p^{2r}e^{-2\pi i(p^{r}-1)/p^{r}}e^{-2\pi i(p^{r-1})/p^{r}}(1-\zeta_{p^{r}}^{-(p^{r}-1)})^{2}(1-\zeta_{p^{r}}^{-p^{r-1}})^{2}}F(q^{1/p^{r}})\\
 & = \frac{\zeta_{p^{r}}^{-1}(1-\zeta_{p^{r}}^{-1})^{2}}{p^{2r}\zeta_{p^{r}}^{1-p^{r}}\zeta_{p^{r}}^{p^{r-1}}(1-\zeta_{p^{r}}^{-(p^{r}-1)})^{2}(1-\zeta_{p^{r}}^{-p^{r-1}})^{2}}F(q^{1/p^{r}}).
\end{align*}
for some explicit $F(q^{1/p^{r}})\in1+\mathbb{Z}[\zeta_{p^{r}}][[q^{1/p^{r}}]]$. Let $\nu_{\pi_{r}}$ denote the normalized valuation in $\mathbb{Z}_{p}[\zeta_{p^{r}}]$ so that $\nu_{\pi_{r}}(1-\zeta_{p^{r}})=1$. Then we have 
\begin{align*}
\nu_{0}(\kappa(\tau)) & = \nu_{\pi_{r}}\left(\frac{\zeta_{p^{r}}^{-1}(1-\zeta_{p^{r}}^{-1})^{2}}{p^{2r}\zeta_{p^{r}}^{1-p^{r}}\zeta_{p^{r}}^{p^{r-1}}(1-\zeta_{p^{r}}^{-(p^{r}-1)})^{2}(1-\zeta_{p^{r}}^{-p^{r-1}})^{2}}\right)\\
 & = \nu_{\pi_{r}}\left(\frac{1}{p^{2r}(1-\zeta_{p^{r}}^{-p^{r-1}})^{2}}\right)\\
 & = -2\nu_{\pi_{r}}(p^{r}(1-\zeta_{p^{r}}^{-p^{r-1}})).
\end{align*}
Note that $(1-\zeta_{p^{r}}^{-p^{r-1}})\mathbb{Z}_{p}[\zeta_{p^{r}}]=(1-\zeta_{p}^{-1})\mathbb{Z}_{p}[\zeta_{p^{r}}]$
and $(1-\zeta_{p}^{-1})^{\varphi(p)}\mathbb{Z}_{p}[\zeta_{p^{r}}]=p\mathbb{Z}_{p}[\zeta_{p^{r}}].$
Hence 
\[
\nu_{\pi_{r}}(1-\zeta_{p^{r}}^{-p^{r-1}})=\varphi(p^{r})/\varphi(p)=p^{r-1}.
\]
Continuing, 
\begin{align*}
\nu_{0}(\kappa(\tau)) & = -2(rp^{r-1}(p-1)+p^{r-1})\\
 & = -2p^{r-1}(pr-r+1)
\end{align*}
as desired.
\end{proof}

\begin{thm}
Let $k\geq1,N\geq3$ and $r\geq1$ be integers and $p\geq2$ be a prime such that $p\nmid N$ and $p^{r}>3$. The exponent $e$ of $\pi$ in the annihilator of $M_{2k}(\Gamma(Np^{r}),\mathbb{Z}_{p}[\zeta_{Np^{r}}])/H^{0}(\mathfrak{X}(Np^{r}),\underline{\omega}_{\mathfrak{X}(Np^{r})}^{\otimes2k})$ is equal to $2kp^{r-1}(pr-r+1)$.
\end{thm}
\begin{proof}
In Theorem \ref{Upper Bound} we showed $2kp^{r-1}(pr-r+1)$ was an upper bound for $e$. Let $\kappa(\tau)\in M_{2}(\Gamma(Np^{r}),\mathbb{Z}_{p}[\zeta_{Np^{r}}])$ be the modular form as in Propositions \ref{Prod of Klein Forms is Hol}, \ref{Level 5 Prod of Klein}, and \ref{v_0(kappa) r2}. We have shown that $\nu_{0}(\kappa)=-2p^{r-1}(pr-r+1)$. In general, we can consider the $k$th power $\kappa(\tau)^{k}\in M_{2k}(\Gamma(Np^{r}),\mathbb{Z}_{p}[\zeta_{Np^{r}}])$ so that $\nu_{0}(\kappa^{k})=-2kp^{r-1}(pr-r+1)$. Then  $\pi^{2kp^{r-1}(pr-r+1)}\kappa(\tau)^{k}$ has integral $q$-expansion at the cusp $0$ so $2kp^{r-1}(pr-r+1)$ is a lower bound for $e$. Hence $e=2kp^{r-1}(pr-r+1)$.
\end{proof}

\appendix

\section{{\label{Ap - mod curve}Appendix \textendash{} The modular curve and modular forms}}

We will summarize the formulation of the regular integral model of the modular curve of full level, as presented in \cite{[KM85]}. Let $R$ be a ring and let ${\rm Ell}_{R}$ denote the category whose objects are elliptic curves $E\rightarrow S$ where $S$ is an $R$-scheme, and whose morphisms $(E\rightarrow S)\rightarrow(E'\rightarrow S')$ are Cartesian squares.

We will concern ourselves with the representability of moduli problems ${\cal F}:{\rm Ell}_{R}\rightarrow{\rm Sets}$.  The functor ${\cal F}$ induces a functor $\tilde{{\cal F}}:{\rm Sch}_{R}\rightarrow{\rm Sets}$ defined by 
\[
S\mapsto\left\{ \left[(E\rightarrow S,\gamma)\right]:\gamma\in{\cal F}(E\rightarrow S)\right\}
\]
where $\left[(E\rightarrow S,\gamma)\right]$ denotes the isomorphism class of the pair $(E\rightarrow S,\gamma)$. If ${\cal F}$ is representable by ${\cal E}\rightarrow\mathfrak{M}({\cal F})$, then $\tilde{{\cal F}}$ is representable by $\mathfrak{M}({\cal F})$. Indeed, from the bijection
\[
\Psi:{\cal F}(E\rightarrow S)\rightarrow{\rm Hom}_{{\rm Ell}_{R}}(E\rightarrow S,{\cal E}\rightarrow\mathfrak{M}({\cal F}))
\]
we can associate to $\gamma\in{\cal F}(E\rightarrow S)$ a morphism
\[
\xymatrix{E\ar@{->}[r]^{\alpha}\ar@{->}[d]_{\pi} & {\cal E}\ar@{->}[d]^{\pi'}\\
S\ar@{->}[r]_{f} & \mathfrak{M}({\cal F})
}
\]
in ${\rm Ell}_{R}$. We define the map $\tilde{\Psi}:\tilde{{\cal F}}(S)\rightarrow{\rm Hom}_{{\rm Sch}_{R}}(S,\mathfrak{M}({\cal F}))$ by sending $[(E\overset{\pi}{\rightarrow}S,\gamma)]$ to the map $\mbox{\ensuremath{f:S\rightarrow\mathfrak{M}({\cal F})}}$. This correspondence is bijective and functorial in $S$ by the properties of $\Psi$. 

\begin{defn}
Let ${\cal F},{\cal F}':{\rm Ell}_{R}\rightarrow{\rm Sets}$ be two moduli problems for elliptic curves. A \textbf{morphism between moduli problems }(\textbf{over ${\rm Ell}_{R}$})\textbf{ }is a natural transformation $\eta:{\cal F}\Rightarrow{\cal F}'$.
\end{defn}

Let $\eta:{\cal F}\Rightarrow{\cal F}'$ be a morphism between two representable moduli problems over ${\rm Ell}_{R}$. We have an induced natural transformation $\tilde{\eta}:\tilde{{\cal F}}\Rightarrow\tilde{{\cal F}}$. The map $\eta$ induces a map $\eta:\mathfrak{M}({\cal F})\rightarrow\mathfrak{M}({\cal F}')$ which on $T$-points coincides with the map $\tilde{\eta}_{T}$ in the following diagram:

\[
\xymatrix{\mathfrak{M}({\cal F})(T)\ar@{->}[r]^{\eta}\ar@{=}[d] & \mathfrak{M}({\cal F}')(T)\ar@{=}[d]\\
{\rm Hom}_{{\rm Sch}_{R}}(T,\mathfrak{M}({\cal F}))\ar@{-}[d]^{\simeq} & {\rm Hom}_{{\rm Sch}_{R}}(T,\mathfrak{M}({\cal F}'))\ar@{-}[d]^{\simeq}\\
\tilde{{\cal F}}(T)\ar@{->}[r]^{\eta_{T}} & \tilde{{\cal F}}'(T)
}
\]

\begin{defn}
Let ${\cal F}$ and ${\cal F}'$ be two moduli problems for elliptic curves over $R$. We define the \textbf{product},\textbf{ }or \textbf{simultaneous},\textbf{ moduli problem} ${\cal F}\times{\cal F}'$ by 
\[(E\rightarrow S)\mapsto{\cal F}(E\rightarrow S)\times{\cal F}'(E\rightarrow S).\]
\end{defn}

Suppose ${\cal F}$ is representable by an elliptic curve ${\cal E}\rightarrow\mathfrak{M}({\cal F})$ and ${\cal F}'$ is relatively representable. According to \cite[{}4.3.4]{[KM85]}, ${\cal F}\times{\cal F}'$ is representable and $\mathfrak{M}({\cal F}\times{\cal F}')={\cal F}'_{{\cal E}/\mathfrak{M}({\cal F})}$ so we naturally have a map $\mathfrak{M}({\cal F}\times{\cal F}')\rightarrow\mathfrak{M}({\cal F})$. Following the notation in \cite{[KM85]}, we will usually denote $\mathfrak{M}({\cal F}\times{\cal F}')$ by $\mathfrak{M}({\cal F},{\cal F}')$.

\begin{defn}
Let $E$ be an elliptic curve over a scheme $S$. A section $s\in E(S)$ corresponds to a morphism\footnote{In general, a section of a morphism $X\rightarrow Y$ will be a locally closed embedding $s:Y\rightarrow X$. When $X\rightarrow Y$ is separated, as in the case of an elliptic curve, the map $s:Y\rightarrow X$ becomes a closed immersion. The associated divisor $\left[s\right]$ is in general defined to be the closure of the scheme-theoretic image of $s$ in $E$.} $s:S\rightarrow E$ whose composition with the structural morphism $E\rightarrow S$ is the identity on $S$. Since $E\rightarrow S$ is separated, $s$ is a closed immersion. We denote by $\left[s\right]$ the scheme-theoretic image of $s$ in $E$ which is an effective Cartier divisor on $E$. Consider the multiplication-by-$N$ isogeny $\left[N\right]:E\rightarrow E$ whose kernel we denote by $E[N]$, an $S$-group scheme. 

A $\Gamma(N)$\textbf{-structure }on $E\rightarrow S$ is a group homomorphism $\phi:(\mathbb{Z}/N\mathbb{Z})^{2}\rightarrow E[N](S)$ such that we have an equality of effective Cartier divisors
\[
\sum_{v\in(\mathbb{Z}/N\mathbb{Z})^{2}}[\phi(v)]=E[N].
\]
In this case, we call $P=\phi(1,0)$ and $Q=\phi(0,1)$ the corresponding \textbf{Drinfeld basis }of $E[N]$.
\end{defn}

We define $\left[\Gamma(N)\right]:{\rm Ell}_{\mathbb{Z}}\rightarrow{\rm Sets}$ to be the moduli problem which assigns to each elliptic curve $E\rightarrow S$ the set of all $\Gamma(N)$-structures on $E\rightarrow S$. Let $(\alpha,f):(E\rightarrow S)\rightarrow(E'\rightarrow S')$ be a morphism in the category ${\rm Ell}_{R}$. Then 
\[
\left[\Gamma(N)\right](\alpha,f):\left[\Gamma(N)\right](E'/S')\rightarrow\left[\Gamma(N)\right](E/S)
\]
is defined by sending a Drinfeld basis $(P,Q)$ of $E'[N]$ to $(\alpha^{*}P,\alpha^{*}Q)$. Indeed, $(\alpha^{*}P,\alpha^{*}Q)$ is a Drinfeld basis of $E[N]$ (see the paragraph proceeding \cite[{}1.4.1.2]{[KM85]}).

We moreover consider the full level $N$ moduli problem over the cyclotomic integers $\mathbb{Z}[\zeta_{N}]:=\mathbb{Z}[X]/(\Phi_{N}(X))$, where $\Phi_{N}(X)\in\mathbb{Z}[X]$ is the $N$th cyclotomic polynomial, by constructing a new moduli problem $\left[\Gamma(N)\right]^{{\rm can}}$ from $\left[\Gamma(N)\right]$ as follows (see also \cite[{}9.4.3.1]{[KM85]}). For a $\mathbb{Z}[\zeta_{N}]$-algebra $A$, let $\zeta_{N}$ denote the image of $X$ mod $\Phi_{N}(X)$ under the map $\mathbb{Z}[\zeta_{N}]\rightarrow A$, following the convention in \cite[{}9.1.5]{[KM85]}. Section 2.8 of \cite{[KM85]} associates to the isogeny $\left[N\right]$ a bilinear pairing $\mbox{\ensuremath{e_{N}:E[N]\times E[N]\rightarrow\mu_{N}}}.$ Define the moduli problem $\left[\Gamma(N)\right]^{{\rm can}}:{\rm Ell}_{\mathbb{Z}[\zeta_{N}]}\rightarrow{\rm Sets}$ which assigns to each elliptic curve $E\rightarrow S$ the set of all Drinfeld bases $(P,Q)$ such that $e_{N}(P,Q)=\zeta_{N}$.

\begin{thm}
\label{KM5.1.1}Suppose $N\geq3$. The moduli problem $\left[\Gamma(N)\right]^{{\rm can}}$ is represented by a regular scheme $\mathfrak{Y}(N):=\mathfrak{M}([\Gamma(N)]^{{\rm can}})$ which is flat over ${\rm Spec}(\mathbb{Z}[\zeta_{N}])$ of dimension 2. Moreover, $\mathfrak{Y}(N)$ is smooth over $\mathbb{Z}[\zeta_{N},1/N]$.
\end{thm}

\begin{proof}
The moduli problem $\left[\Gamma(N)\right]$ is relatively representable and finite, flat over ${\rm Ell}_{\mathbb{Z}}$ by \cite[{}5.1.1]{[KM85]}. Furthermore, $\left[\Gamma(N)\right]$ is rigid whenever $N\geq3$ by \cite[{}2.7.2]{[KM85]}. Therefore, by Proposition \cite[{}4.7.0]{[KM85]} and \cite[{}5.1.1]{[KM85]}, $\left[\Gamma(N)\right]$ is representable by a regular scheme flat over $\mathbb{Z}$ of dimension 2. By \cite[{}9.1.8, {}9.1.9]{[KM85]}, the same things hold true for the associated canonical moduli problem $\left[\Gamma(N)\right]^{{\rm can}}$. Using \cite[§5.1.1]{[KM85]} again, $\left[\Gamma(N)\right]$ is \'{e}tale over ${\rm Ell}_{\mathbb{Z}[1/N]}$, which, along with \cite[{}4.7.1]{[KM85]}, implies smoothness.
\end{proof}

\begin{rmk}
Theorem \ref{KM5.1.1} also holds for the moduli problem $\left[\Gamma(N)\right]^{{\rm can}}$ over ${\rm Ell}_{R}$, for any ring extension $\mathbb{Z}[\zeta_{N}]\rightarrow R$. For any moduli problem ${\cal F}:{\rm Ell}_{\mathbb{Z}}\rightarrow{\rm Sets}$ and ring $R$, let ${\cal F}_{R}:{\rm Ell}_{R}\rightarrow{\rm Sets}$ denote the moduli problem obtained by composing ${\cal F}$ with the forgetful functor ${\rm Ell}_{R}\rightarrow{\rm Ell}_{\mathbb{Z}}$. According to \cite[{}4.13]{[KM85]}, if ${\cal F}$ is relatively representable, then ${\cal F}_{R}$ is also relatively representable by the same morphism ${\cal F}_{E/S}\rightarrow S$ for any $R$-scheme $S$. Furthermore, if ${\cal F}$ is representable by ${\cal E}\rightarrow\mathfrak{M}({\cal F})$, then ${\cal F}_{R}$ is representable by the base change ${\cal E}_{R}\rightarrow\mathfrak{M}({\cal F})_{R}$. In other words, we have $\mathfrak{M}({\cal F}_{R})=\mathfrak{M}({\cal F})_{R}$. 
\end{rmk}

By a process called "normalizing near infinity", as described in \cite[§8.6]{[KM85]}, $\mathfrak{Y}(N)$ extends to a scheme $\mathfrak{X}(N)$ which is proper over $\mathbb{Z}[\zeta_N]$.

\begin{defn}\label{C(X(N)) def}
Let $\mathfrak{C}(\mathfrak{X}(N))$ denote the closed subscheme $\mathfrak{X}(N)-\mathfrak{Y}(N)$ of $\mathfrak{X}(N)$ endowed with the reduced scheme structure, called the \textbf{cuspidal locus} of $\mathfrak{X}(N)$. If the modular curve is clear from context, we denote $\mathfrak{C}(N)=\mathfrak{C}(\mathfrak{X}(N))$. 
\end{defn}

The following definition captures many of the desirable properties of $\mathfrak{X}(N)$. 

\begin{defn}\label{arith surf}
Let $R$ be a Dedekind domain. We call a regular, integral, projective, flat $R$-scheme $X$ of dimension 2 an \textbf{arithmetic surface}.
\end{defn}

\begin{thm}\label{X(N) is arith surf}
Let $N\geq3$. The scheme $\mathfrak{X}(N)$ is an arithmetic surface over $\mathbb{Z}[\zeta_{N}]$. Moreover, $\mathfrak{X}(N)$ is geometrically connected with reduced closed fibers and smooth over $\mathbb{Z}[\zeta_{N},1/N]$.
\end{thm}

\begin{proof}
We have already established that $\mathfrak{Y}(N)$ is regular and flat over $\mathbb{Z}[\zeta_{N}]$ of dimension 2. By \cite[{}10.9.1(2)]{[KM85]}, there exists an open neighborhood of $\mathfrak{C}(N)$ which is smooth over $\mathbb{Z}[\zeta_{N}]$. In particular, $\mathfrak{C}(N)$ is regular and flat over $\mathbb{Z}[\zeta_{N}]$ so $\mathfrak{X}(N)$ is regular and flat over $\mathbb{Z}[\zeta_{N}]$. Since $\mathfrak{X}(N)$
is connected  and regular, $\mathfrak{X}(N)$ is integral.

By construction, $\mathfrak{X}(N)$ is proper over $\mathbb{Z}[\zeta_{N}]$, so by \cite[Theorem 8.3.16]{[Liu02]}, we can conclude $\mathfrak{X}(N)$ is projective over $\mathbb{Z}[\zeta_{N}]$ hence is an arithmetic surface. By \cite[{}10.9.2(2)]{[KM85]}, $\mathfrak{X}(N)$ is geometrically connected and by \cite[{}13.8.4]{[KM85]}, $\mathfrak{X}(N)$ has reduced closed fibers. Smoothness follows from the smoothness of $\mathfrak{Y}(N)$ over $\mathbb{Z}[\zeta_{N},1/N]$.
\end{proof}

Let $N\geq3$ and let $R$ be a noetherian, regular, excellent $\mathbb{Z}[\zeta_{N}]$-algebra. We will construct an invertible sheaf on $\mathfrak{X}(N)_{/R}$ whose global sections will be defined as the space of modular forms. Let $f:{\cal E}\rightarrow\mathfrak{Y}(N)_{/R}$ denote the universal elliptic curve. We define $\underline{\omega}:=f_{*}\Omega_{{\cal E}/\mathfrak{Y}(N)_{/R}}^{1}$ to be the pushforward of the sheaf of Kahler differentials on ${\cal E}$. If $e:\mathfrak{Y}(N)_{/R}\rightarrow{\cal E}$ is the identity section of our elliptic curve, then we have $\underline{\omega}\simeq e^{*}\Omega_{{\cal E}/\mathfrak{Y}(N)_{/R}}^{1}$. Hence $\underline{\omega}$, being the pullback of an invertible sheaf, is an invertible sheaf on $\mathfrak{Y}(N)_{/R}$. According to \cite[{}10.13.2]{[KM85]}, there is a canonical way to extend $\underline{\omega}^{\otimes2}$ to an invertible sheaf on $\mathfrak{X}(N)_{/R}$, which we denote by $\underline{\omega}_{\mathfrak{X}(N)_{/R}}^{\otimes2}$ or simply $\underline{\omega}^{\otimes2}$ if the modular curve is clear from context.

\begin{prop}
The invertible sheaf $\underline{\omega}^{\otimes2}$ canonically extends to an invertible sheaf $\underline{\omega}_{\mathfrak{X}(N)_{/R}}^{\otimes2}$ on $\mathfrak{X}(N)_{/R}$.
\end{prop}

\begin{defn}\label{def:mod form}
We call the invertible sheaf $\underline{\omega}^{\otimes2k}$ on $\mathfrak{X}(N)_{/R}$ the \textbf{modular sheaf }(\textbf{of weight 2k}). The global sections $H^{0}(\mathfrak{X}(N)_{/R},\underline{\omega}^{\otimes2k})$ are known as \textbf{modular forms of weight $2k$ and level $\Gamma(N)$}. The global sections $H^{0}(\mathfrak{X}(N)_{/R},\underline{\omega}^{\otimes2k}(-\mathfrak{C}(N))$ are known as \textbf{cusp forms }of weight 2k and level $\Gamma(N)$.
\end{defn}

The formation of the modular sheaf $\underline{\omega}^{\otimes2}$ behaves well under base change. The following proposition is \cite[{}10.13.6]{[KM85]}.

\begin{prop}\label{Base change}
Let $R\rightarrow R'$ be an extension of noetherian, regular, excellent rings and let 
\[
j:\mathfrak{X}(N)_{/R'}\rightarrow\mathfrak{X}(N)_{/R}
\]
denote the induced base change map. Then we have an isomorphism of invertible sheaves
\[
j^{*}\underline{\omega}_{\mathfrak{X}(N)_{/R}}^{\otimes2}\simeq\underline{\omega}_{\mathfrak{X}(N)_{/R'}}^{\otimes2}.
\]
\end{prop}

Let $f:E\rightarrow S$ be an elliptic curve over a smooth $R$-scheme $S$. According to \cite[{}10.13.10]{[KM85]} (see also \cite[A1.4]{[Kat72]}), we have map $(f_{*}\Omega^{1}{}_{E/S})^{\otimes2}\rightarrow\Omega_{S/R}^{1}$ of ${\cal O}_{S}$-modules, known as the Kodaira-Spencer map, which becomes an isomorphism precisely when $E\rightarrow S$ represents a moduli problem which is \'{e}tale over ${\rm Ell}_{R}$. In particular, if $\mathfrak{Y}(N)_{/R}$ is smooth over $R$ (\emph{e.g. }if $N$ is a unit in $R$), we get an isomorphism $\underline{\omega}_{\mathfrak{Y}(N)}^{\otimes2}\simeq\Omega_{\mathfrak{Y}(N)/R}^{1}$.

\begin{thm}\label{K-S Iso}
Let $R$ be a noetherian, regular, excellent $\mathbb{Z}[\zeta_{N}]$-algebra containing $1/N$. The Kodaira-Spencer isomorphism  $\underline{\omega}_{\mathfrak{Y}(N)}^{\otimes2}\simeq\Omega_{\mathfrak{Y}(N)/R}^{1}$
on $\mathfrak{Y}(N)$ extends to an isomorphism on $\mathfrak{X}(N)$
\[
\underline{\omega}_{\mathfrak{X}(N)}^{\otimes2}\simeq\Omega_{\mathfrak{X}(N)/R}^{1}(\mathfrak{C}(N)).
\]
\end{thm}

\begin{proof}
By \cite[{}5.1.1]{[KM85]}, the moduli problem $\left[\Gamma(N)\right]^{{\rm can}}$ is finite \'{e}tale over ${\rm Ell}_{R}$. The result then follows from \cite[{}10.13.11]{[KM85]}.
\end{proof}

Suppose ${\cal F},{\cal F}'$ are two representable moduli problems that are finite over ${\rm Ell}_{R}$ and normal near infinity. Any morphism $\eta:{\cal F}\Rightarrow{\cal F}'$ of moduli problems is compatible with the usual morphism ${\cal F},{\cal F}'\Rightarrow\left[\Gamma(1)\right]$. Thus the induced map $\eta:\mathfrak{M}({\cal F})\rightarrow\mathfrak{M}({\cal F}')$ fits in the commutative diagram 
\[
\xymatrix{\mathfrak{M}({\cal F})\ar@{->}[rr]^{\eta}\ar@{->}[dr]_{j} &  & \mathfrak{M}({\cal F}')\ar@{->}[dl]^{j}\\
 & \mathbb{A}_{R}^{1}
}
\] Hence $\eta$ extends to a map $\bar{\eta}:\overline{\mathfrak{M}}({\cal F})\rightarrow\overline{\mathfrak{M}}({\cal F}')$
on the normalizations near infinity. 

\begin{prop}\label{Pullback of Modular Sheaf}
Let $R$ be an excellent, noetherian, regular ring and let ${\cal F}$ and ${\cal F}'$ be two representable moduli problems, both finite over ${\rm Ell}_{R}$ and both normal near infinity. Let $\eta:{\cal F}\Rightarrow{\cal F}'$ be a morphism of moduli problems over ${\rm Ell}_{R}$. Under the induced map $\bar{\eta}:\overline{\mathfrak{M}}({\cal F})\rightarrow\overline{\mathfrak{M}}({\cal F}')$ we have 
\[
\bar{\eta}^{*}(\underline{\omega}_{\overline{\mathfrak{M}}({\cal F}')}^{\otimes2})=\underline{\omega}_{\overline{\mathfrak{M}}({\cal F})}^{\otimes2}.
\]
\end{prop}

\begin{proof}
This is \cite[{}10.13.5(2)]{[KM85]}.
\end{proof}

\subsection{\label{sub:cusp=000026specialfiber}Cusps and the Special Fiber}

In this section we will investigate the cusps and special fiber of $\mathfrak{X}(N)$ and recall the theory of $q$-expansions for modular forms. First we establish that the formation of the cuspidal locus $\mathfrak{C}(N)$ and the formation of $\mathfrak{X}(N)$ behaves well under base change. The following is \cite[{}8.6.6, {}8.6.7]{[KM85]}.

\begin{prop}\label{Base Change Modular}
Let $R$ and $R'$ be excellent, noetherian, regular $\mathbb{Z}[\zeta_{N}]$-algebras. For any extension of scalars $R\rightarrow R'$, we have
\[
\mathfrak{Y}(N)_{/R}\times_{{\rm Spec}(R)}{\rm Spec}(R')\simeq\mathfrak{Y}(N)_{/R'}
\]
\[
\mathfrak{C}(\mathfrak{X}(N)_{/R})\times_{{\rm Spec}(R)}{\rm Spec}(R')\simeq\mathfrak{C}(\mathfrak{X}(N)_{/R'})
\]
\[
\mathfrak{X}(N)_{/R}\times_{{\rm Spec}(R)}{\rm Spec}(R')\simeq\mathfrak{X}(N)_{/R'}
\]
\end{prop}

For two groups $G_{1}$ and $G_{2}$, let ${\rm HomSurj}(G_{1},G_{2})$ denote the set of surjective group homomorphisms from $G_{1}$ to $G_{2}$. Any subgroup $\Gamma\leq{\rm GL}_{2}(\mathbb{Z}/N\mathbb{Z})$ acts on the set 
\[
{\rm HomSurj}((\mathbb{Z}/N\mathbb{Z})^{2},\mathbb{Z}/N\mathbb{Z})
\]
via $\gamma\cdot\Lambda=\Lambda\circ\gamma$ for any $\gamma\in\Gamma$. Let 
\[
{\rm HS}(N)={\rm HomSurj}((\mathbb{Z}/N\mathbb{Z})^{2},\mathbb{Z}/N\mathbb{Z})/\left\{ \pm I\right\} 
\]
where $I$ is the identity matrix in ${\rm GL}_{2}(\mathbb{Z}/N\mathbb{Z})$. 

According to \cite[{}10.9.1]{[KM85]}, the cusps $\mathfrak{C}(N)$
of $\mathfrak{X}(N)$ is a disjoint union of $\#{\rm HS}(N)$ many
sections of $\mathfrak{X}(N)$. Furthermore, the formal completion
of $\mathfrak{X}(N)$ along $\mathfrak{C}(N)$ is the disjoint union
of $\#{\rm HS}(N)$ copies of the formal spectrum ${\rm Spf}(\mathbb{Z}[\zeta_{N}]\left\llbracket q^{1/h_{\Lambda}}\right\rrbracket )$
for some integer $h_{\Lambda}\geq1$ dividing $N$, dependent on the
index $\Lambda\in{\rm HS}(N)$. Refer to \cite[{}10.2.5]{[KM85]} which
provides a canonical bijection between a cuspidal section and the
corresponding label in ${\rm HS}(N)$. 

\begin{prop}\label{Description of C(N) as Divisor}
Let $R$ be a Dedekind domain with fraction field $K$. As Weil divisors, we have 
\[
\mathfrak{C}(\mathfrak{X}(N)_{/R})=\sum_{x}\overline{\left\{ x\right\} }
\]
where the sum is indexed over $x\in\mathfrak{C}(\mathfrak{X}(N)_{/K})$.
\end{prop}

\begin{proof}
As a divisor, $\mathfrak{C}(\mathfrak{X}(N)_{/R})$ consists only of horizontal components since $\mathfrak{C}(\mathfrak{X}(N)_{/R})$ is a finite disjoint union of sections of $\mathfrak{X}(N)_{/R}$. By \cite[{}8.3.4]{[Liu02]}, these horizontal components are necessarily the closure, in $\mathfrak{X}(N)_{/R}$, of points closed in the generic fiber $\mathfrak{X}(N)_{/K}$. By Proposition \ref{Base Change Modular}, $\mathfrak{C}(\mathfrak{X}(N)_{/R})_{/K}\simeq\mathfrak{C}(\mathfrak{X}(N)_{/K})$, hence these points are precisely the cusps of $\mathfrak{X}(N)_{/K}$.
\end{proof}

Now we recall how the $q$-expansion map
\[
H^{0}(\mathfrak{X}(N)_{/R},\underline{\omega}^{\otimes2})\rightarrow R[[q^{1/N}]]
\]
is defined for a modular form over a $\mathbb{Z}[\zeta_{N}]$-algebra $R$ at a specified cusp.  Let $c$ be a cusp of $\mathfrak{X}(N)$, which corresponds to the image of a section ${\rm Spec}(\mathbb{Z}[\zeta_{N}])\rightarrow\mathfrak{X}(N)$ and is a connected component $\overline{\left\{ c'\right\} }$ of $\mathfrak{C}(N)$ where $c'$ is a cusp of the generic fiber. According to \cite[§VII, 2.3.2]{[DeRa]}, the cusp $c$ corresponds to a Tate curve ${\rm Tate}(q^{1/N})$ over $\mathbb{Z}[\zeta_{N}][[q^{1/N}]$ (see \cite[§VII, D\'{e}finition 1.16]{[DeRa]}) and a closed immersion $\mbox{\ensuremath{{\rm Spec}(\mathbb{Z}[\zeta_{N}][[q^{1/N}]])\rightarrow\mathfrak{X}(N)}}$ fitting in the diagram
\[
\xymatrix{ & {\rm Tate}(q^{1/N})\ar@{->}[d]\\
 & {\rm Spec}(\mathbb{Z}[\zeta_{N}][[q^{1/N}]]\ar@{->}[r] & \mathfrak{X}(N)
}
\]

By \cite[§VII, Corollaire 2.4]{[DeRa]} (see also \cite[{}10.9.1]{[KM85]}),
this closed immersion induces an isomorphism of formal schemes 
\[
{\rm Spf}(\mathbb{Z}[\zeta_{N}][[q^{1/N}]])\overset{\sim}{\longrightarrow}\hat{\mathfrak{X}}(N)_{c}
\]
where $\hat{\mathfrak{X}}(N)_{c}$ denotes the formal completion of $\mathfrak{X}(N)$ along the cusp $c$. According to \cite[{}8.8(T.2)]{[KM85]} (see also \cite[§VII, 1.16.1]{[DeRa]}) we have an isomorphism 
\[
\phi:\hat{\mathrm{Tate}}(q^{1/N})\overset{\sim}{\rightarrow}\hat{\mathbb{G}}_{m}
\]
between formal Lie groups. Let $\omega_{{\rm can},c}:=\phi^{*}(dX/X)$, where $dX/X$ is the standard invariant differential on $\mathbb{G}_{m}$. Pulling back by ${\rm Tate}(q^{1/N})\rightarrow\mathfrak{X}(N)$ allows us to identify $\hat{\omega}_{c}^{\otimes2}$ with $\mathbb{Z}[\zeta_{N}][[q^{1/N}]]\cdot\omega_{{\rm can},c}$. 
\begin{defn}\label{def: q-exp}
Let $f\in H^{0}(\mathfrak{X}(N),\underline{\omega}^{\otimes2})$ and let $c$ be a cusp of $\mathfrak{X}(N)$. Viewing $f$ in the completed stalk $\underline{\hat{\omega}}{}_{c}^{\otimes2}$, we write $f=f_{c}\cdot\omega_{{\rm can},c}$ for some $f_{c}\in\hat{{\cal O}}_{\mathfrak{X}(N),c}$. We call $f_{c}$ the \textbf{$q$-expansion of $f$ at $c$}. This gives a map 
\[
H^{0}(\mathfrak{X}(N),\underline{\omega}^{\otimes2})\rightarrow\mathbb{Z}[\zeta_{N}][[q^{1/N}].
\]
More generally, for a $\mathbb{Z}[\zeta_{N}]$-algebra $R$, we obtain the $q$-expansion map 
\[
H^{0}(\mathfrak{X}(N)_{/R},\underline{\omega}^{\otimes2})\rightarrow\mathbb{Z}[\zeta_{N}][[q^{1/N}]\otimes_{\mathbb{Z}[\zeta_{N}]}R\subset R[[q^{1/N}]]
\]
by tensoring with $R$. We will always use $\omega_{{\rm can},c}$ as a local generator for $\omega^{\otimes2}$ to obtain the $q$-expansion at $c$. 
\end{defn}

The following result, which is \cite[§VII, Th\'{e}or\'{e}me 3.9]{[DeRa]} and known as the $q$-expansion principle, essentially says $q$-expansions can detect the ``ring of definition''.

\begin{prop}\label{q-exp principle}
Let $A$ be a $\mathbb{Z}[\zeta_{N}]$-algebra, let $B$ be a subalgebra of $A$, and let $f\in H^{0}(\mathfrak{X}(N)_{/A},\underline{\omega}^{\otimes k})$. If the $q$-expansion of $f$ at every cusp of $\mathfrak{X}(N)_{/A}$ lies in $B[[q^{1/N}]]$, then $f$ is a modular form over $B$ i.e. $f\in H^{0}(\mathfrak{X}(N)_{/B},\underline{\omega}^{\otimes k})$.
\end{prop}

Let $N\geq3$ and $r\geq1$ be integers and let $p\geq2$ be a prime such that $p\nmid N$. We will describe the special fiber of $\mathfrak{X}(Np^{r})$. For simplicity, we consider the modular curve $\mathfrak{X}(Np^{r})$ over a $\mathbb{Z}[\zeta_{Np^{r}}]$-algebra $R$ which is a DVR of mixed characteristic $(0,p)$ in which $N$ is invertible and with fraction field $K$ and perfect residue field $k$. Let $\bar{\mathfrak{X}}(Np^{r})$ denote the special fiber of $\mathfrak{X}(Np^{r})$. 

We will recall the theory of Igusa curves, as developed in Section 12 of \cite{[KM85]}, to describe the irreducible components of $\bar{\mathfrak{X}}(Np^{r})$. For $i\in\mathbb{Z}$, let $\sigma^{i}:{\rm Spec}(k)\rightarrow{\rm Spec}(k)$ denote the map induced by the $i$th power of the Frobenius automorphism on $k$. For any $k$-scheme $S$, we let $S^{(\sigma^{i})}$ denote the pullback under $\sigma^{i}$ from which we also obtain a map $F_{{\rm abs}}^{i}:S\rightarrow S$ called the absolute Frobenius. More generally, for any scheme $X\rightarrow S$, we let $X^{(p^{i})}$ denote the pullback under $F_{{\rm abs}}^{i}$ from which we obtain a morphism $F_{X/S}:X\rightarrow X^{(p^{i})}$ of $S$-schemes, called the $i$th-fold relative Frobenius. 

Let ${\cal F}$ be a moduli problem on ${\rm Ell}_{k}$. We define the moduli problem ${\cal F}^{(\sigma^{i})}$ on ${\rm Ell}_{k}$ by extending scalars via $\sigma^{i}$ so that 
\[
{\cal F}(E\rightarrow S)={\cal F}^{(\sigma^{i})}(E^{(\sigma^{i})}\rightarrow S^{(\sigma^{i})}).
\]
If ${\cal F}$ is representable and finite over ${\rm Ell}_{k}$ and normal near infinity, then the same holds for ${\cal F}^{(\sigma^{i})}$ and we have 
\[
\overline{\mathfrak{M}}({\cal F}^{(\sigma^{i})})=\overline{\mathfrak{M}}({\cal F})^{(\sigma^{i})}.
\]

\begin{defn}
Let $E$ be an elliptic curve over an $k$-scheme $S$. An \textbf{Igusa structure of level $p^{r}$} on $E\rightarrow S$ is a point $P\in E^{(p^{r})}(S)$ which generates the kernel of Verschiebung $V^{r}:E^{(p^{r})}\rightarrow E$ in the sense of \cite[{}1.4.1]{[KM85]}. Let $[{\rm Ig}(p^{r})]:{\rm Ell}_{k}\rightarrow{\rm Sets}$ denote the moduli problem which assigns to each elliptic curve $E\rightarrow S$ the set of all Igusa structures of level $p^{r}$ on $E\rightarrow S$.
\end{defn}

According to \cite[{}12.7.1]{[KM85]}, if ${\cal F}$ is a representable moduli problem finite over ${\rm Ell}_{k}$ which is normal near infinity, the simultaneous moduli problem $[{\rm Ig}(p^{r})]\times{\cal F}$ is representable over ${\rm Ell}_{k}$ and normal near infinity. The following result is \cite[{}12.7.2]{[KM85]} applied to the simultaneous moduli problem $\left[{\rm Ig}(p^{r})\right]\times\left[\Gamma(N)\right]^{{\rm can},(\sigma^{-i})}$ over ${\rm Ell}_{k}$ for any $i\in\mathbb{Z}$. We denote \[{\rm Ig}(p^{r},N):=\overline{\mathfrak{M}}([{\rm Ig}(p^{r})],\left[\Gamma(N)\right]^{{\rm can},(\sigma^{-i})}).\]

\begin{prop}\label{Props of Ig}\leavevmode
\begin{enumerate}
    \item[a.] ${\rm Ig}(p^{r},N)$ is a proper smooth curve over $k$.
    \item[b.] The usual projection ${\rm Ig}(p^{r},N)\rightarrow\mathfrak{X}(N)_{/k}^{(\sigma^{-i})}$ is finite and \'{e}tale outside the supersingular points of $\mathfrak{X}(N)_{/k}^{(\sigma^{-i})}$ for all $i\in\mathbb{Z}$.
\end{enumerate}
\end{prop}

Next we define another moduli problem closely related to Igusa structures, which will directly appear in the description of $\bar{\mathfrak{X}}(Np^{r})$. 

\begin{defn}
Let $E$ be an elliptic curve over an $k$-scheme $S$ and fix $1\leq i\leq r$. An \textbf{exotic Igusa structure} \textbf{of level $(p^{r},i)$ }on $E\rightarrow S$ is a point $P\in E(S)$ such that $(O,P)$ is a Drinfeld $p^{i}$-basis of $E\rightarrow S$ along with a point $Q\in E^{(p^{r-i})}(S)$ such that $V^{r-i}(Q)=P$. Let 
\[
[{\rm ExIg}(p^{r},i)]:{\rm Ell}_{k}\rightarrow{\rm Sets}
\]
denote the moduli problem which assigns to each elliptic curve $E\rightarrow S$ the set of all exotic Igusa structures of level $(p^{r},i)$ on $E\rightarrow S$.
\end{defn}

The following result, which is \cite[{}12.10.6]{[KM85]} in the situation ${\cal F}=[\Gamma(N)]^{{\rm can},(\sigma^{-i})}$, relates Igusa structures with exotic Igusa structures and also establishes the representability of $[{\rm ExIg}(p^{r},i)]$. We denote 
\[
{\rm ExIg}(p^{r},i,N)=\overline{\mathfrak{M}}([{\rm ExIg}(p^{r},i)],\left[\Gamma(N)\right]).
\]

\begin{thm}\label{Iso Ig and ExIg}
For $1\leq i\leq r$, we have a canonical isomorphism ${\rm Ig}(p^{r},N)\simeq{\rm ExIg}(p^{r},i,N)$ sitting in the commutative diagram 
\[
\xymatrix{{\rm Ig}(p^{r},N)\ar@{->}[r]^{\simeq}\ar@{->}[d]^{\rho} & {\rm ExIg}(p^{r},i,N)\ar@{->}[d]^{\rho'}\\
\mathfrak{X}(N)_{/k}^{(\sigma^{-i})}\ar@{->}[r]^{F^{i}} & \mathfrak{X}(N)_{/k}
}
\]
where $F^{i}$ is the $i$th-fold relative Frobenius $F_{\mathfrak{X}(N)/k}^{i}$.
\end{thm}

We introduce the following definition from \cite[§13.1]{[KM85]} which will allow us to describe the special fiber of $\mathfrak{X}(Np^{r})$. Consider the general situation with $k$ a field, $Y$ a smooth scheme over $k$, and $X\rightarrow Y$ a finite flat morphism of schemes. Suppose there exists a nonempty finite set $S$ of $k$-rational points of $Y$ such that for each $y_{0}\in S$ there exists a unique closed $k$-rational point $x_{0}\in X$ over $y_{0}$ such that $\hat{{\cal O}}_{X,x_{0}}\simeq k[[x,y]]/(f)$ for some $f\in k[[x,y]]$. The points of $S$ are referred to as the \emph{supersingular points}; indeed in the situation $X$ and $Y$ are modular curves, $S$ will be taken to be the supersingular points of $Y$ which correspond to supersingular elliptic curves. Furthermore, we assume there is a finite collection of $k$-schemes $\left\{ Z_{i}\right\} _{i\in I}$ with a morphism $\coprod_{i\in I}Z_{i}\rightarrow X$ such that 

\begin{itemize}
    \item for each $i\in I$ and $y_{0}\in S$ there exists a unique closed, $k$-rational point $z_{i,0}\in Z_{i}$ over $y_{0}$.
    \item $Z_{i}$ is finite flat over $Y$ and $(Z_{i})^{{\rm red}}$ is smooth over $k$.
    \item $Z_{i}\rightarrow X$ is a closed immersion and $\coprod_{i\in I}Z_{i}\rightarrow X$ is an isomorphism over the complement of $S$ in $Y$.
\end{itemize} By \cite[{}13.1.3]{[KM85]}, known as the ``Crossings Theorem'', if $Y$ is connected, then the $Z_{i}$ are the irreducible components of $X$. Furthermore, if each $Z_{i}$ is reduced, then $X$ is also reduced. 

\begin{defn}
In the situation just discussed, we say $X$ is the \textbf{disjoint union of the $Z_{i}$'s with crossings at the supersingular points}.
\end{defn}

\begin{rmk}
Following \cite[{}13.1.7]{[KM85]}, we can relax the $k$-rationality of $x_{0}$ and $y_{0}$ and instead require $k$-rationality after extending scalars to a separable closure of $k$. In this situation, we still say $X$ is the disjoint union of the $Z_{i}$'s with crossings at the supersingular points.
\end{rmk}
\vspace{.25cm}
Next we describe $\bar{\mathfrak{X}}(Np^{r})$. We identify $(\mathbb{Z}/p^{r}\mathbb{Z})^{\times}$ as the subgroup 
\[
\left\{ \left(\begin{array}{cc}
u & 0\\
0 & u
\end{array}\right):u\in(\mathbb{Z}/p^{r}\mathbb{Z})^{\times}\right\} 
\]
of ${\rm GL}_{2}(\mathbb{Z}/p^{r}\mathbb{Z})$. According to \cite[{}13.7.1]{[KM85]}, the moduli problem $\left[\Gamma(p^{r})\right]^{{\rm can}}\otimes{k}$ on ${\rm Ell}_{k}$ assigns to each elliptic curve $E\rightarrow S$ the set of all Drinfeld $p^{r}$-bases $(P,Q)$ with $e_{p^{r}}(P,Q)=1$. Let $\phi:(\mathbb{Z}/p^{r}\mathbb{Z})^{2}\rightarrow E[p^{r}]$ denote the homomorphism of $S$-schemes corresponding to $(P,Q)$. Consider the diagram
\[
\xymatrix{ &  & (\mathbb{Z}/p^{r}\mathbb{Z})^{2}\ar@{->}[d]^{\phi}\ar@{-->}[dr]^{\Lambda}\\
0\ar@{->}[r] & \ker F^{r}\ar@{->}[r] & E[p^{r}]\ar@{->}[r]^{F^{r}} & \ker(V^{r})\ar@{->}[r] & 0
}
\]
where $\Lambda=F^{r}\circ\phi$. By \cite[{}13.7.2(3)]{[KM85]}, a choice of $\mathbb{Z}/p^{r}\mathbb{Z}$-basis of $(\mathbb{Z}/p^{r}\mathbb{Z})^{2}/\ker\Lambda$ defines an isomorphism $\mathbb{Z}/p^{r}\mathbb{Z}\overset{\sim}{\rightarrow}\ker(V^{r})$, allowing us to view $\Lambda$ as a surjective homomorphism $(\mathbb{Z}/p^{r}\mathbb{Z})^{2}\rightarrow\mathbb{Z}/p^{r}\mathbb{Z}$. 

\begin{defn}
The \textbf{component label} of $\phi$ is the class of $\Lambda$ in 
\[
(\mathbb{Z}/p^{r}\mathbb{Z})^{\times}/{\rm HomSurj}((\mathbb{Z}/p^{r}\mathbb{Z})^{2},\mathbb{Z}/p^{r}\mathbb{Z}).
\]
\end{defn}

By \cite[{}13.7.4, {}13.7.5]{[KM85]}, this establishes a canonical bijection between the irreducible components of the special fiber and the set of component labels. The following is \cite[{}13.7.6]{[KM85]}.

\begin{thm}\label{Special Fiber Comps}
The special fiber of $\mathfrak{X}(Np^{r})$ is the disjoint union, with crossings at the supersingular points of $\mathfrak{X}(N)_{/k}$, of the exotic Igusa curves ${\rm ExIg}(p^{r},r,N)$ over $\mathfrak{X}(N)_{/k}$ indexed by 
\[
(\mathbb{Z}/p^{r}\mathbb{Z})^{\times}/{\rm HomSurj}((\mathbb{Z}/p^{r}\mathbb{Z})^{2},\mathbb{Z}/p^{r}\mathbb{Z}).
\]
Furthermore, $\bar{\mathfrak{X}}(Np^{r})$ is reduced.
\end{thm}

Note that the claim $\bar{\mathfrak{X}}(Np^{r})$ is reduced comes from the fact that ${\rm ExIg}(p^{r},r,N)$ is reduced, together with the ``Crossings Theorem'' of \cite[{}13.1.3]{[KM85]} and \cite[{}13.1.4]{[KM85]}.

Recall the closed subscheme of cusps $\mathfrak{C}(Np^{r})$ of $\mathfrak{X}(Np^{r})$ is a disjoint union of $\#{\rm HS}(Np^{r})$ many sections of $\mathfrak{X}(Np^{r})$ which, by Proposition \ref{Description of C(N) as Divisor}, can be viewed as the closure $\overline{\left\{ x\right\} }$ of points $x\in\mathfrak{C}(\mathfrak{X}(N)_{/K})$. By Corollary \ref{Int of Ratl Pt}, $\overline{\left\{ x\right\} }$ intersects precisely one irreducible component of $\bar{\mathfrak{X}}(Np^{r})$. The following, which is \cite[{}13.9.3]{[KM85]}, tells us precisely the component label of the irreducible component which $\overline{\left\{ x\right\} }$ intersects with, given the index of $\overline{\left\{ x\right\} }$ in ${\rm HS}(Np^{r})$ (see the paragraph preceding Proposition \ref{Description of C(N) as Divisor}).

\begin{thm}\label{Cusp Irr Comp Labels}
The natural projection 
\[
\xymatrix{{\rm HS}(Np^{r})\ar@{=}[d]\\
\left\{ \pm I\right\} \backslash{\rm HomSurj}((\mathbb{Z}/N\mathbb{Z})^{2},\mathbb{Z}/N\mathbb{Z})\times{\rm HomSurj}((\mathbb{Z}/p^{r}\mathbb{Z})^{2},\mathbb{Z}/p^{r}\mathbb{Z})\ar@{->}[d]\\
(\mathbb{Z}/p^{r}\mathbb{Z})^{\times}\backslash{\rm HomSurj}((\mathbb{Z}/p^{r}\mathbb{Z})^{2},\mathbb{Z}/p^{r}\mathbb{Z})
}
\]
assigns to a component $\overline{\left\{ x\right\} }$ of $\mathfrak{C}(Np^{r})$ indexed by ${\rm HS}(Np^{r})$ the irreducible component of $\bar{\mathfrak{X}}(Np^{r})$ that $\overline{\left\{ x\right\} }$ intersects with.
\end{thm}

\begin{cor}\label{Irr Comps Same Num Int}
Each irreducible component of $\bar{\mathfrak{X}}(Np^{r})$ intersects with the same number of cuspidal components of $\mathfrak{C}(Np^{r})$.
\end{cor}

Every group homomorphism $\Lambda:(\mathbb{Z}/p^{r}\mathbb{Z})^{2}\rightarrow\mathbb{Z}/p^{r}\mathbb{Z}$ is uniquely determined by its image on the basis vectors $(1,0)$ and $(0,1)$. Let $\Lambda_{(a,b)}:(\mathbb{Z}/p^{r}\mathbb{Z})^{2}\rightarrow\mathbb{Z}/p^{r}\mathbb{Z}$ denote the map defined by 
\[
\Lambda_{(a,b)}(1,0)=a\mbox{ and }\Lambda_{(a,b)}(0,1)=b.
\]
A complete list of representatives in $(\mathbb{Z}/p^{r}\mathbb{Z})^{\times}/{\rm HomSurj}((\mathbb{Z}/p^{r}\mathbb{Z})^{2},\mathbb{Z}/p^{r}\mathbb{Z})$ is given by 
\[
\begin{cases}
\Lambda_{(1,-a)} & a\in\mathbb{Z}/p^{r}\mathbb{Z}\\
\Lambda_{(-pb,1)} & b\in\mathbb{Z}/p^{r-1}\mathbb{Z}
\end{cases}.
\]
With this labeling, the following is clear:

\begin{cor}\label{Num of Irr Comps}
The special fiber $\bar{\mathfrak{X}}(Np^{r})$ has $p^{r}+p^{r-1}$ many irreducible components.
\end{cor}

We will record the number of cuspidal components of $\mathfrak{C}(Np^{r})$ an irreducible component of $\bar{\mathfrak{X}}(Np^{r})$ intersects with. For convenience, we let $C(Np^{r})=\mathfrak{C}(\mathfrak{X}{}_{/\mathbb{Q}_{p}(\zeta_{Np^{r}}})$, the set of all cusps of the generic fiber. The following is \cite[{}4.2.10]{[Miy76]}.

\begin{lem}\label{Num of Cusps}
For $M\geq3$, we have 
\[
\#C(M)=\frac{1}{2}M^{2}\prod_{p\mid M}(1-1/p^{2})=\frac{1}{2M}\#{\rm SL}_{2}(\mathbb{Z}/MZ).
\]
\end{lem}

\noindent We will also need to know how to compute $\#{\rm SL}_{2}(\mathbb{Z}/M\mathbb{Z})$. The following is from \cite[{}4.2.3, {}4.2.4]{[Miy76]}.

\begin{lem}\label{size SL2(Z/MZ)}\leavevmode
\begin{enumerate}
    \item[a.] Let $A$ and $B$ be coprime integers. Then \[ {\rm SL}_{2}(\mathbb{Z}/AB\mathbb{Z})\simeq{\rm SL}_{2}(\mathbb{Z}/A\mathbb{Z})\times{\rm SL}_{2}(\mathbb{Z}/B\mathbb{Z}).\]
    \item[b.]  Let $p$ be prime and $r\geq1$ an integer. Then \[ \#{\rm SL}_{2}(\mathbb{Z}/p^{r}\mathbb{Z})=p^{3r}-p^{3r-2}. \]
\end{enumerate}
\end{lem}

\begin{prop}\label{num of irr comps int cusp}
Each irreducible component of $\bar{\mathfrak{X}}(Np^{r})$ intersects precisely 
\[
\frac{\#C(Np^{r})}{p^{r}+p^{r-1}}=\varphi(p^{r})\#C(N)
\]
many cuspidal components of $\mathfrak{C}(Np^{r})$. 
\end{prop}

\begin{proof}
By Corollary \ref{Irr Comps Same Num Int}, this quantity is independent of irreducible component $\Lambda$. Hence the number of cuspidal components $\Lambda$ intersects is equal to the total number of cusps of the generic fiber divided by the total number of irreducible components. By Corollary \ref{Num of Irr Comps}, the number of irreducible components is $p^{r}+p^{r-1}$. 

Using Lemma \ref{Num of Cusps} and Lemma \ref{size SL2(Z/MZ)}a, we compute:

\begin{alignat*}{1}
\frac{\#C(Np^{r})}{p^{r}+p^{r-1}} & =\frac{\#{\rm SL}_{2}(\mathbb{Z}/Np^{r}\mathbb{Z})}{2Np^{r}(p^{r}+p^{r-1})}\\
 & =\frac{\#{\rm SL}_{2}(\mathbb{Z}/p^{r}\mathbb{Z})}{p^{r}(p^{r}+p^{r-1})}\cdot\frac{\#{\rm SL}_{2}(\mathbb{Z}/N\mathbb{Z})}{2N}\\
 & =\frac{\#{\rm SL}_{2}(\mathbb{Z}/p^{r}\mathbb{Z})}{p^{r}(p^{r}+p^{r-1})}\#C(N)
\end{alignat*}
Lastly, using Lemma \ref{size SL2(Z/MZ)}b, we have 
\begin{alignat*}{1}
 & =\frac{p^{3r}-p^{3r-2}}{p^{r}(p^{r}+p^{r-1})}\#C(N)\\
 & =\varphi(p^{r})\#C(N).\tag*{\qedhere}
\end{alignat*}
\end{proof}

\section{\label{apx: circulant matrices}Appendix \textendash{} Circulant Matrices}

\begin{defn}
An $n\times n$ \textbf{circulant} matrix $C$ over a field is any matrix of the form
\[
C=\left(\begin{array}{ccccc}
c_{0} & c_{n-1} & \cdots & c_{2} & c_{1}\\
c_{1} & c_{0} & \cdots & c_{3} & c_{2}\\
\vdots & \vdots & \ddots & \vdots & \vdots\\
c_{n-2} & c_{n-3} & \cdots & c_{0} & c_{n-1}\\
c_{n-1} & c_{n-2} & \cdots & c_{1} & c_{0}
\end{array}\right)
\]
where each column is equal to the previous column shifted downward by 1, looping around as appropriate.
\end{defn}

The entries $c_{i,j}$ of an $n\times n$ circulant matrix can be
characterized by the equation 
\[
c_{i,j}=c_{i-j+1,1}
\]
for all $1\leq i,j\le n$ where the index $i-j+1$ is taken modulo $n$ among the residue classes in $\left\{ 0,1,\dots,n-1\right\} $. Let $\zeta_{n}=e^{2\pi i/n}$ denote a primitive $n$th root of unity. The following provides an explicit description of the eigenvalues and corresponding eigenvectors of a circulant matrix. 

\begin{lem}\label{Evals of Circulant}
The eigenvalues of a circulant matrix $C=(c_{ij})$ are precisely
\[
\lambda_{j}=\sum_{k=0}^{n-1}c_{k}\zeta_{n}^{(j-1)(n-k)}
\]
for $j=1,\dots,n$. A corresponding eigenvector of $\lambda_{j}$ is given by 
\[
\vec{v}_{j}=\frac{1}{\sqrt{n}}\left(1,\zeta_{n}^{j-1},\zeta_{n}^{2(j-1)},\dots,\zeta_{n}^{(n-1)(j-1)}\right).
\]
\end{lem}
\begin{proof}
We will verify that $C\vec{v}_{j}=\lambda_{j}\vec{v}_{j}$ for all $j$. For $1\leq i\leq n$, the $i^{{\rm th}}$ entry of $C\vec{v}_{j}$ is 
\[
\frac{1}{\sqrt{n}}\left(c_{i-1}+c_{i-2}\zeta_{n}^{j-1}+c_{i-3}\zeta_{n}^{2(j-1)}+\cdots+c_{0}\zeta^{(i-1)(j-1)}+c_{n-1}\zeta_{n}^{i(j-1)}\cdots+c_{i}\zeta_{n}^{(n-1)(j-1)}\right)
\]
\begin{align*}
 & = \frac{1}{\sqrt{n}}\sum_{k=0}^{i-1}c_{k}\zeta_{n}^{(i-1-k)(j-1)}+\frac{1}{\sqrt{n}}\sum_{\ell=i}^{n-1}c_{\ell}\zeta_{n}^{(n-1+i-\ell)(j-1)}\\
 & = \frac{1}{\sqrt{n}}\sum_{k=0}^{i-1}c_{k}\zeta_{n}^{(n+i-1-k)(j-1)}+\frac{1}{\sqrt{n}}\sum_{\ell=i}^{n-1}c_{\ell}\zeta_{n}^{(n-1+i-\ell)(j-1)}\\
 & = \frac{1}{\sqrt{n}}\sum_{k=0}^{n-1}c_{k}\zeta^{(n+i-1-k)(j-1)}
\end{align*}
while the $i^{{\rm th}}$ entry of $\lambda_{j}\vec{v}_{j}$ is
\begin{align*}
\frac{1}{\sqrt{n}}\zeta_{n}^{(i-1)(j-1)}\sum_{k=0}^{n-1}c_{k}\zeta_{n}^{(j-1)(n-k)} & = \frac{1}{\sqrt{n}}\sum_{k=0}^{n-1}c_{k}\zeta_{n}^{(j-1)(n-k)+(i-1)(j-1)}\\
 & = \frac{1}{\sqrt{n}}\sum_{k=0}^{n-1}c_{k}\zeta_{n}^{(j-1)(n+i-1-k)}
\end{align*}
which agrees with the $i^{{\rm th}}$ entry of $C\vec{v}_{j}$ for all $i$. 
\end{proof}

Circulant matrices are always diagonalizable (see \cite[§2]{[KaWh01]}) so we can write $C=PDP^{-1}$ where the columns of $P$ are the eigenvectors $\vec{v}_{1},\dots,\vec{v}_{n}$ and $D$ is a diagonal matrix whose diagonal entries are the corresponding eigenvalues $\lambda_{1},\dots,\lambda_{n}$. Explicitly, 
\[
P=\frac{1}{\sqrt{n}}\left(\begin{array}{cccccc}
1 & 1 & 1 & \cdots & 1 & 1\\
1 & \zeta_{n} & \zeta_{n}^{2} & \cdots & \zeta_{n}^{n-2} & \zeta_{n}^{n-1}\\
\vdots & \vdots & \vdots & \cdots & \vdots & \vdots\\
1 & \zeta_{n}^{n-2} & \zeta_{n}^{(n-2)2} & \cdots & \zeta_{n}^{(n-2)(n-2)} & \zeta_{n}^{(n-2)(n-1)}\\
1 & \zeta_{n}^{n-1} & \zeta_{n}^{(n-1)2} & \cdots & \zeta_{n}^{(n-1)(n-2)} & \zeta_{n}^{(n-1)(n-1)}
\end{array}\right).
\]
We see the $(i,j)$-entry of $P$ is equal to $\zeta_{n}^{(i-1)(j-1)}/\sqrt{n}$ for $1\leq i,j\leq n$. 

\begin{lem}
The matrix $P$ is unitary i.e. the inverse of $P$ is equal to the conjugate transpose $P^{*}$.
\end{lem}
\begin{proof}
Note that $P$ is symmetric and the conjugate of $\zeta_{n}$ is $\bar{\zeta}_{n}=\zeta_{n}^{-1}$. The $(i,j)$-entry of $PP^{*}$ is equal to 
\[
\frac{1}{n}\sum_{k=0}^{n-1}\zeta_{n}^{k(i-1)}\bar{\zeta}_{n}^{k(j-1)}=\frac{1}{n}\sum_{k=0}^{n-1}\zeta_{n}^{k(i-1)}\zeta_{n}^{-k(j-1)}
\]
\[
=\frac{1}{n}\sum_{k=0}^{n-1}\zeta_{n}^{k(i-j)}=\begin{cases}
1 & \mbox{if }i=j\\
0 & \mbox{if }i\ne j
\end{cases}.
\]
Hence $PP^{*}=I_{n}$ so $P^{*}=P^{-1}$.
\end{proof}

Now that we know $C=PDP^{*}$, we can describe the $(i,j)$ entry of $C^{-1}$, if it exists. 

\begin{prop}\label{inv of circ mat}
Let $C$ be an $n\times n$ invertible circulant matrix whose entries along the first column are $c_{0},\dots,c_{n-1}$ in that order. Then the $(i,j)$-entry of $C^{-1}$ is equal to 
\[
\frac{1}{n}\sum_{k=1}^{n}\lambda_{k-1}^{-1}\zeta_{n}^{(k-1)(i-j)}
\]
where $\lambda_{j}=\sum_{m=0}^{n-1}c_{m}\zeta_{n}^{(j-1)(n-m)}$ are the eigenvalues of $C$ for $1\leq j\leq n$.  
\end{prop}
\begin{proof}
We compute out the product $C^{-1}=PD^{-1}P^{*}$. The matrix $PD^{-1}$ will be of the form
\[
PD^{-1}=\frac{1}{\sqrt{n}}\left(\begin{array}{ccccc}
\lambda_{1}^{-1} & \lambda_{2}^{-1} & \lambda_{3}^{-1} & \cdots &  \lambda_{n}^{-1}\\
\lambda_{1}^{-1} & \lambda_{2}^{-1}\zeta_{n} & \lambda_{3}^{-1}\zeta_{n}^{2} & \cdots  & \lambda_{n}^{-1}\zeta_{n}^{n-1}\\
\vdots & \vdots & \vdots & \cdots &  \vdots\\
\lambda_{1}^{-1} & \lambda_{2}^{-1}\zeta_{n}^{n-2} & \lambda_{3}^{-1}\zeta_{n}^{(n-2)2} & \cdots & \lambda_{n}^{-1}\zeta_{n}^{(n-2)(n-1)}\\
\lambda_{1}^{-1} & \lambda_{2}^{-1}\zeta_{n}^{n-1} & \lambda_{3}^{-1}\zeta_{n}^{(n-1)2} & \cdots & \lambda_{n}^{-1}\zeta_{n}^{(n-1)(n-1)}
\end{array}\right).
\]
Thus the $(i,j)$-entry of $PD^{-1}$ is equal to 
\[
\lambda_{j}^{-1}\zeta_{n}^{(i-1)(j-1)}/\sqrt{n}
\]
and so the $(i,j)$-entry of $PD^{-1}P^{*}$ is equal to
\[
\frac{1}{n}\sum_{k=1}^{n}\lambda_{k}^{-1}\zeta_{n}^{(i-1)(k-1)}\zeta_{n}^{-(k-1)(j-1)}=\frac{1}{n}\sum_{k=1}^{n}\lambda_{k}^{-1}\zeta_{n}^{(k-1)(i-j)}.\qedhere
\] 
\end{proof}

Using our explicit description of the entries of $C^{-1}$, we can establish the following.

\begin{cor}\label{Inv of Circ is Circ}
The inverse of an invertible circulant matrix is circulant.
\end{cor}
\begin{proof}
Let $C=(c_{i,j})$ be an $n\times n$ circulant matrix whose entries
along the first column are ordered $c_{0},c_{1},\dots,c_{n-1}$ so
that $c_{m+1,1}=c_{m}$ for $0\leq m\leq n-1$. By Lemma \ref{inv of circ mat}, the $(i,j)$-entry of $C^{-1}$ is equal to 
\[
c^{i,j}=\frac{1}{n}\sum_{k=1}^{n}\left(\sum_{m=0}^{n-1}c_{m}\zeta_{n}^{(k-1)(n-m)}\right)\zeta_{n}^{(k-1)(i-j)}.
\]
To show $C^{-1}$ is circulant, we will establish the relationship $c^{i,j}=c^{i-j+1,1}$ where the index $i-j+1$ is taken modulo $n$ in the residue class $\{0,\dots,n-1\}$. In the above expression for $c^{i,j}$, the index $(i,j)$ only appears in the term $\zeta_{n}^{(k-1)(i-j)}$. Observe that 
\[
\zeta_{n}^{(k-1)((i-j+1)-1)}=\zeta_{n}^{(k-1)(i-j)}.
\]
Hence $c^{i,j}=c^{i-j+1,1}$ so $C^{-1}$ is circulant. 
\end{proof}

\section{\label{apx: Inv via WMI}Appendix \textendash{} Inverse via the Woodbury Matrix Identity}

Let $A$ be an $n\times n$ invertible matrix, $C$ an invertible
$k\times k$ matrix where $k\leq n$, $U$ an $n\times k$ matrix,
and $V$ a $k\times n$ matrix. The Woodbury matrix identity states
\[
(A+UCV)^{-1}=A^{-1}-A^{-1}U(C^{-1}+VA^{-1}U)^{-1}VA^{-1}.
\]
A proof of this identity can be found in \cite[§1.3]{[HeSe81]}. This identity allows us to compute the inverse of $A+UCV$ provided we can easily compute the inverses of $A$ and $C^{-1}+VA^{-1}U$. 

In general, we can use the Woodbury Matrix Identity to compute the inverse of $A+N$ where $A$ is an $n\times n$ invertible matrix and $N$ is an $n\times n$ matrix of rank $k$, provided $A+N$ is invertible. Let $U$ be the $n\times k$ matrix whose columns $\vec{v}_{1},\dots,\vec{v}_{k}$ are the $k$ linearly independent columns of $N$. Let $\vec{u}_{i}$ denote the vector in the $i$th column of $N$ which can be expressed as 
\[
\vec{u}_{i}=c_{1}\vec{v}_{1}+\cdots+c_{k}\vec{v}_{k}.
\]
We define the $i$th column of the $k\times n$ matrix $V$ to consist of entries $c_{1},\dots,c_{k}$ in that order. Consequently, $N=UCV$ where $C=I_{k}$ is the $k\times k$ identity matrix.

Using the Woodbury Matrix Identity, we will provide a formula for the inverse in the following situation as encountered in Section \ref{sub:Inverting T}. Let $A$ be a block diagonal matrix with two invertible blocks of sizes $n\times n$ and $m\times m$. The inverse $A^{-1}$ is also block diagonal with blocks of the same size. We write 
\[
A^{-1}=\left(\begin{array}{cc}
\begin{array}{ccc}
a_{11} & \cdots & a_{1n}\\
\vdots & \ddots & \vdots\\
a_{n1} & \cdots & a_{nn}
\end{array} & \mathbf{0}\\
\mathbf{0} & \begin{array}{ccc}
b_{11} & \cdots & b_{1m}\\
\vdots & \ddots & \vdots\\
b_{m1} & \cdots & b_{mm}
\end{array}
\end{array}\right).
\]
We also consider the case
\[
U=\left(\begin{array}{cc}
0_{1} & 1_{1}\\
0_{2} & 1_{2}\\
\vdots & \vdots\\
0_{n} & 1_{n}\\
1_{1} & 0_{1}\\
\vdots & \vdots\\
1_{m} & 0_{m}
\end{array}\right)\mbox{ and }V=\left(\begin{array}{ccccccc}
1_{1} & 1_{2} & \cdots & 1_{n} & 0_{1} & \cdots & 0_{m}\\
0_{1} & 0_{2} & \cdots & 0_{n} & 1_{1} & \cdots & 1_{m}
\end{array}\right)
\]
where the subscripts on the $0$ and $1$ entries are there to help keep track of their position. 

\begin{prop}\label{Gen Formula for T^-1}
Let $T=A+UV$ where $A,U,$ and $V$ are the given matrices above. Let $c^{i,j}$ denote the $(i,j)$-entry of $T^{-1}$. We have 
\[
c^{ij}=\begin{cases}
{\displaystyle a_{i,j}+\frac{\beta}{1-\alpha\beta}\left(\sum_{k=1}^{n}a_{i,k}\right)\left(\sum_{k=1}^{n}a_{k,j}\right)} & \mbox{if }1\leq i,j\leq n\\
{\displaystyle \frac{-1}{1-\alpha\beta}\left(\sum_{k=1}^{n}a_{i,k}\right)\left(\sum_{k=1}^{m}b_{k,j}\right)} & \mbox{if }1\leq i\leq n\mbox{ and }n<j\leq n+m\\
{\displaystyle \frac{-1}{1-\alpha\beta}\left(\sum_{k=1}^{m}b_{i,k}\right)\left(\sum_{k=1}^{n}a_{k,j}\right)} & \mbox{if }1\leq j\leq n\mbox{ and }n<i\leq n+m\\
{\displaystyle b_{i,j}+\frac{\alpha}{1-\alpha\beta}\left(\sum_{k=1}^{m}b_{i,k}\right)\left(\sum_{k=1}^{m}b_{k,j}\right)} & \mbox{if }n<i,j\leq n+m
\end{cases}
\]
where $\alpha=\sum a_{i,j}$ is the sum of all entries in the first block in $A^{-1}$ and $\beta=\sum b_{i,j}$ is the sum of all entries in the second block in $A^{-1}$.
\end{prop}
\begin{proof}
By the Woodbury Matrix Identity,
\[
T^{-1}=A^{-1}-A^{-1}U(I_{2}+VA^{-1}U)^{-1}VA^{-1}.
\]
We first compute{\footnotesize{}
\begin{align*}
VA^{-1}U & = V\left(\begin{array}{cc}
\begin{array}{ccc}
a_{11} & \cdots & a_{1n}\\
\vdots & \ddots & \vdots\\
a_{n1} & \cdots & a_{nn}
\end{array} & \mathbf{0}\\
\mathbf{0} & \begin{array}{ccc}
b_{11} & \cdots & b_{1m}\\
\vdots & \ddots & \vdots\\
b_{m1} & \cdots & b_{mm}
\end{array}
\end{array}\right)\left(\begin{array}{cc}
0_{1} & 1_{1}\\
\vdots & \vdots\\
0_{n} & 1_{n}\\
1_{1} & 0_{1}\\
\vdots & \vdots\\
1_{m} & 0_{m}
\end{array}\right)\\
 & = \left(\begin{array}{cccccc}
1_{1} & \cdots & 1_{n} & 0_{1} & \cdots & 0_{m}\\
0_{1} & \cdots & 0_{n} & 1_{1} & \cdots & 1_{m}
\end{array}\right)\left(\begin{array}{cc}
0 & {\displaystyle {\displaystyle \sum_{i=1}^{n}a_{1i}}}\\
\vdots & \vdots\\
0 & {\displaystyle {\displaystyle \sum_{i=1}^{n}a_{n,i}}}\\
{\displaystyle \sum_{i=1}^{m}b_{1,i}} & {\displaystyle 0}\\
\vdots & \vdots\\
{\displaystyle {\displaystyle \sum_{i=1}^{m}b_{m,i}}} & 0
\end{array}\right)\\
 & = \left(\begin{array}{cc}
0 & {\displaystyle {\displaystyle \sum_{i,j=1}^{n}a_{ij}}}\\
{\displaystyle {\displaystyle \sum_{i,j=1}^{m}b_{ij}}} & 0
\end{array}\right)\\
 & = \left(\begin{array}{cc}
0 & \alpha\\
\beta & 0
\end{array}\right).
\end{align*}
} Hence 
\[
(I_{2}+VA^{-1}U)^{-1}=\left(\begin{array}{cc}
1 & \alpha\\
\beta & 1
\end{array}\right)^{-1}=\frac{1}{1-\alpha\beta}\left(\begin{array}{cc}
1 & -\alpha\\
-\beta & 1
\end{array}\right).
\]

Next we will compute
\begin{align*}
U(I_{2}^{-1}+VA^{-1}U)^{-1}V & = \frac{1}{1-\alpha\beta}U\left(\begin{array}{cc}
1 & -\alpha\\
-\beta & 1
\end{array}\right)\left(\begin{array}{cccccc}
1_{1} & \cdots & 1_{n} & 0_{1} & \cdots & 0_{m}\\
0_{1} & \cdots & 0_{n} & 1_{1} & \cdots & 1_{m}
\end{array}\right)\\
 & = \frac{1}{1-\alpha\beta}\left(\begin{array}{cc}
0_{1} & 1_{1}\\
\vdots & \vdots\\
0_{n} & 1_{n}\\
1_{1} & 0_{1}\\
\vdots & \vdots\\
1_{m} & 0_{m}
\end{array}\right)\left(\begin{array}{cccccc}
1_{1} & \cdots & 1_{n} & -\alpha & \cdots & -\alpha\\
-\beta & \cdots & -\beta & 1_{1} & \cdots & 1_{m}
\end{array}\right)\\
 & = \frac{1}{1-\alpha\beta}\left(\begin{array}{cccccc}
-\beta & \cdots & -\beta & 1 & \cdots & 1\\
\vdots & \ddots & \vdots & \vdots & \ddots & \vdots\\
-\beta & \cdots & -\beta & -1 & \cdots & -1\\
1 & \cdots & 1 & -\alpha & \cdots & -\alpha\\
\vdots & \ddots & \vdots & \vdots & \ddots & \vdots\\
1 & \cdots & 1 & -\alpha & \cdots & -\alpha
\end{array}\right).
\end{align*}
Lastly, we compute
\begin{align*}
 &  A^{-1}U(I_{2}^{-1}+VA^{-1}U)^{-1}VA^{-1}\\
 & = \frac{1}{1-\alpha\beta}\left(\begin{array}{cccccc}
a_{11} & \cdots & a_{1,n} & 0 & \cdots & 0\\
\vdots & \ddots & \vdots & \vdots & \ddots & \vdots\\
a_{n,1} & \cdots & a_{n,n} & 0 & \cdots & 0\\
0 & \cdots & 0 & b_{11} & \cdots & b_{1,m}\\
\vdots & \ddots & \vdots & \vdots & \ddots & \vdots\\
0 & \cdots & 0 & b_{m,1} & \cdots & b_{m,m}
\end{array}\right)\left(\begin{array}{cccccc}
-\beta & \cdots & -\beta & 1 & \cdots & 1\\
\vdots & \ddots & \vdots & \vdots & \ddots & \vdots\\
-\beta & \cdots & -\beta & 1 & \cdots & 1\\
1 & \cdots & 1 & -\alpha & \cdots & -\alpha\\
\vdots & \ddots & \vdots & \vdots & \ddots & \vdots\\
1 & \cdots & 1 & -\alpha & \cdots & -\alpha
\end{array}\right)A^{-1}\\
 & = \frac{1}{1-\alpha\beta}\left(\begin{array}{cccccc}
-\beta\sum a_{1,\ell} & \cdots & -\beta\sum a_{1,\ell} & \sum a_{1,\ell} & \cdots & \sum a_{1,\ell}\\
\vdots & \ddots & \vdots & \vdots & \ddots & \vdots\\
-\beta\sum a_{n,\ell} & \cdots & -\beta\sum a_{n,\ell} & \sum a_{n,\ell} & \cdots & \sum a_{n,\ell}\\
\sum b_{1,\ell} & \cdots & \sum b_{1,\ell} & -\alpha\sum b_{1,\ell} & \cdots & -\alpha\sum b_{1,\ell}\\
\vdots & \ddots & \vdots & \vdots & \ddots & \vdots\\
\sum b_{m,\ell} & \cdots & \sum b_{m,\ell} & -\alpha\sum b_{m,\ell} & \cdots & -\alpha\sum b_{m,\ell}
\end{array}\right)A^{-1}.
\end{align*}

We now describe the $(i,j)$-entry of the above product in four different cases.

\textbf{Case 1}($1\leq i,j\leq n$): The $(i,j)$-entry is equal to 
\[
\frac{1}{1-\alpha\beta}\sum_{k=1}^{n}\left(-\beta\sum_{m=1}^{n}a_{i,m}\right)a_{k,j}=-\frac{\beta}{1-\alpha\beta}\left(\sum_{k=1}^{n}a_{i,k}\right)\left(\sum_{k=1}^{n}a_{k,j}\right).
\]

\textbf{Case 2} ($1\leq i\leq n$ and $n<j\leq n+m$): The $(i,j)$-entry
is equal to 
\[
\frac{1}{1-\alpha\beta}\sum_{k=1}^{m}\left(\sum_{m=1}^{n}a_{i,m}\right)b_{k,j}=\frac{1}{1-\alpha\beta}\left(\sum_{k=1}^{n}a_{i,k}\right)\left(\sum_{k=1}^{m}b_{k,j}\right).
\]

\textbf{Case 3} ($1\leq j\leq n$ and $n<i\leq n+m$): The $(i,j)$-entry
is equal to
\[
\frac{1}{1-\alpha\beta}\sum_{k=1}^{n}\left(\sum_{m=1}^{m}b_{i,m}\right)a_{k,j}=\frac{1}{1-\alpha\beta}\left(\sum_{k=1}^{m}b_{i,k}\right)\left(\sum_{k=1}^{n}a_{k,j}\right).
\]

\textbf{Case 4} ($n<i,j\leq n+m$): The $(i,j)$-entry is equal to
\[
\frac{1}{1-\alpha\beta}\sum_{k=1}^{p^{r-1}}\left(-\alpha\sum_{m=1}^{m}b_{i,m}\right)b_{k,j}=-\frac{\alpha}{1-\alpha\beta}\left(\sum_{k=1}^{m}b_{i,k}\right)\left(\sum_{k=1}^{m}b_{k,j}\right).
\]
Subtracting these quantities from corresponding $(i,j)$-entry of
$A^{-1}$ gives our desired result.
\end{proof}

\section*{Acknowledgements}
\addcontentsline{toc}{section}{Acknowledgement}
Most of the results in this paper originated from my PhD thesis.  I would like to thank Bryden Cais for his guidance, support, and incredibly helpful feedback.  I would also like to thank Brandon Levin, Doug Ulmer, and Hang Xue for their helpful comments. This research did not receive any specific grant from funding agencies in the public, commercial, or not-for-profit sectors.

\addcontentsline{toc}{section}{References}

\printbibliography

\end{document}